\documentclass[aps,prb]{revtex4-2}

\usepackage{amsmath,amssymb,amsthm}
\usepackage{mathtools}
\usepackage{xspace}
\usepackage{braket}
\usepackage{hyphenat}
\usepackage{bbold}
\usepackage{enumerate}
\usepackage{tabularx}

\usepackage[british]{babel}

\usepackage{srcltx}
\usepackage{color}

\usepackage{graphicx}
\usepackage{amsmath}
\usepackage{amssymb}
\usepackage{bm}

\begin{document}

\title{Evaluation of eight different families of cubic Euler sums}
\author{J. Braun}
\affiliation{Ludwig Maximilians-Universit{\"a}t, M{\"u}nchen, Germany}
\author{H. J. Bentz}
\affiliation{Institut f\"{u}r Mathematik und Informatik, Universit\"{a}t Hildesheim, Germany}
\date{\today}

\begin{abstract}
We present a study on cubic Euler sums of degree four, five and six, where three different types of denominators $1/k^n$, $1/((2k-1)^n)$ and $1/(k(2k-1))$ will be considered
We demonstrate that for all three orders the complete variety of corresponding nonlinear Euler sums belonging to the eight different families can be explicitly calculated in
terms of zeta values and polylogarithmic values $Li_4(1/2)$, $Li_5(1/2)$, $Li_6(1/2)$, $Li_6(-1/2)$ and $Li_6(-1/8)$.

\end{abstract}
\maketitle

\section{Introduction}
Linear Euler sums like
\begin{eqnarray}
t(n,m) = \sum^{\infty}_{k=1} \frac{H_k^{(n)}}{k^m} 
\end{eqnarray}
or
\begin{eqnarray}
s(n,m) = \sum^{\infty}_{k=1} \frac{h_k^{(n)}}{k^m} 
\end{eqnarray}
had been discovered in the past by many authors \cite{fla97,sit85,zeh07,ade16,ade16a,ce20,bra20,bra21}. An excellent overview about these works is found in the paper of Flajolet and Salvy
\cite{fla97} and in citations therein. Further work on corresponding nonlinear Euler sums can be found, for example, in \cite{ce16,ce17,ce18,bra22,bra24}. 
The aim of the present paper is to present a systematic study on cubic Euler sums of different degrees, where three different types of denominators $1/k^n$, $1/((2k-1)^n)$ and
$1/(k(2k-1))$ will be considered in accordance to our former publications on quadratic Euler sums \cite{bra22,bra24}. Our calculational procedure allows to calculate all
cubic Euler sums belonging to the corresponding eight families in terms zeta values and corresponding polylogarithmic values. In detail we introduce a generalized calculational scheme
that is based on proper two-valued linear and nonlinear integer functions \cite{bra20,bra21,bra22,bra24}. Furthermore, we use extensively the partial fraction decomposition technique for
various Euler sums in combination with multiple zeta value theory (MZV) \cite{ce20,au20}.

\section{The order four case for cubic Euler sums}
In the first section we explicitly calculate all cubic Euler sums of degree four. In this case only cubic members of the third, fifth and sixth family exist.

\subsection{Third family}
We start with the following sum:
\begin{eqnarray}
\sum^{\infty}_{k=1} \frac{H_k^3}{k(2k-1)} &=& -\frac{37}{8}\zeta(4) -7ln(2)\zeta(3) + 7\zeta(3) + 6\left( ln(2) \right)^2\zeta(2) - 12ln(2)\zeta(2) + 8\zeta(2)
+ 16ln(2)- 24\left( ln(2) \right)^2 \nonumber \\ &+& 16\left( ln(2) \right)^3-4\left( ln(2) \right)^4 + \sum^{\infty}_{k=1} \frac{h_k}{k^3}~.
\end{eqnarray}
The corresponding calculational procedure starts with the following expression:
\begin{eqnarray}
\sum^{\infty}_{i=1} \frac{H_i}{(2i-1)} \left( \sum^{\infty}_{k=1} \frac{H_k}{k(i+k)} \right) = \sum^{\infty}_{k=1} \frac{H_k}{k} \left( \sum^{\infty}_{i=1} \frac{H_i}{(i+k)(2i-1)} \right).
\end{eqnarray}
Using the identity \cite{ce16}:
\begin{eqnarray}
\sum^{\infty}_{k=1} \frac{H_k}{k(i+k)} = \frac{1}{2i}H^{(2)}_i + \frac{1}{2i}H^2_i - \frac{H_i}{i^2} + \zeta(2)\frac{1}{i}
\end{eqnarray}
and the identity introduced as lemma 2a in \cite{bra24} 
\begin{eqnarray}
\sum^{\infty}_{k=1} \frac{H_k}{(2k-1)(i+k)} = \frac{1}{2i+1} \left( \zeta(2) - 2 \left(ln(2)\right)^2 + \zeta(2) -\frac{H_i}{i} +\frac{1}{2} H^{(2)}_i +
\frac{1}{2} H^2_i \right)
\end{eqnarray}
the explicit computation results in Eq.~(3).

\subsection{Fifth family}
The only member of the fifth family results to:
\begin{eqnarray}
\sum^{\infty}_{k=1} \frac{h_k^3}{k(2k-1)} &=& \frac{57}{64}\zeta(4) + \frac{7}{8}ln(2)\zeta(3) + \frac{9}{4}\left(ln(2)\right)^2\zeta(2) 
- \frac{1}{8}\sum^{\infty}_{k=1}\frac{h_k}{k^3}~.
\end{eqnarray}
Here the calculational procedure is very similar to that introduced before. We start with the following expression:
\begin{eqnarray}
\sum^{\infty}_{i=1} \frac{h_i}{2i-1} \left( \sum^{\infty}_{k=1} \frac{H_k}{(2k-1)(2i+2k-1)} \right) = \sum^{\infty}_{k=1} \frac{H_k}{2k-1} \left( \sum^{\infty}_{i=1}
\frac{h_i}{(2i-1)(2i+2k-1)} \right).
\end{eqnarray}
The two-valued help function on the left side is defined  by lemma 4a from \cite{bra24} and the two-valued help function on the right side is defined by Eq.~{201} in \cite{bra24}.
With this the explicit computation results in Eq.~(7).

\subsection{Sixth family}
Finally, in the cubic case we have two members of the sixth family. It follows first:
\begin{eqnarray}
\sum^{\infty}_{k=1} \frac{H_k h_k^2}{k(2k-1)} &=& \frac{1}{2}\zeta(4) + \frac{11}{4}ln(2)\zeta(3) + \frac{3}{4}\zeta(3) - 3\left(ln(2)\right)^2\zeta(2) 
+ 3ln(2)\zeta(2) - \frac{1}{2}\sum^{\infty}_{k=1}\frac{h_k}{k^3}~.
\end{eqnarray}
The corresponding calculational procedure is based on the two-valued help function defined by Eq.~(242) from \cite{bra24}:
\begin{eqnarray}
\sum^{\infty}_{i=1} \frac{h_i}{2i-1}  \left( \sum^{\infty}_{k=1} \frac{h_k}{k(i+k)} \right) = 2ln(2)\sum^{\infty}_{i=1} \frac{h^2_i}{i(2i-1)} + 
\sum^{\infty}_{i=1} \frac{H_i h^2_i}{i(2i-1)} - \sum^{\infty}_{i=1} \frac{h_i}{2i-1}  \left( \sum^{i}_{k=1} \frac{h_k}{k}  \right)~.
\end{eqnarray}
The last sum on the right side can be calculated as follows. We define: 
\begin{eqnarray}
\sum^{\infty}_{i=1} \frac{h_i}{2i-1} \left( \sum^{i}_{k=1} \frac{h_k}{k} \right) = -\sum^{\infty}_{i=1} \frac{h_i}{i(2i-1)} \left( \sum^{i}_{k=1} \frac{h_k}{k(2k-1)} \right) +
2\sum^{\infty}_{i=1} \frac{h_i}{i(2i-1)} \left( \sum^{i}_{k=1} \frac{h_k}{2k-1} \right)~.  
\end{eqnarray} 
With
\begin{eqnarray}
\sum^{\infty}_{i=1} \frac{h_i}{i(2i-1)} \left( \sum^{i}_{k=1} \frac{h_k}{k(2k-1)} \right) = \frac{1}{2}\left( \sum^{\infty}_{i=1} \frac{h_i}{2i-1} \right)^2 +
\frac{1}{2} \left( \sum^{\infty}_{i=1} \frac{h_i^2}{i^2(2i-1)^2} \right)
\end{eqnarray}
it follows further: 
\begin{eqnarray}
\sum^{\infty}_{i=1} \frac{h_i}{2i-1} \left( \sum^{i}_{k=1} \frac{h_k}{k} \right) = \sum^{\infty}_{i=1} \frac{h_i^3}{i(2i-1)} + \sum^{\infty}_{i=1} \frac{h_i h_i^{(2)}}{i(2i-1)} 
- \frac{1}{2}\left( \sum^{\infty}_{i=1} \frac{h_i}{2i-1} \right)^2 - \frac{1}{2} \left( \sum^{\infty}_{i=1} \frac{h_i^2}{i^2(2i-1)^2} \right)~.
\end{eqnarray}
Using Eq.~(7) the right side of Eq.~(10) can be calculated explicitly. For the left side of Eq.~(10) we can write
\begin{eqnarray}
\sum^{\infty}_{i=1} \frac{h_i}{2i-1} \left( \sum^{\infty}_{k=1} \frac{h_k}{k(i+k)} \right) = \sum^{\infty}_{k=1} \frac{h_k}{k} \left( \sum^{\infty}_{i=1} \frac{h_i}{(2i-1)(i+k)} \right)~.
\end{eqnarray}
Defining
\begin{eqnarray}
\sum^{\infty}_{i=1} \frac{h_i}{(2i-1)(i+k)} = 2ln(2)\frac{h_k}{2k+1} + \zeta(2)\frac{1}{2k+1} + \frac{H_k h_k}{2k+1} - \frac{1}{2k+1}\sum^{k}_{i=1} \frac{h_i}{i}~,
\end{eqnarray}
where this identity results from the simple solution of a corresponding first order difference equation. This way we get:
\begin{eqnarray}
\sum^{\infty}_{k=1} \frac{h_k}{k} \left( \sum^{\infty}_{i=1} \frac{h_i}{(2i-1)(i+k)} \right) &=& 2ln(2)\sum^{\infty}_{k=1}\frac{h_k^2}{k(2k+1)} + 
\zeta(2)\sum^{\infty}_{k=1} \frac{h_k}{k(2k+1)} + \sum^{\infty}_{k=1} \frac{H_k h_k^2}{k(2k+1)} \nonumber \\ &-& 
\sum^{\infty}_{k=1} \frac{h_k}{k(2k+1)} \left( \sum^{k}_{i=1} \frac{h_i}{i} \right)~.
\end{eqnarray}
Inserting Eq.~(13) and Eq.~(16) into  Eq.~(10) and using the proportionality  
\begin{eqnarray}
-\sum^{\infty}_{k=1} \frac{H_k h_k^2}{k(2k+1)} \sim \sum^{\infty}_{k=1} \frac{H_k h_k^2}{k(2k-1)}
\end{eqnarray}
Eq.~(9) finally results.

The second cubic Euler sum of order four belonging to the sixth family is given by the following expression:
\begin{eqnarray}
\sum^{\infty}_{k=1} \frac{H_k^2 h_k}{k(2k-1)} &=& \frac{1}{4}\zeta(4) - \frac{17}{2}ln(2)\zeta(3) + \frac{17}{2}\zeta(3) + 4\left(ln(2)\right)^2\zeta(2) 
- 8ln(2)\zeta(2) + 4\zeta(2) + \sum^{\infty}_{k=1}\frac{h_k}{k^3}~.
\end{eqnarray}
Again the calculational procedure is similar to that introduced for Eq.~(9). We start with
\begin{eqnarray}
\sum^{\infty}_{i=1} \frac{h_i}{2i-1}  \left( \sum^{\infty}_{k=1} \frac{H_k}{k(i+k)} \right) = \zeta(2) \sum^{\infty}_{i=1} \frac{h_i}{i(2i-1)} - 
\sum^{\infty}_{i=1} \frac{H_i h_i}{i^2(2i-1)} + \frac{1}{2}\sum^{\infty}_{i=1} \frac{H^{(2)}_i h_i}{i(2i-1)} + \frac{1}{2}\sum^{\infty}_{i=1} \frac{H^2_i h_i}{i(2i-1)} ~.
\end{eqnarray}
Rearranging the summations on the left side of Eq.~(19):
\begin{eqnarray}
\sum^{\infty}_{i=1} \frac{h_i}{2i-1} \left( \sum^{\infty}_{k=1} \frac{H_k}{k(i+k)} \right) = \sum^{\infty}_{k=1} \frac{H_k}{k} \left( \sum^{\infty}_{i=1} \frac{h_k}{(2i-1)(i+k)} \right)~, 
\end{eqnarray}
and inserting Eq.~(15) in Eq.~(20) an explicit calculation of the different Euler sums results in Eq.~(18).   
This way all cubic Euler sums of degree four are expressible in terms of zeta values and the polylogarithmic value $Li_4(1/2)$, where the odd-type linear Euler sum $\sum^{\infty}_{k=1}\frac{h_k}{k^3}$
can be expressed in terms of the polylogarithmic value $Li_4(1/2)$ \cite{bra20}.

\section{The order five case for cubic Euler sums}
In the following section we demonstrate that all members of the eight families of cubic Euler sums of degree five can be calculated in terms of zeta values and polylogarithmic values
$Li_4(1/2)$ and $Li_5(1/2)$.

\subsection{First family}
In the order five case two members of the first family exist. We start with: 
\begin{eqnarray}
\sum^{\infty}_{k=1} \frac{H_k h_k^2}{k^2} &=& \frac{155}{32}\zeta(5) + \frac{7}{8}\zeta(2)\zeta(3)~.
\end{eqnarray}
For an explicit calculation it is helpful to introduce first a new identity that results as the solution of the corresponding first order difference equation:
\begin{eqnarray}
\sum^{\infty}_{k=1} \frac{h_k^2}{k(i-k)} = \sum^{\infty}_{k=1} \frac{h_k^2}{k(i+k)} - \frac{h_i^2}{i^2} - 4\frac{h_i h_i^{(2)}}{i}~.
\end{eqnarray}
This nonlinear two-valued identity can now be used as follows:
\begin{eqnarray}
\sum^{\infty}_{i=1} \frac{1}{i} \left( \sum^{\infty}_{k=1} \frac{h_k^2}{k(i-k)} \right) = \sum^{\infty}_{i=1} \frac{1}{i} \left( \sum^{\infty}_{k=1} \frac{h_k^2}{k(i+k)} \right)
- \sum^{\infty}_{i=1} \frac{h_i^2}{i^3} - 4\frac{h_i h_i^{(2)}}{i^2}~.
\end{eqnarray}
Rearranging the summations on both sides of Eq.~(23) we get:
\begin{eqnarray}
\sum^{\infty}_{k=1} \frac{h_k^2}{k} \left( \sum^{\infty}_{i=1} \frac{1}{i(i-k)} \right) = \sum^{\infty}_{k=1} \frac{h_k^2}{k} \left( \sum^{\infty}_{i=1} \frac{1}{i(k+i)} \right)
- \sum^{\infty}_{i=1} \frac{h_i^2}{i^3} - 4\frac{h_i h_i^{(2)}}{i^2}~.
\end{eqnarray}
It follows further
\begin{eqnarray}
\sum^{\infty}_{k=1} \frac{h_k^2}{k} \left( \frac{2}{k^2} - \frac{H_k}{k} \right) = \sum^{\infty}_{k=1} \frac{h_k^2}{k} \left( \frac{H_k}{k} \right)
- \sum^{\infty}_{i=1} \frac{h_i^2}{i^3} - 4\frac{h_i h_i^{(2)}}{i^2}~.
\end{eqnarray}
This way it results
\begin{eqnarray}
\sum^{\infty}_{k=1} \frac{H_k h_k^2}{k^2} =  \frac{3}{2}\sum^{\infty}_{k=1} \frac{h_k^2}{k^3} +  2\sum^{\infty}_{k=1} \frac{h_k h_k^{(2)}}{k^2}~. 
\end{eqnarray}
As both Euler sums on the right side are known explicitly Eq.~(21) follows.

The next Euler sum is defined as: 
\begin{eqnarray}
\sum^{\infty}_{k=1} \frac{H_k^2 h_k}{k^2} &=& \frac{341}{16}\zeta(5) - \frac{53}{4}ln(2)\zeta(4) + \frac{7}{4}\zeta(2)\zeta(3) + 7\left(ln(2)\right)^2\zeta(3) -
\frac{4}{3}\left(ln(2)\right)^3\zeta(2) + \frac{2}{15}\left(ln(2)\right)^5  \nonumber \\ &-& 2ln(2)\sum^{\infty}_{k=1}\frac{h_k}{k^3} - 16 Li_5\left(\frac{1}{2}\right)~.
\end{eqnarray}

The calculational procedure for this sum starts with the following expression:
\begin{eqnarray}
\sum^{\infty}_{i=1} \frac{H_i}{i} \left( \sum^{\infty}_{k=1} \frac{h_k}{k(i-k)} \right) = \sum^{\infty}_{i=1} \frac{H_i}{i} \left( \sum^{\infty}_{k=1} \frac{h_k}{k(i+k)} \right)
- \sum^{\infty}_{i=1} \frac{H_i h_i}{i^3} - 2\frac{H_i h_i^{(2)}}{i^2}~.
\end{eqnarray}
Here we used the identity defined as lemma 1 from \cite{bra21}. Rearranging the summations on both sides of Eq.~(28) we get:
\begin{eqnarray}
-\sum^{\infty}_{k=1} \frac{h_k}{k} \left( \sum^{\infty}_{i=1} \frac{H_i}{i(k-i)} \right) = \sum^{\infty}_{k=1} \frac{h_k}{k} \left( \sum^{\infty}_{i=1} \frac{H_k}{i(k+i)} \right)
- \sum^{\infty}_{k=1} \frac{H_k h_k}{k^3} - 2\frac{H_k h_k^{(2)}}{k^2}~.
\end{eqnarray}
Using lemma 1 from \cite{bra22} and Eq.~(5) it follows further:
\begin{eqnarray}
-\sum^{\infty}_{k=1} \frac{h_k}{k} \left( \sum^{\infty}_{i=1} \frac{H_i}{i(k+i)} - \zeta(2)\frac{1}{k} -2\frac{H_k^{(2)}}{k} \right) = \sum^{\infty}_{k=1} \frac{h_k}{k} 
\left( \sum^{\infty}_{i=1} \frac{H_k}{i(k+i)} \right) - \sum^{\infty}_{k=1} \frac{H_k h_k}{k^3} - 2\frac{H_k h_k^{(2)}}{k^2}~.
\end{eqnarray}
Rearranging the different sums we get:
\begin{eqnarray}
\sum^{\infty}_{k=1} \frac{h_k}{k} \left( \sum^{\infty}_{i=1} \frac{H_i}{i(k+i)}\right) = \frac{1}{2}\zeta(2)\sum^{\infty}_{k=1} \frac{h_k}{k^2} + \sum^{\infty}_{k=1} \frac{H_k^{(2)} h_k}{k^2}
+  \sum^{\infty}_{k=1} \frac{H_k h_k^{(2)}}{k^2} + \frac{1}{2}\sum^{\infty}_{k=1} \frac{H_k h_k}{k^3}~.
\end{eqnarray}
This is equal to
\begin{eqnarray}
\sum^{\infty}_{k=1} \frac{h_k}{k} \left( \zeta(2) \frac{1}{k} -  \frac{H_k}{k^2} + \frac{1}{2} \frac{H_k^{(2)}}{k} +  \frac{1}{2}  \frac{H_k^2}{k} \right) = 
\frac{1}{2}\zeta(2)\sum^{\infty}_{k=1} \frac{h_k}{k^2} + \sum^{\infty}_{k=1} \frac{H_k^{(2)} h_k}{k^2}
+  \sum^{\infty}_{k=1} \frac{H_k h_k^{(2)}}{k^2} +  \frac{1}{2}\sum^{\infty}_{k=1} \frac{H_k h_k}{k^3}~.
\end{eqnarray}
Finally we arrive at:
\begin{eqnarray}
\sum^{\infty}_{k=1} \frac{H_k^2 h_k}{k^2} = - \zeta(2) \frac{h_k}{k^2} + \sum^{\infty}_{k=1} \frac{H_k^{(2)} h_k}{k^2} + 2 \sum^{\infty}_{k=1} \frac{H_k h_k^{(2)}}{k^2} +
3\sum^{\infty}_{k=1} \frac{H_k h_k}{k^3}~.
\end{eqnarray}
All Euler sums on the right side of  Eq.~(33) are known explicitly. Inserting the corresponding values Eq.~(27) follows.

\subsection{Fourth family}
It is favorably to continue with the members of the forth and fifth family before we evaluate the members of the second family. First we introduce the only member of the fourth family. 
It follows
\begin{eqnarray}
\sum^{\infty}_{k=1} \frac{h_k^3}{(2k-1)^2} = - \frac{713}{256}\zeta(5) + \frac{249}{64}ln(2)\zeta(4) + \frac{3}{32}\zeta(2)\zeta(3) + \left(ln(2)\right)^3\zeta(2) -
\frac{1}{40}\left(ln(2)\right)^5 + 3Li_5\left(\frac{1}{2}\right)~.
\end{eqnarray}
This result is known from literature \cite{zeh07}~.

\subsection{Fifth family}
In the order five case two cubic sums exist in the fifth family. The first one is:
\begin{eqnarray}
\sum^{\infty}_{k=1} \frac{H_k^2 h_k}{(2k-1)^2} &=& \frac{651}{32}\zeta(5) - \frac{257}{16}ln(2)\zeta(4) + \frac{45}{8}\zeta(4)- \frac{7}{16}\zeta(2)\zeta(3) + \frac{7}{4}\left(ln(2)\right)^2\zeta(3)
+ \frac{7}{2}ln(2)\zeta(3) - \frac{5}{2}\zeta(3) \nonumber \\ &+& \frac{5}{3}\left(ln(2)\right)^3\zeta(2) - 6\left(ln(2)\right)^2\zeta(2) + 7ln(2)\zeta(2) - 4\zeta(2) +
\frac{2}{15}\left(ln(2)\right)^5 - \sum^{\infty}_{k=1}\frac{h_k}{k^3} - 16 Li_5\left(\frac{1}{2}\right)~.
\end{eqnarray}

The above result can be obtained from:
\begin{eqnarray}
\sum^{\infty}_{i=1} \frac{h_i}{2i-1} \left( \sum^{\infty}_{k=1} \frac{H_k}{(2k+1)(i+k)} \right) = \sum^{\infty}_{k=1} \frac{H_k}{2k+1} \left( \sum^{\infty}_{i=1} \frac{h_k}{(2i-1)(i+k)} \right)~.
\end{eqnarray}
The two-valued help functions on the left and right side of Eq.~(36) had been defined before. This way we get:
\begin{eqnarray}
\left( \zeta(2) - 2\left(ln(2)\right)^2 \right) \sum^{\infty}_{i=1} \frac{h_i}{(2i-1)^2} + \frac{1}{2} \sum^{\infty}_{i=1} \frac{H_{i-1}^{(2)} h_i}{(2i-1)^2} +
\frac{1}{2} \sum^{\infty}_{i=1} \frac{H_{i-1}^2 h_i}{(2i-1)^2} &=& 2ln(2)\sum^{\infty}_{i=1} \frac{H_i h_i}{i(2i+1)} - \zeta(2)\sum^{\infty}_{i=1} \frac{H_i}{(2i+1)^2} \nonumber \\ &+&
\sum^{\infty}_{i=1} \frac{H_i}{(2i+1)^2} \left( \sum^{i}_{k=1} \frac{H_{k-1}}{2k-1} \right)~.
\end{eqnarray}
For the last sum on the right side of Eq.~(37) we can write
\begin{eqnarray}
\sum^{\infty}_{i=1} \frac{H_i}{(2i+1)^2} \left( \sum^{i}_{k=1} \frac{H_{k-1}}{2k-1} \right) = \sum^{\infty}_{i=1} \frac{H_i^2 h_i}{(2i+1)^2} - \sum^{\infty}_{i=1} \frac{H_i}{(2i+1)^2} 
\left( \sum^{i}_{k=1} \frac{h_k}{k} \right)~.
\end{eqnarray}

This way it remains to calculate 
\begin{eqnarray}
\sum^{\infty}_{i=1} \frac{H_i}{(2i+1)^2} \left( \sum^{i}_{k=1} \frac{h_k}{k} \right) = \sum^{\infty}_{i=1} \frac{H_i}{2i+1} \left( \sum^{\infty}_{k=1} \frac{h_k}{k(2i+2k-1)} \right) - 
\frac{1}{2}\zeta(2) \sum^{\infty}_{i=1} \frac{H_i}{(2i+1)^2}
\end{eqnarray}
or
\begin{eqnarray}
\sum^{\infty}_{i=1} \frac{H_i}{(2i+1)^2} \left( \sum^{i}_{k=1} \frac{h_k}{k} \right) = \sum^{\infty}_{i=1} \frac{h_i}{i} \left( \sum^{\infty}_{k=1} \frac{H_k}{(2k+1)(2i+2k-1)} \right) - 
\frac{1}{2}\zeta(2) \sum^{\infty}_{i=1} \frac{H_i}{(2i+1)^2}
\end{eqnarray}
or
\begin{eqnarray}
\sum^{\infty}_{i=1} \frac{H_i}{(2i+1)^2} \left( \sum^{i}_{k=1} \frac{h_k}{k} \right) = \sum^{\infty}_{i=1} \frac{h_i}{i} \left( \frac{h_i^{(2)}}{2i} +  \frac{h_i^2}{2i} 
- ln(2)\frac{h_i}{i} \right) - \frac{1}{2}\zeta(2) \sum^{\infty}_{i=1} \frac{H_i}{(2i+1)^2}~.
\end{eqnarray}
The result is:
\begin{eqnarray}
\sum^{\infty}_{i=1} \frac{H_i}{(2i+1)^2} \left( \sum^{i}_{k=1} \frac{h_k}{k} \right) = \frac{31}{16}\zeta(5) - \frac{15}{16}ln(2)\zeta(4)~. 
\end{eqnarray}
With this an explicit computation of Eq.~(37) is possible, and this way Eq.~(35) results.

The second Euler sum belonging to this family is
\begin{eqnarray}
\sum^{\infty}_{k=1} \frac{H_k h_k^2}{(2k-1)^2} &=& -\frac{155}{64}\zeta(5) + \frac{83}{16}ln(2)\zeta(4) + \frac{15}{16}\zeta(4) - \frac{7}{32}\zeta(2)\zeta(3) - \frac{7}{4}\left(ln(2)\right)^2\zeta(3)
+ \frac{7}{4}ln(2)\zeta(3) - \frac{3}{8}\zeta(3) \nonumber \\ &-& \frac{7}{6}\left(ln(2)\right)^3\zeta(2) + \frac{3}{2}\left(ln(2)\right)^2\zeta(2) - \frac{3}{2}ln(2)\zeta(2) -
\frac{1}{30}\left(ln(2)\right)^5 + \frac{1}{4}ln(2)\sum^{\infty}_{k=1}\frac{h_k}{k^3} - \frac{1}{4}\sum^{\infty}_{k=1}\frac{h_k}{k^3} \nonumber \\ &+& 4Li_5\left(\frac{1}{2}\right)~.
\end{eqnarray}

The corresponding calculational scheme is based an the following expression:
\begin{eqnarray}
\sum^{\infty}_{k=1} \frac{H_k}{2k-1} \left( \sum^{\infty}_{i=1} \frac{H_i}{i(2i+2k-1)} \right) &=& 2\left(ln(2)\right)^2 \sum^{\infty}_{k=1}\frac{H_k}{(2k-1)^2} - 
4\sum^{\infty}_{k=1}\frac{H_k h_k}{(2k-1)^3} - 4ln(2)\sum^{\infty}_{k=1}\frac{H_k h_k}{(2k-1)^2} \nonumber \\ &+& 2\sum^{\infty}_{k=1}\frac{H_k h_k^{(2)}}{(2k-1)^2}
+ 4ln(2)\sum^{\infty}_{k=1}\frac{H_k}{(2k-1)^3} + 2\sum^{\infty}_{k=1}\frac{H_k h_k^2}{(2k-1)^2}
\end{eqnarray}
or
\begin{eqnarray}
\sum^{\infty}_{k=1} \frac{H_k}{k} \left( \sum^{\infty}_{i=1} \frac{H_i}{(2i-1)(2i+2k-1)} \right)  &=& 2ln(2) \sum^{\infty}_{k=1}\frac{H_k}{k(2k-1)} - ln(2)\sum^{\infty}_{k=1}\frac{H_k h_k}{k^2}-
\sum^{\infty}_{k=1}\frac{H_k h_k}{k^2(2k-1)} + \frac{1}{2}\sum^{\infty}_{k=1}\frac{H_k h_k^{(2)}}{k^2} \nonumber \\ &+& \frac{1}{2}\sum^{\infty}_{k=1}\frac{H_k h_k^2}{k^2}~.
\end{eqnarray}
All Euler sums are known explicitly. This way Eq.~(43) results, where one member of the first family is needed. 

\subsection{Second family}
The only member of the second family is:
\begin{eqnarray}
\sum^{\infty}_{k=1} \frac{H_k^3}{(2k-1)^2} &=& -\frac{1271}{32}\zeta(5) + \frac{183}{8}ln(2)\zeta(4) + \frac{135}{8}\zeta(4) + \frac{7}{4}\zeta(2)\zeta(3) - 42ln(2)\zeta(3) + \frac{35}{2}\zeta(3) -
2\left(ln(2)\right)^3\zeta(2) \nonumber \\ &+& 18\left(ln(2)\right)^2\zeta(2) - 12ln(2)\zeta(2) - 2\zeta(2) - 24ln(2) + 24\left(ln(2)\right)^2 - 8\left(ln(2)\right)^3 - 
\frac{2}{5}\left(ln(2)\right)^5 \nonumber \\ &+& 6\sum^{\infty}_{k=1}\frac{h_k}{k^3} + 48Li_5\left(\frac{1}{2}\right)~.
\end{eqnarray}

In this case we start the calculation with the decomposition of the well known cubic Euler sum
\begin{eqnarray}
\sum^{\infty}_{k=1} \frac{H_k^3}{k^2} = \sum^{\infty}_{k=1} \frac{H_{2k-1}^3}{(2k-1)^2} + \frac{1}{4}\sum^{\infty}_{k=1} \frac{H_{2k}^3}{k^2}
\end{eqnarray}
or
\begin{eqnarray}
\sum^{\infty}_{k=1} \frac{H_k^3}{k^2} &=& \sum^{\infty}_{k=1} \frac{(h_k + \frac{1}{2}H_k)^3}{(2k-1)^2} - \frac{3}{2}\sum^{\infty}_{k=1} \frac{(h_k + \frac{1}{2}H_k)^2}{k(2k-1)^2} +
\frac{3}{4}\sum^{\infty}_{k=1} \frac{h_k + \frac{1}{2}H_k}{k^2(2k-1)^2} + \frac{1}{4}\sum^{\infty}_{k=1} \frac{(h_k + \frac{1}{2}H_k)^3}{k^2} \nonumber \\ &-& \frac{1}{8}
\sum^{\infty}_{k=1} \frac{1}{k^3(2k-1)^2}~.
\end{eqnarray}
From this expression follows after some simple algebraic manipulations the corresponding member of the second family as a function of linear, quadratic and cubic Euler sums which are known explicitly,
or have been calculated before. This way Eq.~(46) results.

\subsection{Third family}
Again, for the third family only one member exists:
\begin{eqnarray}
\sum^{\infty}_{k=1} \frac{H_k^2H_k^{(2)}}{k(2k-1)} &=& \frac{93}{4}\zeta(5) - 9ln(2)\zeta(4) - \frac{35}{2}\zeta(4) - 2\zeta(2)\zeta(3) + 28ln(2)\zeta(3) - 8\zeta(3) -
8\left(ln(2)\right)^2\zeta(2) \nonumber \\ &+& 8\zeta(2) + 32ln(2) - 32\left(ln(2)\right)^2  + \frac{32}{3}\left(ln(2)\right)^3 + 
\frac{4}{15}\left(ln(2)\right)^5 - 4\sum^{\infty}_{k=1}\frac{h_k}{k^3} - 32Li_5\left(\frac{1}{2}\right)~.
\end{eqnarray}
We start the corresponding calculation with an expression discussed by Adegoke \cite{ade16}:
\begin{eqnarray}
\sum^{k}_{i=1} \frac{H_i^{(2)}}{i} = H_k H_k^{(2)} + H_k^{(3)} - \sum^{k}_{i=1}  \frac{H_i}{i^2}~.
\end{eqnarray}
From this it follows:
\begin{eqnarray}
\sum^{\infty}_{k=1} \frac{H_k}{k(2k-1)} \left( \sum^{k}_{i=1} \frac{H_i^{(2)}}{i} \right) = \sum^{\infty}_{k=1} \frac{H_k^2 H_k^{(2)}}{k(2k-1)} + \sum^{\infty}_{k=1} \frac{H_k H_k^{(3)}}{k(2k-1)} - 
\sum^{\infty}_{k=1} \frac{H_k}{k(2k-1)} \left( \sum^{k}_{i=1} \frac{H_i}{i^2} \right)~.
\end{eqnarray}
Now we calculate the double sum on the right side first. It follows:
\begin{eqnarray}
\sum^{\infty}_{k=1} \frac{H_k}{k(2k-1)} \left( \sum^{k}_{i=1} \frac{H_i}{i^2} \right) = \sum^{\infty}_{k=1} \frac{H_k^2}{k^3(2k-1)} - \sum^{\infty}_{k=1} \frac{H_k}{k^2}
\left( \sum^{k}_{i=1} \frac{H_i}{i(2i-1)} \right) 
\end{eqnarray}
or
\begin{eqnarray}
\sum^{\infty}_{k=1} \frac{H_k}{k(2k-1)} \left( \sum^{k}_{i=1} \frac{H_i}{i^2} \right) = \sum^{\infty}_{k=1} \frac{H_k^2}{k^3(2k-1)} - 2\sum^{\infty}_{k=1} \frac{H_k}{k^2}
\left( \sum^{k}_{i=1} \frac{H_i}{2i-1} \right)  + \sum^{\infty}_{k=1} \frac{H_k}{k^2} \left( \sum^{k}_{i=1} \frac{H_i}{i} \right)
\end{eqnarray}
or
\begin{eqnarray}
\sum^{\infty}_{k=1} \frac{H_k}{k(2k-1)} \left( \sum^{k}_{i=1} \frac{H_i}{i^2} \right) = \sum^{\infty}_{k=1} \frac{H_k^2}{k^3(2k-1)} - 2\sum^{\infty}_{k=1} \frac{H_k}{k^2}
\left( \sum^{k}_{i=1} \frac{H_i}{2i-1} \right)  + \frac{1}{2}\sum^{\infty}_{k=1} \frac{H_k^3}{k^2} + \frac{1}{2}\sum^{\infty}_{k=1} \frac{H_k  H_k^{(2)}}{k^2}~.
\end{eqnarray}
It remains to calculate the double sum
\begin{eqnarray}
\sum^{\infty}_{k=1} \frac{H_k}{k^2} \left( \sum^{k}_{i=1} \frac{H_i}{2i-1} \right) = \sum^{\infty}_{k=1} \frac{H_k}{k^2} \left( \sum^{k}_{i=1} \frac{H_{i-1}}{2i-1} \right) + 
2\sum^{\infty}_{k=1} \frac{H_k h_k}{k^2} - \sum^{\infty}_{k=1} \frac{H_k^2}{k^2}
\end{eqnarray}
or
\begin{eqnarray}
\sum^{\infty}_{k=1} \frac{H_k}{k^2} \left( \sum^{k}_{i=1} \frac{H_i}{2i-1} \right) = 2\sum^{\infty}_{k=1} \frac{H_k h_k}{k^2} - \sum^{\infty}_{k=1} \frac{H_k^2}{k^2} 
+ \sum^{\infty}_{k=1} \frac{H_k}{k} \left( \sum^{\infty}_{i=1} \frac{h_i}{i(k+i)} \right) - 2ln(2)\sum^{\infty}_{k=1} \frac{H_k h_k}{k^2} 
\end{eqnarray}
or
\begin{eqnarray}
\sum^{\infty}_{k=1} \frac{H_k}{k^2} \left( \sum^{k}_{i=1} \frac{H_i}{2i-1} \right) = 2\sum^{\infty}_{k=1} \frac{H_k h_k}{k^2} - 2ln(2)\sum^{\infty}_{k=1} \frac{H_k h_k}{k^2} 
- \sum^{\infty}_{k=1} \frac{H_k^2}{k^2} + \sum^{\infty}_{i=1} \frac{h_i}{i} \left( \sum^{\infty}_{k=1} \frac{H_k}{k(i+k)} \right)~. 
\end{eqnarray}
Using Eq.~(5) we finally get:
\begin{eqnarray}
\sum^{\infty}_{k=1} \frac{H_k}{k^2} \left( \sum^{k}_{i=1} \frac{H_i}{2i-1} \right) &=& 2\sum^{\infty}_{k=1} \frac{H_k h_k}{k^2} - 2ln(2)\sum^{\infty}_{k=1} \frac{H_k h_k}{k^2} 
- \sum^{\infty}_{k=1} \frac{H_k^2}{k^2} + \zeta(2) \sum^{\infty}_{k=1} \frac{h_k}{k^2} - \sum^{\infty}_{k=1} \frac{H_k h_k}{k^3} \nonumber \\ &+& \frac{1}{2} \sum^{\infty}_{k=1} \frac{H_k^{(2)} h_k}{k^2} +
\frac{1}{2} \sum^{\infty}_{k=1} \frac{H_k^2 h_k}{k^2}~.
\end{eqnarray}
This results in
\begin{eqnarray}
\sum^{\infty}_{k=1} \frac{H_k}{k^2} \left( \sum^{k}_{i=1} \frac{H_i}{2i-1} \right) &=& \frac{155}{32}\zeta(5) - \frac{45}{8}ln(2)\zeta(4) + \frac{11}{8}\zeta(4) + \frac{7}{4}\zeta(2)\zeta(3)
- ln(2)\sum^{\infty}_{k=1}\frac{h_k}{k^3} + \sum^{\infty}_{k=1}\frac{h_k}{k^3} \nonumber \\ &-& \frac{1}{2} \sum^{\infty}_{k=1} \frac{H_k h_k}{k^3}~. 
\end{eqnarray}

The last calculational step consists in the computation of the double sum appearing on the left side of Eq.~(51). We start with Eq.~(126) from \cite{bra20}: 
\begin{eqnarray}
\sum^{k-1}_{i=1} \frac{H_i^{(2)}}{i} = \zeta(3) - \zeta(2)\frac{1}{k} + \zeta(2)H_k - \sum^{k-1}_{i=1} \frac{H_i}{i^2} - \sum^{\infty}_{i=1} \frac{H_i}{(i+k)^2}~.
\end{eqnarray} 
From this we get:
\begin{eqnarray}
\sum^{\infty}_{k=1} \frac{H_k}{k(2k-1)} \left( \sum^{k-1}_{i=1} \frac{H_i^{(2)}}{i} \right) &=& \zeta(3)\sum^{\infty}_{k=1} \frac{H_k}{k(2k-1)} - \zeta(2)\sum^{\infty}_{k=1}\frac{H_k}{k^2(2k-1)} +
\zeta(2)\sum^{\infty}_{k=1} \frac{H_k^2}{k(2k-1)} - \sum^{\infty}_{k=1} \frac{H_k}{k(2k-1)} \left( \sum^{k-1}_{i=1} \frac{H_i}{i^2}  \right) \nonumber \\ &-& \sum^{\infty}_{k=1} \frac{H_k}{k(2k-1)} \left(
\sum^{\infty}_{i=1} \frac{H_i}{(i+k)^2} \right)~. 
\end{eqnarray}
It remains to calculate the last double sum on the right side of Eq.~(61) as the fourth double sum on the right side is defined by Eq.~(54). It follows first:
\begin{eqnarray}
\sum^{\infty}_{k=1} \frac{H_k}{k(2k-1)} \left( \sum^{\infty}_{i=1} \frac{H_i}{(i+k)^2} \right) = \sum^{\infty}_{i=1} H_i \left( \sum^{\infty}_{k=1} \frac{H_k}{k(2k-1)(i+k)^2} \right) 
\end{eqnarray}
and by partial fraction decomposition we get:
\begin{eqnarray}
\sum^{\infty}_{k=1} \frac{H_k}{k(2k-1)} \left( \sum^{\infty}_{i=1} \frac{H_i}{(i+k)^2} \right) &=& 4\sum^{\infty}_{i=1} \frac{H_i}{2i+1} \left( \sum^{\infty}_{k=1} \frac{H_k}{(2k-1)(i+k)} \right) -
\sum^{\infty}_{i=1} \frac{H_i}{i} \left( \sum^{\infty}_{k=1} \frac{H_k}{k(i+k)} \right) \nonumber \\ &+& \sum^{\infty}_{i=1} \frac{H_i}{i} \left( \sum^{\infty}_{k=1} \frac{H_k}{(i+k)^2} \right)
- 2\sum^{\infty}_{i=1} \frac{H_i}{2i+1} \left( \sum^{\infty}_{k=1} \frac{H_k}{(i+k)^2} \right)~.
\end{eqnarray}
The first and the second double sum can be calculated straightforwardly by use of the corresponding two-valued help functions. The third sum can be calculated in the following way: 
\begin{eqnarray}
\sum^{\infty}_{i=1} \frac{H_i}{i} \left( \sum^{\infty}_{k=1} \frac{H_k}{(i+k)^2} \right) = \sum^{\infty}_{k=1} H_k  \left( \sum^{\infty}_{i=1} \frac{H_i}{i(k+i)^2} \right)
\end{eqnarray}
or
\begin{eqnarray}
\sum^{\infty}_{i=1} \frac{H_i}{i} \left( \sum^{\infty}_{k=1} \frac{H_k}{(i+k)^2} \right) = \sum^{\infty}_{k=1} \frac{H_k}{k}  \left( \sum^{\infty}_{i=1} \frac{H_i}{i(k+i)} \right) -
\sum^{\infty}_{k=1} \frac{H_k}{k}  \left( \sum^{\infty}_{i=1} \frac{H_i}{(k+i)^2} \right)
\end{eqnarray}
or
\begin{eqnarray}
\sum^{\infty}_{i=1} \frac{H_i}{i} \left( \sum^{\infty}_{k=1} \frac{H_k}{(i+k)^2} \right) = \frac{1}{2}\sum^{\infty}_{k=1} \frac{H_k}{k} \left( \sum^{\infty}_{i=1} \frac{H_i}{i(k+i)} \right)~.
\end{eqnarray}
Again, the double sum on the right side of Eq.~(66) can be calculated straightforwardly by use of the corresponding two-valued help function. The fourth sum on the right side of  Eq.~(63) can be
recalculated by an appropriate shift of the summation indices. It follows after some tedious algebra:
\begin{eqnarray}
\sum^{\infty}_{i=1} \frac{H_i}{2i+1} \left( \sum^{\infty}_{k=1} \frac{H_k}{(i+k)^2} \right) &=& 5\zeta(4) - \zeta(3) -4ln(2)\zeta(2) + \frac{16}{3}\left(ln(2)\right)^3 - 
2\sum^{\infty}_{k=1}\frac{h_k}{k^3} \nonumber \\ &+& \sum^{\infty}_{i=1} \frac{H_i}{2i-1} \left( \sum^{\infty}_{k=1} \frac{H_k}{(i+k)^2} \right)~. 
\end{eqnarray}
This way we get:
\begin{eqnarray}
\sum^{\infty}_{k=1} \frac{H_k}{k(2k-1)} \left( \sum^{\infty}_{i=1} \frac{H_i}{(i+k)^2} \right) &=& -\frac{473}{8}\zeta(5) + \frac{61}{2}ln(2)\zeta(4) + \frac{25}{4}\zeta(4) + 5\zeta(2)\zeta(3)
-14\left(ln(2)\right)^2\zeta(3) - 14ln(2)\zeta(3) \nonumber \\ &+& \frac{16}{3}\left(ln(2)\right)^3\zeta(2) - \frac{8}{15}\left(ln(2)\right)^5  + 6\sum^{\infty}_{k=1}\frac{h_k}{k^3} +
64Li_5\left(\frac{1}{2}\right)
\end{eqnarray}
and Eq.~(49) finally results.

\subsection{Sixth family}
Also, for the sixth family we have only one member:
\begin{eqnarray}
\sum^{\infty}_{k=1} \frac{h_k^2h_k^{(2)}}{k(2k-1)} &=& -\frac{93}{64}\zeta(5) + \frac{49}{16}ln(2)\zeta(4) + \frac{3}{8}\zeta(2)\zeta(3) - \frac{1}{3} \left(ln(2)\right)^3\zeta(2)
- \frac{1}{60}\left(ln(2)\right)^5 + 2Li_5\left(\frac{1}{2}\right)~.
\end{eqnarray}
The calculational procedure starts with:
\begin{eqnarray}
\sum^{\infty}_{k=1} \frac{h_k^{(2)}}{k} \left( \sum^{\infty}_{i=1} \frac{H_i}{i(2k+2i-1)} \right) &=& 2\left(ln(2)\right)^2 \sum^{\infty}_{k=1} \frac{h_k^{(2)}}{k(2k-1)} +
4ln(2) \sum^{\infty}_{k=1} \frac{h_k^{(2)}}{k(2k-1)^2} -4ln(2) \sum^{\infty}_{k=1} \frac{h_k h_k^{(2)}}{k(2k-1)} \nonumber \\ &-& 4\sum^{\infty}_{k=1} \frac{h_k h_k^{(2)}}{k(2k-1)^2} + 
2\sum^{\infty}_{k=1} \frac{h_k^{(2)} h_k^{(2)}}{k(2k-1)} + \sum^{\infty}_{k=1} \frac{h_k^2h_k^{(2)}}{k(2k-1)}~.
\end{eqnarray}
Rearranging the summations on the left side of Eq.~(70) and using the corresponding two-valued function 
\begin{eqnarray}
\sum^{\infty}_{i=1} \frac{h_i^{(2)}}{i(2k+2i-1)} = \left( \frac{3}{2}ln(2)\zeta(2) - \frac{7}{8}\zeta(3) \right) \frac{1}{2k-1} + \frac{3}{4}\zeta(2) \frac{H_{k-1}}{2k-1} - \frac{1}{2}
\frac{1}{2k-1}  \sum^{k-1}_{i=1}  \frac{h_i}{i^2}
\end{eqnarray}
it follows after some simple algebra:
\begin{eqnarray}
\sum^{\infty}_{k=1} \frac{h_k^2 h_k^{(2)}}{k(2k-1)} &=& 2ln(2) \sum^{\infty}_{k=1} \frac{h_k h_k^{(2)}}{k(2k-1)} - \left(ln(2)\right)^2 \sum^{\infty}_{k=1} \frac{h_k^{(2)}}{k(2k-1)} -
2ln(2) \sum^{\infty}_{k=1} \frac{h_k^{(2)}}{k(2k-1)^2} + 2\sum^{\infty}_{k=1} \frac{h_k h_k^{(2)}}{k(2k-1)^2} \nonumber \\ &-& \sum^{\infty}_{k=1} \frac{h_k^{(2)} h_k^{(2)}}{k(2k-1)} +
\left( \frac{3}{4}ln(2)\zeta(2) - \frac{7}{16}\zeta(3) \right) \frac{H_k}{k(2k-1)} + \frac{3}{8}\zeta(2) \frac{H_k H_{k-1}}{k(2k-1)} \nonumber \\ &-& \frac{1}{4} 
\sum^{\infty}_{k=1} \frac{H_k}{k(2k-1)} \left( \sum^{k-1}_{i=1} \frac{h_i}{i^2} \right)~.
\end{eqnarray}
In Eq.~(72) only the last Euler sum on the right side is not known explicitly. The calculational procedure starts with: 
\begin{eqnarray}
\sum^{\infty}_{k=1} \frac{H_k}{k(2k-1)} \left( \sum^{k-1}_{i=1} \frac{h_i}{i^2} \right) = \sum^{\infty}_{k=1} \frac{H_k}{k(2k-1)} \left( \sum^{\infty}_{i=1} \frac{h_i}{i^2} \right) -
\sum^{\infty}_{k=1} \frac{h_k}{k^2} \left( \sum^{k}_{i=1} \frac{H_i}{i(2i-1)} \right)
\end{eqnarray}
or
\begin{eqnarray}
\sum^{\infty}_{k=1} \frac{H_k}{k(2k-1)} \left( \sum^{k-1}_{i=1} \frac{h_i}{i^2} \right) = \sum^{\infty}_{k=1} \frac{H_k}{k(2k-1)} \left( \sum^{\infty}_{i=1} \frac{h_i}{i^2} \right) -
2\sum^{\infty}_{k=1} \frac{h_k}{k^2} \left( \sum^{k}_{i=1} \frac{H_i}{2i-1} \right) + \frac{1}{2} \sum^{\infty}_{k=1} \frac{H_k^2 h_k}{k^2} + \frac{1}{2} \sum^{\infty}_{k=1} \frac{H^{(2)}_k h_k}{k^2}~.
\end{eqnarray}
This way it remains to calculate the second sum:
\begin{eqnarray}
\sum^{\infty}_{k=1} \frac{h_k}{k^2} \left( \sum^{k}_{i=1} \frac{H_i}{2i-1} \right) = \sum^{\infty}_{k=1} \frac{H_k h_k^2}{k^2} - \sum^{\infty}_{k=1} \frac{h_k}{k^2}
\left( \sum^{k}_{i=1} \frac{h_i}{i} \right) - \sum^{\infty}_{k=1} \frac{H_k h_k}{k^2} + 2\sum^{\infty}_{k=1} \frac{h_k^2}{k^2}~.
\end{eqnarray}
Again the second sum on the right side of Eq.~(75) is unknown. Starting with:
\begin{eqnarray}
\sum^{\infty}_{k=1} \frac{h_k}{k} \left( \sum^{\infty}_{i=1} \frac{h_i}{i(k-i)} \right) =  \sum^{\infty}_{k=1} \frac{h_k}{k} \left( \sum^{\infty}_{i=1} \frac{h_i}{i(k+i)} \right) - 
\sum^{\infty}_{k=1} \frac{h_k^2}{k^3} - 2\sum^{\infty}_{k=1} \frac{h_k h_k^{(2)}}{k^2}
\end{eqnarray}
and using lemma 1 from \cite{bra20} we get:
\begin{eqnarray}
\sum^{\infty}_{k=1} \frac{h_k}{k} \left( \sum^{\infty}_{i=1} \frac{h_i}{i(k+i)} \right) = \sum^{\infty}_{k=1} \frac{h_k^2}{k^3} + 2\sum^{\infty}_{k=1} \frac{h_k h_k^{(2)}}{k^2}
\end{eqnarray}
or
\begin{eqnarray}
\sum^{\infty}_{k=1} \frac{h_k}{k} \left( 2ln(2)\frac{h_k}{k} + \frac{H_k h_k}{k} - \frac{1}{k} \sum^{k}_{i=1} \frac{h_i}{i} \right)  = \sum^{\infty}_{k=1} \frac{h_k^2}{k^3} +
2\sum^{\infty}_{k=1} \frac{h_k h_k^{(2)}}{k^2}~.
\end{eqnarray}
From this it follows:
\begin{eqnarray}
\sum^{\infty}_{k=1} \frac{h_k}{k^2} \left( \sum^{k}_{i=1} \frac{h_i}{i} \right) =  -\frac{31}{32}\zeta(5) + \frac{45}{8}ln(2)\zeta(4) + \frac{7}{8}\zeta(2)\zeta(3)~.
\end{eqnarray}
Inserting Eq.~(79) in Eq.~(75) an explicit calculation of Eq.~(73) is possible, and this way Eq.~(69) results.

\subsection{Seventh family}
The seventh family consists of four members. First we have:
\begin{eqnarray}
\sum^{\infty}_{k=1} \frac{H_k^2 h_k^{(2)}}{k(2k-1)} &=& -\frac{403}{16}\zeta(5) + \frac{31}{4}ln(2)\zeta(4) + \frac{105}{8}\zeta(4) - \frac{1}{4}\zeta(2)\zeta(3)
+ \frac{21}{2}\left(ln(2)\right)^2\zeta(3) - 35ln(2)\zeta(3) + \frac{21}{2}\zeta(3)  \nonumber \\ &-& \frac{10}{3}\left(ln(2)\right)^3\zeta(2) +  12\left(ln(2)\right)^2\zeta(2) -
6ln(2)\zeta(2)	- \frac{4}{15}\left(ln(2)\right)^5 - 2ln(2)\sum^{\infty}_{k=1}\frac{h_k}{k^3} + 4\sum^{\infty}_{k=1}\frac{h_k}{k^3} \nonumber \\ &+& 32Li_5\left(\frac{1}{2}\right)~.
\end{eqnarray}
Eq.~(80) results from the following approach:
\begin{eqnarray}
\sum^{\infty}_{k=1} \frac{h_k^{(2)}}{2k-1} \left( \sum^{\infty}_{i=1} \frac{H_i}{i(k+i)} \right) &=& \sum^{\infty}_{k=1} \frac{h_k^{(2)}}{k(2k-1)} - \sum^{\infty}_{k=1} \frac{H_k h_k^{(2)}}{k^2(2k-1)}
+ \frac{1}{2}\sum^{\infty}_{k=1} \frac{H_k^{(2)} h_k^{(2)}}{k(2k-1)} + \frac{1}{2}\sum^{\infty}_{k=1} \frac{H_k^2 h_k^{(2)}}{k(2k-1)}~.
\end{eqnarray}
Rearranging the summations on the left side it follows:
\begin{eqnarray}
\sum^{\infty}_{k=1} \frac{H_k}{k} \left( \sum^{\infty}_{i=1} \frac{h_i^{(2)}}{(2i-1)(k+i)} \right) &=& \left( \frac{21}{8}\zeta(3) - \frac{3}{2}ln(2)\zeta(2) \right) \sum^{\infty}_{k=1} \frac{H_k}{k(2k+1)}
+ \frac{3}{2}\zeta(2) \sum^{\infty}_{k=1} \frac{H_k h_k}{k(2k+1)} - 2ln(2)\sum^{\infty}_{k=1} \frac{H_k h_k^{(2)}}{k(2k+1)} \nonumber \\ &-& \sum^{\infty}_{k=1} \frac{H_k}{k(2k+1)}
\left( \sum^{k}_{i=1} \frac{H_{i-1}}{(2i-1)^2} \right)~. 
\end{eqnarray}
As the left side of Eq.~(81) is equal to the left side of Eq.~(82) we get Eq.~(80) by a rearrangement and an explicit calculation of the different Euler sums. Only the last sum in  Eq.~(82) is not
known explicitly. Here it follows:
\begin{eqnarray}
\sum^{\infty}_{k=1} \frac{H_k}{k(2k+1)} \left( \sum^{k}_{i=1} \frac{H_{i-1}}{(2i-1)^2} \right)  &=& \sum^{\infty}_{k=1} \frac{H_k}{k(2k+1)} \left( \sum^{\infty}_{i=1} \frac{H_{i-1}}{(2i-1)^2} \right)
+ \sum^{\infty}_{k=1} \frac{H_k H_{k-1}}{k(2k+1)(2k-1)^2} \nonumber \\ &-& \sum^{\infty}_{k=1} \frac{H_{k-1}}{(2k-1)^2} \left( \sum^{k}_{i=1} \frac{H_i}{i(2i+1)} \right)~.
\end{eqnarray} 
From this it follows for the last sum in Eq.~(83):
\begin{eqnarray}
\sum^{\infty}_{k=1} \frac{H_{k-1}}{(2k-1)^2} \left( \sum^{k}_{i=1} \frac{H_i}{i(2i+1)} \right)  &=& \sum^{\infty}_{k=1} \frac{H_{k-1}}{(2k-1)^2} \left( \sum^{k}_{i=1} \frac{H_i}{i} \right) -
2\sum^{\infty}_{k=1} \frac{H_{k-1}}{(2k-1)^2} \left( \sum^{k}_{i=1} \frac{H_i}{2i+1} \right)~.
\end{eqnarray}
Again for the last sum on the right side of Eq.~(84) we get:
\begin{eqnarray}
\sum^{\infty}_{k=1} \frac{H_{k-1}}{(2k-1)^2} \left( \sum^{k}_{i=1} \frac{H_i}{2i+1} \right) &=& \sum^{\infty}_{k=1} \frac{H_k H_{k-1}}{(2k+1)(2k-1)^2} + \sum^{\infty}_{k=1} \frac{H_{k-1}}{(2k-1)^2}
\left( \sum^{k-1}_{i=1} \frac{H_i}{2i+1} \right)~.
\end{eqnarray}
Shifting the inner summation index i in the double sum it follows:
\begin{eqnarray}
\sum^{\infty}_{k=1} \frac{H_{k-1}}{(2k-1)^2} \left( \sum^{k}_{i=1} \frac{H_i}{2i+1} \right) &=& \sum^{\infty}_{k=1} \frac{H_k H_{k-1}}{(2k+1)(2k-1)^2} + \sum^{\infty}_{k=1} \frac{H_{k-1}}{(2k-1)^2}
\left( \sum^{k}_{i=1} \frac{H_{i-1}}{2i-1} \right)~.
\end{eqnarray}
Shifting now the outer summation index k in the double sum we get:
\begin{eqnarray}
\sum^{\infty}_{k=1} \frac{H_{k-1}}{(2k-1)^2} \left( \sum^{k}_{i=1} \frac{H_i}{2i+1} \right) &=& \sum^{\infty}_{k=1} \frac{H_k H_{k-1}}{(2k+1)(2k-1)^2} + \sum^{\infty}_{k=1} \frac{H_k}{(2k+1)^2}
\left( \sum^{k+1}_{i=1} \frac{H_{i-1}}{2i-1} \right)
\end{eqnarray}
or
\begin{eqnarray}
\sum^{\infty}_{k=1} \frac{H_{k-1}}{(2k-1)^2} \left( \sum^{k}_{i=1} \frac{H_i}{2i+1} \right) &=& \sum^{\infty}_{k=1} \frac{H_k H_{k-1}}{(2k+1)(2k-1)^2} + \sum^{\infty}_{k=1} \frac{H_k^2}{(2k+1)^3} 
+ \sum^{\infty}_{k=1} \frac{H_k}{(2k+1)^2} \left( \sum^{k}_{i=1} \frac{H_{i-1}}{2i-1} \right)~.
\end{eqnarray}
The last sum has been calculated explicitly before (Eq.~(38)). Thus all sums are known explicitly and this way Eq.~(80) results.

The second member follows to:
\begin{eqnarray}
\sum^{\infty}_{k=1} \frac{H_k H_k^{(2)} h_k}{k(2k-1)} &=& -\frac{403}{16}\zeta(5) + \frac{201}{8}ln(2)\zeta(4) - \frac{95}{8}\zeta(4) + \frac{5}{4}\zeta(2)\zeta(3)
- \frac{21}{2}\left(ln(2)\right)^2\zeta(3) + 7ln(2)\zeta(3) + 5\zeta(3)  \nonumber \\ &+& \frac{4}{3}\left(ln(2)\right)^3\zeta(2) -
8ln(2)\zeta(2) + 8\zeta(2) - \frac{2}{15}\left(ln(2)\right)^5 + ln(2)\sum^{\infty}_{k=1}\frac{h_k}{k^3} + \sum^{\infty}_{k=1}\frac{h_k}{k^3} + 16Li_5\left(\frac{1}{2}\right)~.
\end{eqnarray}

Using the identity introduced in \cite{ade16}:
\begin{eqnarray}
\sum^{k}_{i=1} \frac{H_i}{i^2} + \sum^{k}_{i=1} \frac{H_i^{(2)}}{i} - H_k^{(3)} =  H_k H_k^{(2)}
\end{eqnarray}
it follows
\begin{eqnarray}
\sum^{\infty}_{k=1} \frac{H_k H_k^{(2)} h_k}{k(2k-1)}  &=& \sum^{\infty}_{k=1} \frac{h_k}{k(2k-1)} \left( \sum^{k}_{i=1} \frac{H_i}{i^2} \right) +
\sum^{\infty}_{k=1} \frac{h_k}{k(2k-1)} \left( \sum^{k}_{i=1} \frac{H_i^{(2)}}{i} \right) - \sum^{\infty}_{k=1} \frac{H_k^{(3)} h_k}{k(2k-1)}
\end{eqnarray}
or
\begin{eqnarray}
\sum^{\infty}_{k=1} \frac{H_k H_k^{(2)} h_k}{k(2k-1)}  &=& \sum^{\infty}_{k=1} \frac{h_k}{k(2k-1)} \left( \sum^{k-1}_{i=1} \frac{H_i}{i^2} \right) +
\sum^{\infty}_{k=1} \frac{h_k}{k(2k-1)} \left( \sum^{k-1}_{i=1} \frac{H_i^{(2)}}{i} \right) - \sum^{\infty}_{k=1} \frac{H_k^{(3)} h_k}{k(2k-1)} + \sum^{\infty}_{k=1} \frac{H_k^{(2)} h_k}{k^2(2k-1)}
\nonumber \\ &+& \sum^{\infty}_{k=1} \frac{H_k h_k}{k^3(2k-1)}~.
\end{eqnarray}
Using Eq.~(126) discussed in \cite{bra20}:
\begin{eqnarray}
\sum^{k-1}_{i=1} \frac{H_i}{i^2} + \sum^{k-1}_{i=1} \frac{H_i^{(2)}}{i} = \zeta(3) - \zeta(2)\frac{1}{k} + \zeta(2)\frac{H_k}{k} - \sum^{\infty}_{i=1}
\frac{H_i}{(i+k)^2}
\end{eqnarray}
we get 
\begin{eqnarray}
\sum^{\infty}_{k=1} \frac{H_k H_k^{(2)} h_k}{k(2k-1)}  &=& \zeta(2)\zeta(3) - \zeta(2)\sum^{\infty}_{k=1} \frac{h_k}{k^2(2k-1)} + \zeta(2)\sum^{\infty}_{k=1} \frac{H_k h_k}{k(2k-1)} +
\sum^{\infty}_{k=1} \frac{H_k h_k}{k^3(2k-1)} + \sum^{\infty}_{k=1} \frac{H_k^{(2)} h_k}{k^2(2k-1)} \nonumber \\ &-& \sum^{\infty}_{k=1} \frac{H_k^{(3)} h_k}{k(2k-1)} -
\sum^{\infty}_{k=1} \frac{h_k}{k(2k-1)} \left( \sum^{\infty}_{i=1} \frac{H_i}{(i+k)^2} \right)~.
\end{eqnarray}
Only the last double sum on the right side of Eq.~(94) is not known explicitly and must be computed here. After rearranging the summations and after a quite long algebraic procedure
we arrive at the following result:
\begin{eqnarray}
\sum^{\infty}_{k=1} \frac{H_k H_k^{(2)} h_k}{k(2k-1)}  &=& \zeta(2)\zeta(3) - \zeta(2)\sum^{\infty}_{k=1} \frac{h_k}{k^2(2k-1)} + \zeta(2)\sum^{\infty}_{k=1} \frac{H_k h_k}{k(2k-1)} - 
\sum^{\infty}_{k=1} \frac{H_k^{(3)} h_k}{k(2k-1)} + \sum^{\infty}_{k=1} \frac{H_k^{(2)} h_k}{k^2(2k-1)}  \nonumber \\ &+& \sum^{\infty}_{k=1} \frac{H_k h_k}{k^3(2k-1)} -
8ln(2)\sum^{\infty}_{k=1} \frac{H_k h_k}{(2k+1)^2} -4\zeta(2)\sum^{\infty}_{k=1} \frac{H_k}{(2k+1)^2} + 2ln(2)\sum^{\infty}_{k=1} \frac{H_k h_k}{k^2} + \sum^{\infty}_{k=1} \frac{H_k^2 h_k}{k^2}
\nonumber \\ &-& \frac{7}{4}\zeta(3)\sum^{\infty}_{k=1} \frac{H_k}{k(2k+1)} - \left( \zeta(2) - 4 \right) \sum^{\infty}_{k=1} \frac{H_k h_k}{k(2k+1)} - 2 \sum^{\infty}_{k=1} \frac{H_k^2 }{k(2k+1)}
- \left( 4 - 4ln(2) \right) \sum^{\infty}_{k=1} \frac{H_k h_k^{(2)}}{k(2k+1)}  \nonumber \\ &+& \sum^{\infty}_{k=1} \frac{H_k H_k^{(2)} h_k}{k(2k+1)} - 4\sum^{\infty}_{k=1} \frac{H_k}{(2k+1)^2} 
\left( \sum^{k}_{i=1} \frac{H_{i-1}}{2i-1} \right) - \sum^{\infty}_{k=1} \frac{H_k}{k^2} \left( \sum^{k}_{i=1} \frac{h_i}{i} \right) - \sum^{\infty}_{k=1} \frac{H_k}{k(2k+1)}
\left( \sum^{k}_{i=1} \frac{h_i}{i^2} \right) \nonumber \\ &+& 2\sum^{\infty}_{k=1} \frac{H_k}{k(2k+1)} \left( \sum^{k}_{i=1} \frac{H_i}{(2i-1)^2} \right)~.
\end{eqnarray}
Besides the last four double sums on the right side of Eq.~(95) all Euler sums are known explicitly. The sum next to the last can be calculated by use of Eq.~(75) and the last sum results from Eq.~(83).
Therefore it remains to compute the other two double sums in order to arrive at Eq.~(89).
Starting with:
\begin{eqnarray}
\sum^{\infty}_{k=1} \frac{h_k}{k} \left( \sum^{\infty}_{i=1} \frac{H_i}{i(k-i)}\right) &=& \sum^{\infty}_{k=1} \frac{h_k}{k} \left( \sum^{\infty}_{i=1} \frac{H_i}{i(k+i)}\right) - \zeta(2)\sum^{\infty}_{k=1} \frac{h_k}{k^2}
- 2\sum^{\infty}_{k=1} \frac{H_k^{(2)} h_k}{k^2}
\end{eqnarray}
it follows first by rearranging the summations:
\begin{eqnarray}
-\sum^{\infty}_{i=1} \frac{H_i}{i} \left( \sum^{\infty}_{k=1} \frac{h_k}{k(i-k)}\right) &=& \sum^{\infty}_{i=1} \frac{H_i}{i} \left( \sum^{\infty}_{k=1} \frac{h_k}{k(i+k)}\right) - \zeta(2)\sum^{\infty}_{k=1} \frac{h_k}{k^2}
- 2\sum^{\infty}_{k=1} \frac{H_k^{(2)} h_k}{k^2}
\end{eqnarray}
and from this we get 
\begin{eqnarray}
2\sum^{\infty}_{i=1} \frac{H_i}{i} \left( \sum^{\infty}_{k=1} \frac{h_k}{k(i+k)}\right) &=& \zeta(2)\sum^{\infty}_{k=1} \frac{h_k}{k^2} + 2\sum^{\infty}_{k=1} \frac{H_k^{(2)} h_k}{k^2} + 
\sum^{\infty}_{k=1} \frac{H_k h_k}{k^3} + 2\sum^{\infty}_{k=1} \frac{H_k h_k^{(2)}}{k^2}
\end{eqnarray}
or
\begin{eqnarray}
\sum^{\infty}_{i=1} \frac{H_i}{i} \left( 2ln(2)\frac{h_i}{i} + \frac{H_i h_i}{i} - \frac{1}{i}\sum^{i}_{k=1} \frac{h_k}{k} \right) &=&  \frac{1}{2}\zeta(2)\sum^{\infty}_{k=1} \frac{h_k}{k^2} + 
\sum^{\infty}_{k=1} \frac{H_k^{(2)} h_k}{k^2} + \frac{1}{2} \sum^{\infty}_{k=1} \frac{H_k h_k}{k^3} + \sum^{\infty}_{k=1} \frac{H_k h_k^{(2)}}{k^2}
\end{eqnarray}
or
\begin{eqnarray}
\sum^{\infty}_{i=1} \frac{H_i}{i^2} \left( \sum^{i}_{k=1} \frac{h_k}{k} \right) &=&  2ln(2)\sum^{\infty}_{i=1} \frac{H_i h_i}{i^2} + \sum^{\infty}_{i=1} \frac{H_i^2 h_i}{i^2} -
\frac{1}{2}\zeta(2)\sum^{\infty}_{i=1} \frac{h_i}{i^2} - \sum^{\infty}_{i=1} \frac{H_i^{(2)} h_i}{k^i} - \frac{1}{2} \sum^{\infty}_{i=1} \frac{H_i h_i}{i^3} - \sum^{\infty}_{i=1} \frac{H_i h_i^{(2)}}{i^2}~.
\end{eqnarray}
From this it follows
\begin{eqnarray}
\sum^{\infty}_{k=1} \frac{H_k}{k^2} \left( \sum^{k}_{i=1} \frac{h_i}{i} \right) &=& \frac{403}{16}\zeta(5) - ln(2)\zeta(4) - \frac{7}{2}\left(ln(2)\right)^2\zeta(3) + 2\left(ln(2)\right)^3\zeta(2) - 
\frac{7}{15}\left(ln(2)\right)^5 \nonumber \\ &-&  16ln(2)Li_4\left(\frac{1}{2}\right) - 24Li_5\left(\frac{1}{2}\right)~.
\end{eqnarray}
With the help of Eq.~(38) we are able to apply a similar procedure to the remaining double sum. The result is: 
\begin{eqnarray}
\sum^{\infty}_{k=1} \frac{H_k}{(2k+1)^2} \left( \sum^{k}_{i=1} \frac{H_{i-1}}{2i-1} \right) &=& \frac{465}{32}\zeta(5) - \frac{19}{2}ln(2)\zeta(4) + \frac{7}{16}\zeta(2)\zeta(3) - \frac{7}{4}\left(ln(2)\right)^2\zeta(3)
+ \frac{5}{3}\left(ln(2)\right)^3\zeta(2) + \frac{2}{15}\left(ln(2)\right)^5  \nonumber \\ &-& 16Li_5\left(\frac{1}{2}\right)~. 
\end{eqnarray}
This way Eq.~(89) follows.

The third member is given by the following expression: 
\begin{eqnarray}
\sum^{\infty}_{k=1} \frac{H_k^{(2)} h_k^2}{k(2k-1)} &=& \frac{403}{16}\zeta(5) - \frac{91}{4}ln(2)\zeta(4) - \frac{45}{8}\zeta(4) + \frac{5}{4}\zeta(2)\zeta(3)
+ 14\left(ln(2)\right)^2\zeta(3) + \frac{3}{2}\zeta(3) - \frac{8}{3}\left(ln(2)\right)^3\zeta(2) \nonumber \\ &+& 6ln(2)\zeta(2) + \frac{4}{15}\left(ln(2)\right)^5 - 
2ln(2) \sum^{\infty}_{k=1} \frac{h_k}{k^3} - 32Li_5\left(\frac{1}{2}\right)~.
\end{eqnarray}
This cubic Euler sum results from the ansatz:
\begin{eqnarray}
\sum^{\infty}_{i=1} \frac{H_i^{(2)}}{2i-1} \left( \sum^{\infty}_{k=1} \frac{H_k}{(2k-1)(2i+2k-1)} \right) &=& 2ln(2)\sum^{\infty}_{i=1} \frac{H_i^{(2)}}{(2i-1)^2} - ln(2)\sum^{\infty}_{i=1} \frac{H_i^{(2)} h_i}{i(2i-1)} 
- \sum^{\infty}_{i=1} \frac{H_i^{(2)} h_i}{i(2i-1)^2} \nonumber \\ &+& \frac{1}{2}\sum^{\infty}_{i=1} \frac{H_i^{(2)} h_i^{(2)}}{i(2i-1)} + \frac{1}{2}\sum^{\infty}_{i=1} \frac{H_i^{(2)} h_i^2}{i(2i-1)}~.
\end{eqnarray}
Rearranging the summations we get:
\begin{eqnarray}
\sum^{\infty}_{k=1} \frac{H_k}{2k-1} \left( \sum^{\infty}_{k=1} \frac{H_i^{(2)}}{(2i-1)(2k+2i-1)} \right) &=& \left( 4ln(2) -\zeta(2) \right)\sum^{\infty}_{i=1} \frac{H_k}{(2k-1)^2} - 
4ln(2)\sum^{\infty}_{i=1} \frac{H_k}{(2k-1)^3} +  \frac{1}{2}\zeta(2)\sum^{\infty}_{i=1} \frac{H_k h_k}{k(2k-1)}  \nonumber \\ &+& 2ln(2)\sum^{\infty}_{i=1} \frac{H_k h_k^{(2)}}{k(2k-1)} - 
2\sum^{\infty}_{k=1} \frac{H_k}{k(2k-1)} \left( \sum^{k-1}_{i=1} \frac{h_i}{(2i-1)^2} \right)~.
\end{eqnarray}
Comparing both sides of Eq.~(104) and Eq.~(105) it shows up that only the last double sum is not known explicitly. It follows first:
\begin{eqnarray}
\sum^{\infty}_{k=1} \frac{H_k}{k(2k-1)} \left( \sum^{k-1}_{i=1} \frac{h_i}{(2i-1)^2} \right) &=& 2\sum^{\infty}_{k=1} \frac{H_k}{k(2k-1)} \left( \sum^{\infty}_{i=1} \frac{h_i}{(2i-1)^2} \right) -
\sum^{\infty}_{k=1} \frac{h_k}{(2k-1)^2}  \left( \sum^{k}_{i=1} \frac{H_i}{i(2i-1)} \right)~. 
\end{eqnarray}
By partial fraction decomposition of the second term it follows:
\begin{eqnarray}
\sum^{\infty}_{k=1} \frac{h_k}{(2k-1)^2} \left( \sum^{k}_{i=1} \frac{H_i}{i(2i-1)} \right) &=& 4\sum^{\infty}_{k=1} \frac{h_k^2}{(2k-1)^2} -2\sum^{\infty}_{k=1} \frac{H_k h_k}{(2k-1)^2} -
 \frac{1}{2}\sum^{\infty}_{k=1} \frac{H_k^{(2)} h_k}{(2k-1)^2} - \frac{1}{2}\sum^{\infty}_{k=1} \frac{H_k^2 h_k}{(2k-1)^2} \nonumber \\ &+& 2\sum^{\infty}_{k=1} \frac{H_k h_k^2}{(2k-1)^2} -
2\sum^{\infty}_{k=1} \frac{h_k}{(2k-1)^2}  \left( \sum^{k}_{i=1} \frac{h_i}{i} \right)~.
\end{eqnarray}
Again, only the last double sum on the right side of Eq.~(107) is unknown. It follows:
\begin{eqnarray}
\sum^{\infty}_{k=1} \frac{h_k}{(2k-1)} \left( \sum^{\infty}_{i=1} \frac{h_i}{i(2i+2k-1)} \right) &=& \frac{1}{2}\sum^{\infty}_{k=1} \frac{h_k}{(2k-1)^2} - \sum^{\infty}_{k=1} \frac{h_k^2}{k(2k-1)^2} +
\sum^{\infty}_{k=1} \frac{h_k}{(2k-1)^2} \left( \sum^{k}_{i=1} \frac{h_i}{i} \right)
\end{eqnarray}
or
\begin{eqnarray}
\sum^{\infty}_{i=1} \frac{h_i}{i} \left( \sum^{\infty}_{k=1} \frac{h_k}{(2k-1)(2i+2k-1)} \right) &=& \frac{3}{8}\zeta(2)\sum^{\infty}_{i=1} \frac{h_i}{i^2} - \frac{1}{4}\sum^{\infty}_{k=1} \frac{h_i^2}{i^3} +
\frac{1}{4}\sum^{\infty}_{i=1} \frac{h_i}{i^2} \left( \sum^{i}_{k=1} \frac{h_k}{k} \right)~.
\end{eqnarray}
As the double sum in Eq.~(109) is known explicitly (Eq.~(79)) the corresponding double sum follows by comparing  Eq.~(108) and Eq.~(109):
\begin{eqnarray}
\sum^{\infty}_{k=1} \frac{h_k}{(2k-1)^2} \left( \sum^{k}_{i=1} \frac{h_i}{i} \right) &=& \frac{31}{128}\zeta(5) + \frac{15}{32}ln(2)\zeta(4) + \frac{15}{16}\zeta(4) + \frac{7}{32}\zeta(2)\zeta(3)
+ \frac{7}{4}ln(2)\zeta(3) - \frac{3}{8}\zeta(3) + \frac{3}{2}\left(ln(2)\right)^2\zeta(2) \nonumber \\ &-& \frac{3}{2}ln(2)\zeta(2) - \frac{1}{4} \sum^{\infty}_{k=1} \frac{h_k}{k^3}~.
\end{eqnarray}
With this Eq.~(107) follows and finally Eq.~(103) results.

The last cubic Euler sum belonging to the seventh family follows to:
\begin{eqnarray}
\sum^{\infty}_{k=1} \frac{H_k h_k h_k^{(2)}}{k(2k-1)} &=& \frac{403}{64}\zeta(5) - \frac{211}{32}ln(2)\zeta(4) + \frac{105}{32}\zeta(4) - \frac{1}{16}\zeta(2)\zeta(3)
- \frac{7}{4}\left(ln(2)\right)^2\zeta(3) + \frac{7}{4}ln(2)\zeta(3) + \frac{7}{6}\left(ln(2)\right)^3\zeta(2) \nonumber \\ &-& \frac{3}{2}\left(ln(2)\right)^2\zeta(2) + \frac{1}{30}\left(ln(2)\right)^5 + 
\frac{1}{4}ln(2) \sum^{\infty}_{k=1} \frac{h_k}{k^3} - \frac{1}{4}\sum^{\infty}_{k=1} \frac{h_k}{k^3} - 4Li_5\left(\frac{1}{2}\right)~.
\end{eqnarray}
In order to proof this result we start with
\begin{eqnarray}
\sum^{\infty}_{k=1} \frac{h_k^{(2)}}{2k-1} \left( \sum^{\infty}_{i=1} \frac{h_i}{i(k+i)} \right) &=& 2ln(2)\sum^{\infty}_{k=1} \frac{h_k h_k^{(2)}}{k(2k-1)} + \sum^{\infty}_{k=1} \frac{H_k h_k h_k^{(2)}}{k(2k-1)}
- \sum^{\infty}_{k=1} \frac{h_k^{(2)}}{k(2k-1)} \left( \sum^{k}_{i=1} \frac{h_i}{i} \right)
\end{eqnarray}
or
\begin{eqnarray}
\sum^{\infty}_{k=1} \frac{h_k^{(2)}}{2k-1} \left( \sum^{\infty}_{i=1} \frac{h_i}{i(k+i)} \right) &=& 2ln(2)\sum^{\infty}_{k=1} \frac{h_k h_k^{(2)}}{k(2k-1)} +
\sum^{\infty}_{k=1} \frac{h_k^{(2)}}{k(2k-1)} \left( \sum^{k}_{i=1} \frac{H_{i-1}}{2i-1} \right)
\end{eqnarray}
or
\begin{eqnarray}
\sum^{\infty}_{i=1} \frac{h_i}{i} \left( \sum^{\infty}_{k=1} \frac{h_k^{(2)}}{(2k-1)(i+k)} \right) &=& 2ln(2)\sum^{\infty}_{k=1} \frac{h_k h_k^{(2)}}{k(2k-1)} -
\sum^{\infty}_{k=1} \frac{h_k^{(2)}}{k(2k-1)} \left( \sum^{k}_{i=1} \frac{1}{i(2i-1)} \right) \nonumber \\ &+& \frac{1}{2} \sum^{\infty}_{k=1} \frac{h_k^{(2)}}{k(2k-1)} \left( \sum^{k}_{i=1} \frac{H_{i}}{i(2i-1)} \right) +
\frac{1}{2} \sum^{\infty}_{k=1} \frac{h_k^{(2)}}{k(2k-1)} \left( \sum^{k}_{i=1} \frac{H_{i}}{i} \right)~.
\end{eqnarray}
Next we rearrange the summations in the last to the next double sum on the right side of Eq.~(114):
\begin{eqnarray}
\sum^{\infty}_{k=1} \frac{h_k^{(2)}}{k(2k-1)} \left( \sum^{k}_{i=1} \frac{H_{i}}{i(2i-1)} \right)  &=& \sum^{\infty}_{k=1} \frac{h_k^{(2)}}{k(2k-1)} \left( \sum^{\infty}_{i=1} \frac{H_{i}}{i(2i-1)} \right)
+ \sum^{\infty}_{k=1} \frac{H_k h_k^{(2)}}{k^2(2k-1)^2} \nonumber \\ &-& \sum^{\infty}_{k=1} \frac{H_{k}}{k(2k-1)} \left( 2\sum^{k}_{i=1} \frac{h_i^{(2)}}{2i-1} - \sum^{k}_{i=1} \frac{h_i^{(2)}}{i} \right)~.
\end{eqnarray}
Inserting now the two well known identities \cite{ade16}:
\begin{eqnarray}
\sum^{k}_{i=1} \frac{h_i^{(2)}}{2i-1} = h_k^{(3)} +  h_k h_k^{(2)} - \sum^{k}_{i=1} \frac{h_i}{(2i-1)^2}
\end{eqnarray}
and
\begin{eqnarray}
\sum^{k}_{i=1} \frac{h_i^{(2)}}{i} =  H_k h_k^{(2)} - \sum^{k}_{i=1} \frac{H_{i-1}}{(2i-1)^2}
\end{eqnarray}
in Eq.~(115) we arrive at the following expression for Eq.~(114):
\begin{eqnarray}
\sum^{\infty}_{k=1} \frac{h_k}{k}\hspace{-0.4cm} && \left( \sum^{\infty}_{i=1} \frac{h_i^{(2)}}{(2i-1)(k+i)} \right) = 2ln(2)\sum^{\infty}_{k=1} \frac{h_k h_k^{(2)}}{k(2k-1)} - 2\sum^{\infty}_{k=1} \frac{h_k h_k^{(2)}}{k(2k-1)}
+ \sum^{\infty}_{k=1} \frac{H_k h_k^{(2)}}{k(2k-1)} + \frac{1}{4}\sum^{\infty}_{k=1} \frac{H_k^{(2)} h_k^{(2)}}{k(2k-1)} \nonumber \\ &+& \frac{3}{4}\sum^{\infty}_{k=1} \frac{H_k^2 h_k^{(2)}}{k(2k-1)} +
\frac{1}{2} \sum^{\infty}_{k=1} \frac{h_k^{(2)}}{k(2k-1)} \left( \sum^{\infty}_{i=1} \frac{H_{i}}{i(2i-1)} \right) + \frac{1}{2}\sum^{\infty}_{k=1} \frac{H_k h_k^{(2)}}{k^2(2k-1)^2}  - 
\sum^{\infty}_{k=1} \frac{H_k h_k^{(3)}}{k(2k-1)} \nonumber \\ &-& \sum^{\infty}_{k=1} \frac{H_k h_k h_k^{(2)}}{k(2k-1)} - \frac{1}{2} \sum^{\infty}_{k=1} \frac{H_k}{k(2k-1)} \left( \sum^{k}_{i=1} \frac{H_{i-1}}{(2i-1)^2} \right)
+ \sum^{\infty}_{k=1} \frac{H_k}{k(2k-1)} \left( \sum^{k}_{i=1} \frac{h_{i}}{(2i-1)^2} \right)~.
\end{eqnarray}
On the other hand we have
\begin{eqnarray}
\sum^{\infty}_{k=1} \frac{h_k}{k}\hspace{-0.4cm} && \left( \sum^{\infty}_{i=1} \frac{h_i^{(2)}}{(2i-1)(k+i)} \right) = \left( \frac{21}{8}\zeta(3) - 
\frac{3}{2}ln(2)\zeta(2) \right) \sum^{\infty}_{k=1} \frac{h_k}{k(2k+1)} + \frac{3}{2}\zeta(2)\sum^{\infty}_{k=1} \frac{h_k^2}{k(2k+1)} \nonumber \\ &-& 2ln(2)\sum^{\infty}_{k=1} \frac{h_k h_k^{(2)}}{k(2k+1)} -
\sum^{\infty}_{k=1} \frac{h_k}{k(2k+1)} \left( \sum^{\infty}_{i=1} \frac{H_{i-1}}{(2i-1)^2} \right)~.
\end{eqnarray}
This way Eq.~(111) follows by comparison of Eq.~(118) and Eq.~(119) if all corresponding double sums are known explicitly. The last double sum on the right side of Eq.~(118) results from Eq.~(106) by
rearranging the summations in Eq.~(106). The calculational scheme for the last double sum on the right side of Eq.~(119) appears analogously to the explicit computation of Eq.~(83) where $H_k$ has to be interchanged
in the outer sum by $h_k$. The final result is:
\begin{eqnarray}
\sum^{\infty}_{k=1} \frac{h_k}{k(2k+1)} \left( \sum^{\infty}_{i=1} \frac{H_{i-1}}{(2i-1)^2} \right) &=& -\frac{31}{4}\zeta(5) + \frac{211}{32}ln(2)\zeta(4) + \frac{21}{16}\zeta(2)\zeta(3) - 
\frac{7}{4}\left(ln(2)\right)^2\zeta(3) - \frac{7}{6}\left(ln(2)\right)^3\zeta(2) \nonumber \\ &-& \frac{1}{30}\left(ln(2)\right)^5 + \frac{1}{4}ln(2) \sum^{\infty}_{k=1} \frac{h_k}{k^3} + Li_5\left(\frac{1}{2}\right)~.
\end{eqnarray}
This way it remains to calculate the first double sum on the right side of Eq.~(118). We start with Eq.~(88)
\begin{eqnarray}
\sum^{\infty}_{k=1} \frac{H_{k-1}}{(2k-1)^2} \left( \sum^{k}_{i=1} \frac{H_i}{2i+1} \right) &=& \sum^{\infty}_{k=1} \frac{H_{k-1}}{(2k-1)^2} \left( \sum^{k}_{i=1} \frac{H_i}{2i-1} \right) + 
\sum^{\infty}_{k=1} \frac{H_k H_{k-1}}{(2k+1)(2k-1)^2} \nonumber \\ &-& \sum^{\infty}_{k=1} \frac{H_{k-1}}{(2k-1)^2} \left( \sum^{k}_{i=1} \frac{1}{i(2i-1)} \right)~.
\end{eqnarray}
From this it follows: 
\begin{eqnarray}
\frac{1}{2}\sum^{\infty}_{k=1} \frac{H_{k-1}}{(2k-1)^2} \left( \sum^{k}_{i=1} \frac{H_i}{i(2i-1)} \right) &=& - \frac{1}{2}\sum^{\infty}_{k=1} \frac{H_{k-1}}{(2k-1)^2} \left( \sum^{k}_{i=1} \frac{H_i}{i} \right) +
\sum^{\infty}_{k=1} \frac{H_{k-1}}{(2k-1)^2} \left( \sum^{k}_{i=1} \frac{H_i}{2i+1} \right) \nonumber \\ &+& \sum^{\infty}_{k=1} \frac{H_{k-1}}{(2k-1)^2} \left( \sum^{k}_{i=1} \frac{1}{i(2i-1)} \right) -
\sum^{\infty}_{k=1} \frac{H_k H_{k-1}}{(2k+1)(2k-1)^2}~.
\end{eqnarray}
Rearranging the summations on the left side of Eq.~(122) the corresponding double sum follows explicitly as all other terms are known.
\begin{eqnarray}
\sum^{\infty}_{k=1} \frac{H_k}{k(2k-1)} \left( \sum^{\infty}_{i=1} \frac{H_{i-1}}{(2i-1)^2} \right) &=& -\frac{465}{8}\zeta(5) + 38ln(2)\zeta(4) + \frac{45}{8}\zeta(4) + \frac{7}{8}\zeta(2)\zeta(3) - 
\frac{7}{2}\left(ln(2)\right)^2\zeta(3) - \frac{35}{2}ln(2)\zeta(3) \nonumber \\ &-& \frac{2}{3}\left(ln(2)\right)^3\zeta(2) + 6\left(ln(2)\right)^2\zeta(2) - \frac{8}{15}\left(ln(2)\right)^5 + 
3\sum^{\infty}_{k=1} \frac{h_k}{k^3} + 64Li_5\left(\frac{1}{2}\right)~.
\end{eqnarray}
This completes the computation of Eq.~(111).

\subsection{Eighth family}
The eighth family consists of only one member which is well known from literature \cite{zeh07}:
\begin{eqnarray}
\sum^{\infty}_{k=1} \frac{h_k^3}{k^2} &=& \frac{21}{8}\zeta(2)\zeta(3)~.
\end{eqnarray}
This way all cubic Euler sums of order five which are members of the eight different families have been calculated explicitly in terms of zeta values and polylogarithmic values $Li_4(1/2)$ and $Li_5(1/2)$.

\section{The order six case for cubic Euler sums}
In this section we demonstrate that all members of the eight families of cubic Euler sums of degree six can be calculated in terms of zeta values and polylogarithmic values $Li_4(1/2)$,
$Li_5(1/2)$, $Li_6(1/2)$, $Li_6(-1/2)$ and $Li_6(-1/8)$. We start with the following identity \cite{ce20}:
\begin{eqnarray}
\sum^{\infty}_{k=1} (-)^{k+1}\frac{H_k^3}{k^3} &=& \frac{9}{2}\zeta({\overline 5},1) - 12Li_6\left(\frac{1}{2}\right) - 12ln(2)Li_5\left(\frac{1}{2}\right) + 3\zeta(2)Li_4\left(\frac{1}{2}\right)
- 6\left(ln(2)\right)^2Li_4\left(\frac{1}{2}\right) + \frac{771}{64}\zeta(6)  \nonumber \\ &-& \frac{15}{8}\left(ln(2)\right)^2\zeta(4) 
- \frac{207}{64}\zeta(3)^2 + \frac{21}{8}ln(2)\zeta(2)\zeta(3) - \frac{7}{4}\left(ln(2)\right)^3\zeta(3) + \frac{7}{8}\left(ln(2)\right)^4\zeta(2) \nonumber \\ &-&  \frac{1}{6}\left(ln(2)\right)^6~.
\end{eqnarray}
With the corresponding multiple zeta value:
\begin{eqnarray}
\zeta({\overline 5},1) &=& -\frac{93}{128}\zeta(6) + \frac{31}{64}\zeta(3)^2 + \frac{1}{16}\sum^{\infty}_{k=1}\frac{h_k}{k^5}
\end{eqnarray}
it follows
\begin{eqnarray}
\sum^{\infty}_{k=1} (-)^{k+1}\frac{H_k^3}{k^3} &=& \frac{2247}{256}\zeta(6) - \frac{15}{8}\left(ln(2)\right)^2\zeta(4) - \frac{135}{128}\zeta(3)^2 + \frac{21}{8}ln(2)\zeta(2)\zeta(3) -
\frac{7}{4}\left(ln(2)\right)^3\zeta(3) + \frac{7}{8}\left(ln(2)\right)^4\zeta(2) \nonumber \\ &-& \frac{1}{6}\left(ln(2)\right)^6 - 12Li_6\left(\frac{1}{2}\right) - 12ln(2)Li_5\left(\frac{1}{2}\right) +
3\zeta(2)Li_4\left(\frac{1}{2}\right) - 6\left(ln(2)\right)^2Li_4\left(\frac{1}{2}\right) \nonumber \\ &+& \frac{9}{32}\sum^{\infty}_{k=1}\frac{h_k}{k^5}~.
\end{eqnarray}
On the other hand we have:
\begin{eqnarray}
\sum^{\infty}_{k=1} (-)^{k+1}\frac{H_k^3}{k^3} &=& \frac{31}{32} \sum^{\infty}_{k=1} \frac{H_k^3}{k^3} - \frac{1}{4} \sum^{\infty}_{k=1} \frac{h_k^3}{k^3} -
\frac{3}{4} \sum^{\infty}_{k=1} \frac{H_k h_k^2}{k^3} - \frac{3}{16} \sum^{\infty}_{k=1} \frac{H_k^2 h_k}{k^3}
\end{eqnarray}
and similar for the corresponding cubic Euler sums with odd-type denominators. This provides us with to basic identities between different members of the corresponding eight families  
of cubic Euler sums of order six. The next step consists in the formulation of new two-valued help functions. It results first:
\begin{eqnarray}
F(k) = \sum^{\infty}_{i=1} \frac{H_i^2}{(2i-1)(2i+2k-1)} &=& \frac{1}{k}\sum^{\infty}_{i=1} \frac{H_i^2}{(2i-1)(2i+1)} + \left( 2\left( ln(2) \right)^2 - \frac{1}{2}\zeta(2) \right) \frac{h_{k-1}}{k}
+ 2\frac{h_{k-1} h_{k-1}^{(2)}}{k} + \frac{4}{3}\frac{h_{k-1}^{(3)}}{k} \nonumber \\ &-& 2ln(2)\frac{h_{k-1}^2}{k} - \frac{2}{k}\sum^{k-1}_{i=1}\frac{h_i}{(2i-1)^2} + \frac{h_{k-1}^3}{k}~. 
\end{eqnarray}
This two-valued nonlinear help function represents the solution of the corresponding inhomogeneous first order difference equation.
With this we get:
\begin{eqnarray}  
\sum^{\infty}_{k=1} \frac{1}{k^2} \left( \sum^{\infty}_{i=1} \frac{H_i^2}{(2i-1)(2i+2k-1)} \right) &=& \sum^{\infty}_{i=1} \frac{H_i^2}{2i-1} \left( \sum^{\infty}_{k=1} \frac{1}{k^2(2i+2k-1)} \right)
\nonumber \\ &=& \zeta(2) \sum^{\infty}_{i=1}\frac{H_i^2}{(2i-1)^2} - 2\left( \sum^{\infty}_{i=1}\frac{H_i^2}{(2i-1)^2} \left[ \frac{2h_i}{2i-1} - \frac{2ln(2)}{2i-1} \right] \right) 
\nonumber \\ &=& \zeta(2)\sum^{\infty}_{i=1}\frac{H_i^2}{(2i-1)^2} + 4ln(2)\sum^{\infty}_{i=1}\frac{H_i^2}{(2i-1)^3} - 4\sum^{\infty}_{i=1}\frac{H_i^2 h_i}{(2i-1)^3}
\end{eqnarray}
or
\begin{eqnarray} 
\zeta(2)\sum^{\infty}_{i=1}\frac{H_i^2}{(2i-1)^2} &+& 4ln(2)\sum^{\infty}_{i=1}\frac{H_i^2}{(2i-1)^3} - 4\sum^{\infty}_{i=1}\frac{H_i^2 h_i}{(2i-1)^3} = \zeta(3) \left( \frac{1}{2}\zeta(2) +
2ln(2) - 2\left( ln(2) \right)^2 \right) \nonumber \\ &+& \left( 2\left( ln(2) \right)^2 - \frac{1}{2}\zeta(2) \right) \sum^{\infty}_{i=1}\frac{h_{i-1}}{i^3} + 2\sum^{\infty}_{i=1}
\frac{h_{i-1} h_{i-1}^{(2)}}{i^3} + \frac{4}{3} \sum^{\infty}_{i=1}\frac{h_{i-1}^{(3)}}{i^3} - 2ln(2)\sum^{\infty}_{i=1}\frac{h_{i-1}^2}{i^3} \nonumber \\ &+& \frac{2}{3} \sum^{\infty}_{i=1}\frac{h_{i-1}^3}{i^3}
-2\zeta(3)\sum^{\infty}_{i=1}\frac{h_i}{(2i-1)^2} - 2\sum^{\infty}_{i=1}\frac{h_i}{i^3(2i-1)^2} + 2\sum^{\infty}_{i=1}\frac{H_i^{(3)} h_i}{(2i-1)^2}~.
\end{eqnarray}
With this we have an identity between the two cubic Euler sums:
\begin{eqnarray}
\sum^{\infty}_{i=1}\frac{H_i^2 h_i}{(2i-1)^3} \sim \sum^{\infty}_{i=1}\frac{h_i^3}{i^3}
\end{eqnarray}
established as all other Euler sums in Eq.~(131) are known explicitly.
Analogously we were able to find 
\begin{eqnarray}
\sum^{\infty}_{i=1}\frac{H_i h_i^2}{(2i-1)^3} \sim \sum^{\infty}_{i=1}\frac{H_i h_i^2}{i^3}
\end{eqnarray}
by calculating the following expression
\begin{eqnarray}
\sum^{\infty}_{i=1}\frac{H_i}{(2i-1)^2} \left( \sum^{\infty}_{k=1} \frac{H_k}{k(2i+2k-1)} \right) &=& \frac{1}{2} \sum^{\infty}_{k=1} \frac{H_k}{k^2} \left( \sum^{\infty}_{i=1}\frac{H_i}{(2i-1)^2} \right) -
\frac{1}{2} \sum^{\infty}_{k=1} \frac{H_k}{k^2} \left( \sum^{\infty}_{i=1}\frac{H_i}{(2i-1)(2i+2k-1)} \right)
\end{eqnarray}
where both two valued help functions appearing in Eq.~(134) have been defined before in Eq.~(44) and Eq.~(45).
A third identity can be found performing the following computation:
\begin{eqnarray}
\sum^{\infty}_{i=1}\frac{H_i}{i^2} \left( \sum^{\infty}_{i=1}\frac{h_i}{k(k+i)} \right) &=& 2ln(2) \frac{H_k h_k}{k^3} + \frac{H_k^2 h_k}{k^3} - \sum^{\infty}_{k=1} \frac{H_k}{k^3} 
\left( \sum^{k}_{i=1}\frac{h_i}{i} \right)
\end{eqnarray}
or
\begin{eqnarray}
\sum^{\infty}_{k=1}\frac{h_k}{k}\left( \sum^{\infty}_{i=1} \frac{H_i}{i^2(k+i)} \right)  &=& \sum^{\infty}_{k=1}\frac{h_k}{k^2} \left( \sum^{\infty}_{i=1} \frac{H_i}{i^2} \right) -
\sum^{\infty}_{k=1}\frac{h_k}{k^2} \left( \sum^{\infty}_{i=1} \frac{H_i}{i(k+i)} \right)~.
\end{eqnarray}
From this we get:
\begin{eqnarray}
\sum^{\infty}_{k=1} \frac{H_k}{k^3} \left( \sum^{k}_{i=1}\frac{h_i}{i} \right)  &=& -\frac{7}{2}\zeta(3)^2 + 2ln(2) \frac{H_k h_k}{k^3} + \zeta(2)\sum^{\infty}_{k=1}\frac{h_k}{k^3} -
2\sum^{\infty}_{k=1} \frac{H_k h_k}{k^4} + \frac{1}{2} \sum^{\infty}_{k=1} \frac{H_k^{(2)} h_k}{k^3} + \frac{3}{2} \sum^{\infty}_{k=1} \frac{H_k^2 h_k}{k^3}~. 
\end{eqnarray}
Calculating now
\begin{eqnarray}
\sum^{\infty}_{i=1} \frac{h_i}{2i-1} \left( \sum^{\infty}_{k=1} \frac{H_k}{k^2(2k+2i-1)} \right) &=& \sum^{\infty}_{i=1} \frac{h_i}{(2i-1)^2} \left( \sum^{\infty}_{k=1}\frac{H_k}{k^2} \right) -
4\left(ln(2)\right)^2 \sum^{\infty}_{i=1} \frac{h_i}{(2i-1)^3} + 8\sum^{\infty}_{i=1} \frac{h_i^2}{(2i-1)^4} \nonumber \\ &+& 8ln(2)\sum^{\infty}_{i=1} \frac{h_i^2}{(2i-1)^3} - 4 \sum^{\infty}_{i=1}
\frac{h_i h_i^{(2)}}{(2i-1)^3} - 8ln(2)\sum^{\infty}_{i=1} \frac{h_i}{(2i-1)^4}  \nonumber \\ &-& 4\sum^{\infty}_{i=1} \frac{h_i^3}{(2i-1)^3}
\end{eqnarray}
or
\begin{eqnarray}
\sum^{\infty}_{k=1} \frac{H_k}{k^2} \left( \sum^{\infty}_{i=1} \frac{h_i}{(2i-1)(2k+2i-1)} \right) &=& \frac{3}{8}\zeta(2)\sum^{\infty}_{k=1} \frac{H_k}{k^3} - \frac{1}{4} \sum^{\infty}_{k=1}
\frac{H_k h_k}{k^4} + \frac{1}{4}\sum^{\infty}_{k=1} \frac{H_k}{k^3} \left( \sum^{k}_{i=1}\frac{h_i}{i} \right)~.
\end{eqnarray}
Comparing Eq.~(137) with Eq.~(139) a third identity
\begin{eqnarray}
\sum^{\infty}_{k=1} \frac{H_k^2 h_k}{k^3} \sim \sum^{\infty}_{k=1} \frac{h_k^3}{(2k-1)^3}
\end{eqnarray}
has been found.

\subsection{First family}
In order to calculate explicitly the first member of this family another identity \cite{ce20} is essential
\begin{eqnarray}
\sum^{\infty}_{k=1} \frac{(-1)^{k+1}}{k^2} \left( \sum^{k}_{i=1} \frac{(-1)^{i+1}}{i} \right)^2 \left( \sum^{k}_{i=1} \frac{(-1)^{i+1}}{i^2} \right) &=& 4\zeta({\overline 5},1) - 8Li_6\left(\frac{1}{2}\right) -
8ln(2)Li_5\left(\frac{1}{2}\right) - \zeta(2))Li_4\left(\frac{1}{2}\right) \nonumber \\ &-& 4\left(ln(2)\right)^2Li_4\left(\frac{1}{2}\right) + \frac{4495}{384}\zeta(6) + \frac{5}{8}\left(ln(2)\right)^2\zeta(4) -
\frac{43}{32}\zeta(3)^2 \nonumber \\ &-& \frac{7}{8}ln(2)\zeta(2)\zeta(3) - \frac{7}{6}\left(ln(2)\right)^3\zeta(3) + \frac{11}{24}\left(ln(2)\right)^4\zeta(2) \nonumber \\ &-& \frac{1}{9}\left(ln(2)\right)^6~.
\end{eqnarray}
With
\begin{eqnarray}
\sum^{2k}_{i=1} \frac{(-1)^{i+1}}{i} = h_k - \frac{1}{2}H_k
\end{eqnarray}
it follows:
\begin{eqnarray}
\sum^{\infty}_{k=1} \frac{(-1)^{k+1}}{k^2} \left( \sum^{k}_{i=1} \frac{(-1)^{i+1}}{i} \right)^2 \left( \sum^{k}_{i=1} \frac{(-1)^{i+1}}{i^2} \right) &=& \frac{1}{2}\sum^{\infty}_{k=1} \frac{h_k^2 h_k^{(2)}}{k^2}
- \frac{1}{8} \sum^{\infty}_{k=1} \frac{H_k^{(2)} h_k^2}{k^2} + \frac{1}{8} \sum^{\infty}_{k=1} \frac{H_k^2 h_k^{(2)}}{k^2} - \frac{1}{32} \sum^{\infty}_{k=1} \frac{H_k^2 H_k^{(2)}}{k^2} \nonumber \\ &-&
\frac{1}{2}\sum^{\infty}_{k=1} \frac{H_k h_k h_k^{(2)}}{k^2} + \frac{1}{8}\sum^{\infty}_{k=1} \frac{H_k H_k^{(2)} h_k}{k^2}~. 
\end{eqnarray}
As a next step we introduce the two-valued help function:
\begin{eqnarray}
\sum^{\infty}_{i=1} \frac{H_i^2}{(2i+2k-1)^2} &=& \sum^{\infty}_{i=1} \frac{H_i^2}{(2i+1)^2} + \left( 7\zeta(3) - 6ln(2)\zeta(2) \right) h_{k-1} + \left( \zeta(2) -4\left(ln(2)\right)^2 \right) h_{k-1}^{(2)}
- 6h_{k-1}^{(4)} \nonumber \\ &+& 8h_{k-1} h_{k-1}^{(2)} - 2 h_{k-1}^{(2)} h_{k-1}^{(2)} - 8h_{k-1}h_{k-1}^{(3)} + 3\zeta(2)h_{k-1}^2 + 8\sum^{k-1}_{i=1}\frac{h_i}{(2i-1)^3} - 4h_{k-1}^2h_{k-1}^{(2)} 
\end{eqnarray}
which again represents the solution of the corresponding inhomogeneous first order difference equation. With this we get:
\begin{eqnarray}
\sum^{\infty}_{k=1} \frac{1}{k^2} \left(\sum^{\infty}_{i=1} \frac{H_i^2}{(2i+2k-1)^2} \right) &=& \zeta(2)\sum^{\infty}_{i=1} \frac{H_i^2}{(2i+1)^2} + \left( 7\zeta(3) - 6ln(2)\zeta(2) \right) 
\sum^{\infty}_{k=1} \frac{h_{k-1}}{k^2} + \left( \zeta(2) -4\left(ln(2)\right)^2 \right) \sum^{\infty}_{k=1} \frac{h_{k-1}^{(2)}}{k^2} \nonumber \\ &-& 6\sum^{\infty}_{k=1}\frac{h_{k-1}^{(4)}}{k^2}
+ 8\sum^{\infty}_{k=1}\frac{h_{k-1} h_{k-1}^{(2)}}{k^2} - 2\sum^{\infty}_{k=1}\frac{h_{k-1}^{(2)} h_{k-1}^{(2)}}{k^2} - 8\sum^{\infty}_{k=1}\frac{h_{k-1}h_{k-1}^{(3)}}{k^2}
\nonumber \\ &+& 3\zeta(2)\sum^{\infty}_{k=1} \frac{h_{k-1}^2}{k^2} + 8\sum^{\infty}_{k=1} \frac{1}{k^2} \left( \sum^{k-1}_{i=1}\frac{h_i}{(2i-1)^3} \right) -4\sum^{\infty}_{k=1} \frac{h_{k-1}^2h_{k-1}^{(2)}}{k^2} 
\end{eqnarray}
or
\begin{eqnarray}
\sum^{\infty}_{i=1} H_i^2 \left( \sum^{\infty}_{k=1} \frac{1}{k^2(2i+2k-1)^2} \right) &=& 4 \zeta(2)\sum^{\infty}_{k=1} \frac{H_k^2}{(2k-1)^2} - 8\sum^{\infty}_{k=1} \frac{H_k^2 h_k}{(2k-1)^3} + 8ln(2) 
\frac{H_k^2}{(2k-1)^3} \nonumber \\ &-& 4\sum^{\infty}_{k=1} \frac{H_k^2 h_{k}^{(2)}}{(2k-1)^2}~.
\end{eqnarray}
The last cubic Euler sum on the right side can be calculated in a direct way. This is demonstrated by the explicit calculation of Eq.~(199). Using furthermore the identity defined in Eq.~(132) it follows after
a tedious computational procedure:
\begin{eqnarray}
\sum^{\infty}_{k=1} \frac{h_k^{(2)} h_k^2}{k^2} = \frac{315}{128}\zeta(6) + \frac{49}{48}\zeta(3)^2 - \frac{1}{3}\sum^{\infty}_{k=1} \frac{h_{k}^3}{k^3}~. 
\end{eqnarray}
With this we are able to calculate an identity that connects two cubic Euler sums appearing in Eq.~(128) and in Eq.~(143). Another helpful identity can be obtained from the following calculational procedure. We start with:
\begin{eqnarray}
\sum^{\infty}_{k=1} \frac{h_k}{k^2} \left( \sum^{\infty}_{i=1} \frac{h_i}{i(i+k)} \right) &=& \frac{49}{32}\zeta(3)^2 = 2ln(2)\sum^{\infty}_{k=1} \frac{h_k^2}{k^3} +
\sum^{\infty}_{k=1} \frac{H_k h_k^2}{k^3} - \sum^{\infty}_{k=1} \frac{h_k}{k^3} \left( \sum^{k}_{i=1} \frac{h_i}{i} \right)~.
\end{eqnarray} 
This way the corresponding nested double sum results to
\begin{eqnarray}
\sum^{\infty}_{k=1} \frac{h_k}{k^3} \left( \sum^{k}_{i=1} \frac{h_i}{i} \right) = -\frac{49}{32}\zeta(3)^2 + 2ln(2)\sum^{\infty}_{k=1} \frac{h_k^2}{k^3} + \sum^{\infty}_{k=1} \frac{H_k h_k^2}{k^3}~.
\end{eqnarray}
We continue with
\begin{eqnarray}
\sum^{\infty}_{k=1} \frac{h_k}{k} \left( \sum^{\infty}_{i=1} \frac{h_i}{(2i-1)^2(2i+2k-1)} \right) &=& \frac{1}{2}\zeta(2) \sum^{\infty}_{i=1} \frac{h_i}{(2i-1)^3} + 
\sum^{\infty}_{i=1} \frac{h_i}{(2i-1)^3} \left( \sum^{i-1}_{k=1} \frac{h_k}{k} \right)  \nonumber  
\end{eqnarray}
or
\begin{eqnarray}
\sum^{\infty}_{i=1} \frac{h_i}{(2i-1)^2} \left( \sum^{\infty}_{k=1} \frac{h_k}{k(2i+2k-1)} \right) &=& \frac{1}{2} \sum^{\infty}_{i=1}\frac{h_k}{k^2} \left( \sum^{\infty}_{i=1} \frac{h_i}{(2i-1)^2} \right) -
\frac{3}{16}\zeta(2)\sum^{\infty}_{k=1} \frac{h_k}{k^3} - \frac{1}{8} \sum^{\infty}_{k=1} \frac{h_k}{k^3} \left( \sum^{k-1}_{i=1} \frac{h_i}{i} \right)~.
\end{eqnarray}
This way it follows by using Eq.~(149):
\begin{eqnarray}
\sum^{\infty}_{k=1} \frac{h_k}{(2k-1)^3} \left( \sum^{k}_{i=1} \frac{h_i}{i} \right) &=& \frac{49}{256}\zeta(3)^2 + \frac{1}{2} \sum^{\infty}_{k=1} \frac{h_k}{k^2} \left( \sum^{\infty}_{i=1} \frac{h_i}{(2i-1)^2} \right)
- \frac{3}{16}\zeta(2)\sum^{\infty}_{k=1} \frac{h_k}{k^3} + \frac{1}{8} \sum^{\infty}_{k=1} \frac{h_k^2}{k^4} \nonumber \\ &-& \frac{1}{2}\zeta(2) \sum^{\infty}_{k=1} \frac{h_k}{(2k-1)^3} +
\sum^{\infty}_{k=1} \frac{h_k^2}{k(2k-1)^3} - \frac{1}{4}ln(2) \sum^{\infty}_{k=1} \frac{h_k^2}{k^3} -  \frac{1}{8}\sum^{\infty}_{k=1} \frac{H_k h_k^2}{k^3}~.
\end{eqnarray}
Now we proceed with \cite{ade16}
\begin{eqnarray}
\sum^{\infty}_{k=1} \frac{h_k}{(2k-1)^3} \left( \sum^{k}_{i=1} \frac{H_{i-1}}{i} \right) &=& \sum^{\infty}_{k=1} \frac{H_k h_k^2}{(2k-1)^3} - \sum^{\infty}_{k=1} \frac{h_k}{(2k-1)^3} \left( \sum^{k}_{i=1} \frac{h_i}{i} \right)~.
\end{eqnarray}
Using explicitly Eq.~(133) it follows after some algebraic manipulations:
\begin{eqnarray}
\sum^{\infty}_{k=1} \frac{h_k}{(2k-1)^3} \left( \sum^{k}_{i=1} \frac{H_{i-1}}{2i-1} \right) &=& \frac{195}{256}\zeta(6) + \frac{465}{64}ln(2)\zeta(5) - \frac{121}{32}\left(ln(2)\right)^2\zeta(4)
- \frac{133}{256}\zeta(3)^2 - \frac{43}{32}ln(2)\zeta(2)\zeta(3) \nonumber \\ &-& \frac{1}{3}\left(ln(2)\right)^4\zeta(2) + \frac{1}{30}\left(ln(2)\right)^6 + \frac{1}{8}\zeta(2)\sum^{\infty}_{k=1}\frac{h_k}{k^3}
- \frac{1}{8}\sum^{\infty}_{k=1}\frac{h_k}{k^5} - 4ln(2)Li_5\left(\frac{1}{2}\right)~. 
\end{eqnarray}
In a next step we calculate
\begin{eqnarray}
\sum^{\infty}_{k=1} \frac{H_k}{k(2k-1)} \left( \sum^{k}_{i=1} \frac{h_i}{(2i-1)^3} \right) &=& \sum^{\infty}_{k=1} \frac{H_k}{k(2k-1)} \left( \sum^{\infty}_{i=1} \frac{h_i}{(2i-1)^3} \right) + 
\frac{H_k h_k}{k(2k-1)^4} \nonumber \\ &-& \sum^{\infty}_{k=1} \frac{h_k}{(2k-1)^3} \left( \sum^{k}_{i=1} \frac{H_{i}}{i(2i-1)} \right)~. 
\end{eqnarray}
From this it follows by explicit use of Eq.~(132):
\begin{eqnarray}
\sum^{\infty}_{k=1} \frac{H_k}{k(2k-1)} \left( \sum^{k}_{i=1} \frac{h_i}{(2i-1)^3} \right) &=& \frac{285}{512}\zeta(6) - \frac{1085}{64}ln(2)\zeta(5) + \frac{837}{64}\zeta(5) + \left(ln(2)\right)^2\zeta(4)
- \frac{61}{16}ln(2)\zeta(4) + \frac{203}{384}\zeta(3)^2 \nonumber \\ &+& \frac{77}{16}ln(2)\zeta(2)\zeta(3) - \frac{7}{4}\zeta(2)\zeta(3) + \frac{2}{3}\left(ln(2)\right)^4\zeta(2) - \frac{2}{3}\left(ln(2)\right)^3\zeta(2)
- \frac{1}{15}\left(ln(2)\right)^6 \nonumber \\ &+& \frac{1}{15}\left(ln(2)\right)^5 - \frac{1}{4}\zeta(2)\sum^{\infty}_{k=1}\frac{h_k}{k^3} + \frac{3}{8}\sum^{\infty}_{k=1}\frac{h_k}{k^5} +
8ln(2)Li_5\left(\frac{1}{2}\right) - 8Li_5\left(\frac{1}{2}\right) \nonumber \\ &-& \frac{1}{12}\frac{h_k^3}{k^3}~.
\end{eqnarray}
Using now the two-valued help function:
\begin{eqnarray}
\sum^{\infty}_{k=1} \frac{H_k}{2k-1} \left( \sum^{\infty}_{i=1} \frac{H_i^{(3)}}{(2i-1)(2i+2k-1)} \right) &=& \sum^{\infty}_{i=1} \frac{1}{i^3(2i-1)} \left(\sum^{\infty}_{k=1}\frac{H_k}{(2k-1)(2i+2k-1)}\right) +
\frac{1}{2}\zeta(3) \sum^{\infty}_{k=1} \frac{H_k h_k}{k(2k-1)} \nonumber \\ &-& \zeta(2)\sum^{\infty}_{k=1} \frac{H_k h_k^{(2)}}{k(2k-1)} - 4ln(2)\sum^{\infty}_{k=1} \frac{H_k h_k^{(3)}}{k(2k-1)} \nonumber \\ &+&
4\sum^{\infty}_{k=1} \frac{H_k}{k(2k-1)} \left( \sum^{k}_{i=1} \frac{h_i}{(2i-1)^3} \right)
\end{eqnarray}
or
\begin{eqnarray}
\sum^{\infty}_{i=1} \frac{H_i^{(3)}}{2i-1} \left( \sum^{\infty}_{k=1}\frac{H_k}{(2k-1)(2i+2k-1)} \right) &=& 2ln(2)\sum^{\infty}_{i=1} \frac{H_i^{(3)}}{(2i-1)^2} - ln(2)\sum^{\infty}_{i=1} \frac{H_i^{(3)} h_i}{i(2i-1)} -
\sum^{\infty}_{i=1} \frac{H_i^{(3)} h_i}{i(2i-1)^2} + \frac{1}{2} \sum^{\infty}_{i=1} \frac{H_i^{(3)} h_i^{(2)}}{i(2i-1)} \nonumber \\ &+& \frac{1}{2} \sum^{\infty}_{i=1} \frac{H_i^{(3)} h_i^2}{i(2i-1)}
\end{eqnarray}
we are able to define an identity between the two cubic Euler sums:
\begin{eqnarray}
\sum^{\infty}_{i=1} \frac{H_i^{(3)} h_i^2}{i(2i-1)} \sim \sum^{\infty}_{i=1}\frac{h_i^3}{i^3}~.
\end{eqnarray}
Using the ansatz
\begin{eqnarray}
\sum^{\infty}_{i=1} \frac{1}{i^2} \left( \sum^{\infty}_{k=1} \frac{h_k^2}{k(i-k)} \right) &=& \sum^{\infty}_{k=1} \frac{h_k^2}{k} \left( \sum^{\infty}_{i=1} \frac{1}{i^2(i-k)} \right) =
\sum^{\infty}_{k=1} \frac{h_k^2}{k} \left( \frac{3}{k^3} - \zeta(2)\frac{1}{k} - \frac{H_k}{k^2 } \right) \nonumber \\ &=& 3\frac{h_k^2}{k^4} - \zeta(2)\frac{h_k^2}{k} - \frac{H_k h_k^2}{k^3}
\end{eqnarray}
and inserting the corresponding two-valued function which results from the corresponding inhomogeneous first order difference equation in Eq.~(159)
\begin{eqnarray}
\sum^{\infty}_{k=1} \frac{h_k^2}{k(i-k)} = \frac{1}{2}\zeta(2)\frac{h_i}{i} + 2ln(2)\frac{h_i^2}{i} - \frac{2}{3}\frac{h_i^{(3)}}{i} - \frac{h_i^2}{i^2} - \frac{4}{3}\frac{h_i^3}{i} - 2\frac{h_i h_i^{(2)}}{i} +
\frac{1}{i}\sum^{i-1}_{k=1} \frac{h_k^2}{(i-k)}
\end{eqnarray}
we get by rearranging the summations in the corresponding nested double sum:
\begin{eqnarray}
\frac{4}{3}\sum^{\infty}_{k=1}\frac{h_k^3}{k^3} &=& \frac{1}{2}\zeta(2)\sum^{\infty}_{k=1}\frac{h_k}{k^3} - 3\sum^{\infty}_{k=1}\frac{h_k^2}{k^4} + 2ln(2)\sum^{\infty}_{k=1}\frac{h_k^2}{k^3} - \frac{2}{3}\sum^{\infty}_{k=1}
\frac{h_k^{(3)}}{k^3} - 2\sum^{\infty}_{k=1}\frac{h_k h_k^{(2)}}{k^3} + 2\sum^{\infty}_{k=1}\frac{H_k h_k^2}{k^3} + \sum^{\infty}_{k=1}\frac{H_k^{(2)} h_k^2}{k^2}  \nonumber \\ &-& \sum^{\infty}_{k=1}\frac{1}{k^3} 
\left( \sum^{k-1}_{i=1} \frac{h_i^2}{i} \right)~.
\end{eqnarray}
Following Adegoke \cite{ade16} it stands 
\begin{eqnarray}
\sum^{\infty}_{k=1}\frac{1}{k^3} \left( \sum^{k}_{i=1} \frac{h_i^2}{i} \right) &=& \frac{2}{3}\sum^{\infty}_{k=1}\frac{h_k^3}{k^3} - \frac{2}{3}\sum^{\infty}_{k=1}\frac{h_k^{(3)}}{k^3} + 2\sum^{\infty}_{k=1}\frac{1}{k^3}
\left( \sum^{k}_{i=1} \frac{h_i}{(2i-1)^2} \right) - \sum^{\infty}_{k=1}\frac{1}{k^3} \left( \sum^{k}_{i=1} \frac{h_i^2}{i(2i-1)} \right)~.
\end{eqnarray}
Rearranging the summations in the corresponding nested double sum it follows:
\begin{eqnarray}
\sum^{\infty}_{k=1}\frac{1}{k^3} \left( \sum^{k}_{i=1} \frac{h_i^2}{i} \right) &=& \frac{2}{3}\sum^{\infty}_{k=1}\frac{h_k^3}{k^3} - \frac{2}{3}\sum^{\infty}_{k=1}\frac{h_k^{(3)}}{k^3} + 2\zeta(3)
\sum^{\infty}_{k=1} \frac{h_k}{(2k-1)^2} + 2\sum^{\infty}_{k=1} \frac{h_k}{k^3(2k-1)^2} - 2\sum^{\infty}_{k=1} \frac{H_k^{(3)} h_k}{(2k-1)^2}  \nonumber \\ &-& \zeta(3)\sum^{\infty}_{k=1} \frac{h_k^2}{k(2k-1)} - 
\sum^{\infty}_{k=1} \frac{h_k^2}{k^4(2k-1)} + \sum^{\infty}_{k=1} \frac{H_k^{(3)} h_k^2}{k(2k-1)}~.
\end{eqnarray}
By inserting Eq.~(163) in Eq.~(161) and by use of the explicitly known Euler sums it follows after some algebraic manipulations:
\begin{eqnarray}
\sum^{\infty}_{k=1}\frac{H_k^{(2)} h_k^2}{k^2} &=& \frac{225}{16}\zeta(6) - \frac{259}{48}\zeta(3)^2 - \zeta(2)\sum^{\infty}_{k=1}\frac{h_k}{k^3} + \sum^{\infty}_{k=1}\frac{h_k}{k^5} - 2\sum^{\infty}_{k=1}\frac{H_k h_k^2}{k^3}
+ \frac{4}{3} \sum^{\infty}_{k=1}\frac{h_k^3}{k^3}~. 
\end{eqnarray}
This way we have defined a second type of an essential identity between the unknown cubic Euler sums appearing in Eq.~(128) and Eq.~(143). As a consequence we need a third identity of this type. It is advantageous to use the
two-valued help function:
\begin{eqnarray}
\sum^{\infty}_{i=1} \frac{H_i}{i} \left(\sum^{\infty}_{k=1} \frac{h_k^{(2)}}{k(i+k)} \right) &=& \frac{3}{2}\zeta(2)\sum^{\infty}_{k=1} \frac{H_k h_k}{k^2} - 2ln(2)\sum^{\infty}_{k=1}\frac{H_k h_k^{(2)}}{k^2} -
\sum^{\infty}_{i=1} \frac{H_i}{i^2} \left( \sum^{i}_{k=1}\frac{H_{k-1}}{(2k-1)^2} \right)
\end{eqnarray}
or
\begin{eqnarray}
\sum^{\infty}_{k=1} \frac{h_k^{(2)}}{k} \left( \sum^{\infty}_{i=1}\frac{H_i}{i(k+i)} \right) &=& \zeta(2)\sum^{\infty}_{k=1}\frac{h_k^{(2)}}{k^2} - \sum^{\infty}_{k=1}\frac{H_k h_k^{(2)}}{k^3} + \frac{1}{2}
\sum^{\infty}_{k=1}\frac{H_k^{(2)} h_k^{(2)}}{k^2} + \frac{1}{2}\sum^{\infty}_{k=1}\frac{H_k^2 h_k^{(2)}}{k^2}~.
\end{eqnarray}
This way it follows
\begin{eqnarray}
\sum^{\infty}_{k=1} \frac{H_k}{k^2} \left( \sum^{k}_{i=1}\frac{H_{i-1}}{(2i-1)^2} \right) &=& \frac{3}{2}\zeta(2)\sum^{\infty}_{k=1}\frac{H_k h_k}{k^2} - 2ln(2)\sum^{\infty}_{k=1}\frac{H_k h_k^{(2)}}{k^2} - 
\zeta(2)\sum^{\infty}_{k=1}\frac{h_k^{(2)}}{k^2} + \sum^{\infty}_{k=1}\frac{H_k h_k^{(2)}}{k^3} - \frac{1}{2} \sum^{\infty}_{k=1}\frac{H_k^{(2)} h_k^{(2)}}{k^2} \nonumber \\ &-& \frac{1}{2}
\sum^{\infty}_{k=1}\frac{H_k^2 h_k^{(2)}}{k^2}~.
\end{eqnarray}
Furthermore we calculate
\begin{eqnarray}
\sum^{\infty}_{k=1} \frac{h_k}{k^2} \left( \sum^{\infty}_{i=1}\frac{H_i}{i(k+i)}\right) &=& \zeta(2)\sum^{\infty}_{k=1} \frac{h_k}{k^3} + \sum^{\infty}_{k=1} \frac{h_k}{k^3} \left( \sum^{k-1}_{i=1}\frac{H_i}{i} \right)
\nonumber \\ &=& \zeta(2)\sum^{\infty}_{k=1} \frac{h_k}{k^3} - \sum^{\infty}_{k=1} \frac{H_k h_k}{k^4} + \frac{1}{2}\sum^{\infty}_{k=1} \frac{H_k^{(2)} h_k}{k^3} + \frac{1}{2}\sum^{\infty}_{k=1}\frac{H_k^2 h_k}{k^3}
\end{eqnarray}
or
\begin{eqnarray}
\sum^{\infty}_{i=1}\frac{H_i}{i} \left( \sum^{\infty}_{k=1} \frac{h_k}{k^2(i+k)} \right) &=& \sum^{\infty}_{i=1}\frac{H_i}{i^2} \left( \sum^{\infty}_{k=1} \frac{h_k}{k^2} \right) -2ln(2)\sum^{\infty}_{k=1} \frac{H_k h_k}{k^3} 
- \sum^{\infty}_{k=1} \frac{H_k^2 h_k}{k^3} + \sum^{\infty}_{k=1} \frac{H_k}{k^3} \left( \sum^{k}_{i=1}\frac{h_{i}}{i} \right)~.
\end{eqnarray}
This way it follows:
\begin{eqnarray}
\sum^{\infty}_{k=1} \frac{H_k}{k^3} \left( \sum^{k}_{i=1}\frac{h_{i}}{i} \right) &=& \zeta(2)\sum^{\infty}_{k=1} \frac{h_k}{k^3} - \sum^{\infty}_{i=1}\frac{H_i}{i^2} \left( \sum^{\infty}_{k=1} \frac{h_k}{k^2} \right)
- \sum^{\infty}_{k=1} \frac{H_k h_k}{k^4} + 2ln(2)\sum^{\infty}_{k=1} \frac{H_k h_k}{k^3} + \frac{1}{2}\sum^{\infty}_{k=1}\frac{H_k^{(2)} h_k}{k^3} \nonumber \\ &+& \frac{3}{2}\sum^{\infty}_{k=1}\frac{H_k^2 h_k}{k^3}~.
\end{eqnarray}
We proceed further in calculating the following expressions:
\begin{eqnarray}
\sum^{\infty}_{i=1}\frac{H_i}{i^2} \left( \sum^{\infty}_{k=1} \frac{h_k}{(i+k)^2} \right) &=& \frac{7}{4}\zeta(3)\sum^{\infty}_{k=1} \frac{H_k}{k^2} + \zeta(2)\sum^{\infty}_{k=1} \frac{H_k h_k}{k^2} - 4ln(2)\sum^{\infty}_{k=1}
\frac{H_k h_k^{(2)}}{k^2} - 2\sum^{\infty}_{k=1} \frac{H_k}{k^2} \left( \sum^{k}_{i=1}\frac{H_{i-1}}{(2i-1)^2} \right) \nonumber \\ &-& \sum^{\infty}_{k=1} \frac{H_k}{k^2} \left( \sum^{k}_{i=1}\frac{H_{i-1}^{(2)}}{2i-1} \right)
\end{eqnarray}
and
\begin{eqnarray}
\sum^{\infty}_{i=1}\frac{h_i}{i^2} \left( \sum^{\infty}_{k=1} \frac{H_k}{(i+k)^2} \right) &=& \zeta(3)\sum^{\infty}_{k=1} \frac{h_k}{k^2} - \zeta(2)\sum^{\infty}_{k=1}\frac{h_k}{k^3} + \zeta(2)\sum^{\infty}_{k=1}
\frac{H_k h_k}{k^2} - \sum^{\infty}_{k=1} \frac{H_k^{(3)} h_k}{k^2} - \sum^{\infty}_{k=1} \frac{H_k H_k^{(2)} h_k}{k^2} + \sum^{\infty}_{k=1} \frac{H_k h_k}{k^4} \nonumber \\ &+& \sum^{\infty}_{k=1} \frac{H_k^{(2)} h_k}{k^3}~. 
\end{eqnarray}
The difference between these two expressions results in
\begin{eqnarray}
\sum^{\infty}_{i=1}\frac{H_i}{i^2} \left( \sum^{\infty}_{k=1} \frac{h_k}{(i+k)^2} \right) - \sum^{\infty}_{i=1}\frac{h_i}{i^2} \left( \sum^{\infty}_{k=1} \frac{H_k}{(i+k)^2} \right) &=& 
\sum^{\infty}_{k=1}\frac{H_k}{k^2} \left( \sum^{\infty}_{i=1} \frac{h_i}{i(k+i)} \right) - \sum^{\infty}_{k=1}\frac{h_k}{k^2} \left( \sum^{\infty}_{i=1} \frac{H_i}{i(k+i)} \right)~.
\end{eqnarray}
The explicit calculation gives us:
\begin{eqnarray}
\frac{7}{2}\zeta(3)\sum^{\infty}_{k=1} \frac{H_k}{k^2} &-& 2\zeta(3)\sum^{\infty}_{k=1} \frac{h_k}{k^2} + 4\zeta(2)\sum^{\infty}_{k=1} \frac{h_k}{k^3} - 4ln(2)\sum^{\infty}_{k=1} \frac{H_k h_k}{k^3} -
4\sum^{\infty}_{k=1} \frac{H_k h_k}{k^4} - 8ln(2)\sum^{\infty}_{k=1} \frac{H_k^{(2)} h_k}{k^2} \nonumber \\ &-& 2\sum^{\infty}_{k=1} \frac{H_k^{(2)} h_k}{k^3} + 2\sum^{\infty}_{k=1} \frac{H_k^{(3)} h_k}{k^2} +
\sum^{\infty}_{k=1} \frac{H_k^2 h_k}{k^3} + 2\sum^{\infty}_{k=1} \frac{H_k H_k^{(2)} h_k}{k^2}  = 4\sum^{\infty}_{k=1} \frac{H_k}{k^2} \left( \sum^{k}_{i=1}\frac{H_{i-1}}{(2i-1)^2} \right) \nonumber \\ &+&
2\sum^{\infty}_{k=1} \frac{H_k}{k^2} \left( \sum^{k}_{i=1}\frac{H_{i-1}^{(2)}}{2i-1} \right) + 2\sum^{\infty}_{k=1} \frac{H_k}{k^3} \left( \sum^{k}_{i=1}\frac{H_{i-1}}{2i-1} \right)~.
\end{eqnarray}
Next we need to calculate the last cubic Euler sum on the left side. It follows:
\begin{eqnarray}
\sum^{\infty}_{k=1} \frac{H_k H_k^{(2)} h_k}{k^2} &=& \sum^{\infty}_{k=1} \frac{H_k}{k^2} \left( \sum^{k}_{i=1}\frac{h_i}{i^2} + \sum^{k}_{i=1}\frac{H_{i-1}^{(2)}}{2i-1} \right) = \sum^{\infty}_{k=1} \frac{H_k}{k^2}
\left( \sum^{k}_{i=1} \frac{h_i}{i^2} \right) + \sum^{\infty}_{k=1} \frac{H_k}{k^2} \left( \sum^{k}_{i=1}\frac{H_{i-1}^{(2)}}{2i-1} \right)~.
\end{eqnarray}
By use of the following identity we get:
\begin{eqnarray}
\sum^{\infty}_{k=1} \frac{H_k}{k^2} \left( \sum^{k}_{i=1} \frac{h_i}{i(k+i)} \right) = \sum^{\infty}_{k=1} \frac{H_k}{k^2} \left( \frac{H_k h_k}{2k} \right) = \frac{1}{2}\sum^{\infty}_{k=1} \frac{H_k^2 h_k}{k^3}~. 
\end{eqnarray}
Rearranging the summations in the nested double sum it follows:
\begin{eqnarray}
\sum^{\infty}_{k=1} \frac{H_k}{k^2} \left( \sum^{k}_{i=1} \frac{h_i}{i^2} \right) &=& \sum^{\infty}_{k=1} \frac{H_k h_k^2}{k^3} - \sum^{\infty}_{k=1} \frac{h_k^2}{k^4} + \zeta(2)\sum^{\infty}_{k=1} \frac{h_k}{k^3}
- \frac{1}{2}\sum^{\infty}_{k=1} \frac{H_k^{(2)} h_k}{k^3} + \sum^{\infty}_{k=1} \frac{H_k h_k}{k^4}~.
\end{eqnarray}
Inserting first Eq.~(177) in Eq.~(175) and then inserting this result in Eq.~(174) we get:
\begin{eqnarray}
&&\frac{7}{2}\zeta(3)\sum^{\infty}_{k=1} \frac{H_k}{k^2} - 2\zeta(3)\sum^{\infty}_{k=1} \frac{h_k}{k^2} + 6\zeta(2))\sum^{\infty}_{k=1} \frac{h_k}{k^3} -4ln(2)\sum^{\infty}_{k=1} \frac{H_k h_k}{k^3} -
2\sum^{\infty}_{k=1} \frac{H_k h_k}{k^4} - 8ln(2)\sum^{\infty}_{k=1} \frac{H_k h_k^{(2)}}{k^2} - 2\sum^{\infty}_{k=1} \frac{h_k^2}{k^4} \nonumber \\ &-& 2\sum^{\infty}_{k=1} \frac{H_k^{(2)} h_k}{k^3} + 2\sum^{\infty}_{k=1}
\frac{H_k^{(3)} h_k}{k^2} + \sum^{\infty}_{k=1} \frac{H_k^2 h_k}{k^3} + 2\sum^{\infty}_{k=1} \frac{H_k h_k^2}{k^3} = 4\frac{H_k}{k^2} \left( \sum^{k}_{i=1}\frac{H_{i-1}}{(2i-1)^2} \right) \nonumber \\ &+&
2\sum^{\infty}_{k=1} \frac{H_k}{k^3} \left( \sum^{k}_{i=1}\frac{H_{i-1}}{2i-1} \right)~.
\end{eqnarray}
Finally, using the well known identity 
\begin{eqnarray}
\sum^{k}_{i=1}\frac{H_{i-1}}{2i-1} = H_k h_k - \sum^{k}_{i=1} \frac{h_i}{i}
\end{eqnarray}
Eq.~(167), Eq.~(170) and the explicit form of all known Euler sums the third essential identity between two cubic Euler sums appearing in Eq.~(128) and Eq.~(143) follows after some algebraic manipulations: 
\begin{eqnarray}
\sum^{\infty}_{k=1} \frac{H_k^2 h_k}{k^3} + \sum^{\infty}_{k=1} \frac{H_k h_k^2}{k^3} + \sum^{\infty}_{k=1} \frac{H_k^2 h_k^{(2)} }{k^2} &=& \frac{585}{64}\zeta(6) + \frac{49}{16}\zeta(3)^2 - 
2\zeta(2)\sum^{\infty}_{k=1}\frac{h_k}{k^3} - \frac{1}{4}\sum^{\infty}_{k=1}\frac{h_k}{k^5}~.
\end{eqnarray}
What remains is to found an identity where Eq.~(175) and corresponding members of Eq.~(128) are involved. This can be done by inserting Eq.~(167) in Eq.~(171) and by rearranging the summations on the left side of Eq.~(171). 
From this it follows after some algebraic manipulations:
\begin{eqnarray}
\sum^{\infty}_{k=1} \frac{H_k}{k^2} \left( \sum^{k}_{i=1}\frac{H_{i-1}^{(2)}}{2i-1} \right) &=& \zeta(3)\sum^{\infty}_{k=1} \frac{h_k}{k^2} - \zeta(2)\sum^{\infty}_{k=1} \frac{h_k}{k^3} - 2ln(2)\sum^{\infty}_{k=1}
\frac{H_k^{(2)} h_k}{k^2} + \zeta(2)\sum^{\infty}_{k=1}\frac{H_k h_k}{k^2} - \sum^{\infty}_{k=1} \frac{h_k}{k^2} \left( \sum^{\infty}_{k=i} \frac{H_i}{i^2} \right) \nonumber \\ &+& \sum^{\infty}_{k=1} \frac{H_k^{(2)}}{k^2}
\left( \sum^{k}_{i=1} \frac{h_i}{i} \right)~,
\end{eqnarray}
where the last nested sum on the right side of Eq.~(181) explicitly results from the ansatz:
\begin{eqnarray}
\sum^{\infty}_{k=1} \frac{h_k}{2k-1} \left( \sum^{k}_{i=1}\frac{H_i^{(2)}}{i(2i+2k-1)} \right) &=& \left( \frac{5}{2}\zeta(3)-2ln(2)\zeta(2) \right) \sum^{\infty}_{k=1}\frac{h_k}{(2k-1)^2} + 2\zeta(2)\sum^{\infty}_{k=1}
\frac{h_k h_{k-1}}{(2k-1)^2} + 8ln(2)\sum^{\infty}_{k=1} \frac{h_k h_{k-1}^{(2)}}{(2k-1)^2} \nonumber \\ &-& 8\sum^{\infty}_{k=1} \frac{h_k}{(2k-1)^2} \left( \sum^{k-1}_{i=1} \frac{h_i}{(2i-1)^2} \right)~.
\end{eqnarray}
This way it follows explicitly after some algebraic manipulations
\begin{eqnarray}
\sum^{\infty}_{k=1} \frac{H_k}{k^2} \left( \sum^{k}_{i=1}\frac{H_{i-1}^{(2)}}{2i-1} \right) &=& \frac{315}{64}\zeta(6) - \frac{49}{16}\zeta(3)^2 + \frac{1}{2} \zeta(2)\sum^{\infty}_{k=1}\frac{h_k}{k^3}~.
\end{eqnarray}
Inserting this result in Eq.~(175) the corresponding identity follows to:
\begin{eqnarray}
\sum^{\infty}_{k=1} \frac{H_k H_k^{(2)} h_k}{k^2} &=& \frac{45}{64}\zeta(6) + \zeta(2)\sum^{\infty}_{k=1}\frac{h_k}{k^3} + \frac{1}{4}\sum^{\infty}_{k=1}\frac{h_k}{k^5} + \sum^{\infty}_{k=1} \frac{H_k h_k^2}{k^3}~. 
\end{eqnarray}
The last identity follows immediately by inserting Eq.~(147), Eq.~(163), Eq.~(179) and Eq.~(184) in Eq.~(143). We get:
\begin{eqnarray}
\sum^{\infty}_{k=1} \frac{H_k h_k h_k^{(2)}}{k^2} &=& \frac{45}{32}\zeta(6) + \frac{203}{96}\zeta(3)^2 + \frac{1}{4}\zeta(2)\sum^{\infty}_{k=1}\frac{h_k}{k^3} -  \frac{1}{4}\sum^{\infty}_{k=1}\frac{h_k}{k^5}
- \frac{1}{3}\sum^{\infty}_{k=1}\frac{h_k^3}{k^3}~.
\end{eqnarray}
Inserting now in Eq.~(143) the cubic Euler sum 
\begin{eqnarray}
\sum^{\infty}_{k=1} \frac{H_k^2 H_k^{(2)}}{k^2} = \frac{41}{12}\zeta(6) + 2\zeta(3)^2
\end{eqnarray}
which is known from literature \cite{ce20} together with the five identities calculated before it stands: 
\begin{eqnarray}
\frac{1}{4}\sum^{\infty}_{k=1}\frac{H_k h_k^2}{k^3} &-& \frac{1}{8}\sum^{\infty}_{k=1}\frac{H_k^2 h_k}{k^3} - \frac{1}{6}\sum^{\infty}_{k=1}\frac{h_k^3}{k^3} = \frac{5727}{768}\zeta(6) - \frac{53}{16}\left(ln(2)\right)^2\zeta(4)
+ \frac{7}{48}\zeta(3)^2 + \frac{7}{3}\left(ln(2)\right)^3\zeta(3) - \frac{1}{2}\left(ln(2)\right)^4\zeta(2) \nonumber \\ &+& \frac{1}{18}\left(ln(2)\right)^6 - \frac{1}{2}\left(ln(2)\right)^2\sum^{\infty}_{k=1}\frac{h_k}{k^3}
+ \frac{1}{4}\sum^{\infty}_{k=1}\frac{h_k}{k^5} - 8Li_6\left(\frac{1}{2}\right) - 8ln(2)Li_5\left(\frac{1}{2}\right)~.
\end{eqnarray}
Substituting now the next to the last cubic Euler sum in Eq.~(128) by use of Eq.~(187) and doing some algebra Eq.~(188) results. In analogy one may substitute the last cubic Euler sum in Eq.~(128) by use of Eq.~(187). Then
Eq.~(188) results, and in consequence also the first member of the eighth family can be calculated explicitly. Due to Eq.~(147) also the second member of the the eighth family is known. Furthermore Eq.~(190) - Eq.~(193)  
result. With this all members of the first and eighth family are known explicitly. Finally, because of Eq.~(132) and  Eq.~(133) the first and second member of the fifth family are known explicitly. Using Eq.~(140) the
first member of the fourth family can be calculated explicitly and by use of Eq.~(204) the second member too. This way all members of the fourth family are known explicitly. 

\begin{eqnarray}
\sum^{\infty}_{k=1} \frac{H_k^2 h_k}{k^3} &=& -\frac{1639}{32}\zeta(6) + \frac{53}{2}\left(ln(2)\right)^2\zeta(4) - \frac{77}{16}\zeta(3)^2 - \frac{56}{3}\left(ln(2)\right)^3\zeta(3) +
4\left(ln(2)\right)^4\zeta(2) - \frac{4}{9}\left(ln(2)\right)^6  \nonumber \\ &-& \zeta(2)\sum^{\infty}_{k=1}\frac{h_k}{k^3} + 4\left(ln(2)\right)^2\sum^{\infty}_{k=1}\frac{h_k}{k^3} -
\frac{7}{4}\sum^{\infty}_{k=1}\frac{h_k}{k^5} + 64Li_6\left(\frac{1}{2}\right) + 64ln(2)Li_5\left(\frac{1}{2}\right)
\end{eqnarray}
and
\begin{eqnarray}
\sum^{\infty}_{k=1} \frac{H_k h_k^2}{k^3} =  \frac{315}{64}\zeta(6) - \frac{49}{32}\zeta(3)^2 - \frac{1}{2}\zeta(2)\sum^{\infty}_{k=1}\frac{h_k}{k^3} +
\frac{1}{8}\sum^{\infty}_{k=1}\frac{h_k}{k^5}
\end{eqnarray}
and
\begin{eqnarray}
\sum^{\infty}_{k=1} \frac{H_k^{(2)} h_k^2}{k^2} =  \frac{45}{8}\zeta(6) - \frac{7}{4}\zeta(3)^2 + \frac{3}{4}\sum^{\infty}_{k=1}\frac{h_k}{k^5}
\end{eqnarray}
and
\begin{eqnarray}
\sum^{\infty}_{k=1} \frac{H_k^2 h_k^{(2)}}{k^2} &=& \frac{887}{16}\zeta(6) - \frac{53}{2}\left(ln(2)\right)^2\zeta(4) + \frac{301}{32}\zeta(3)^2 + \frac{56}{3}\left(ln(2)\right)^3\zeta(3) -
4\left(ln(2)\right)^4\zeta(2) + \frac{4}{9}\left(ln(2)\right)^6  \nonumber \\ &-& \frac{1}{2}\zeta(2)\sum^{\infty}_{k=1}\frac{h_k}{k^3} - 4\left(ln(2)\right)^2\sum^{\infty}_{k=1}\frac{h_k}{k^3} +
\frac{11}{8}\sum^{\infty}_{k=1}\frac{h_k}{k^5} - 64Li_6\left(\frac{1}{2}\right) - 64ln(2)Li_5\left(\frac{1}{2}\right)
\end{eqnarray}
and
\begin{eqnarray}
\sum^{\infty}_{k=1} \frac{H_k H_k^{(2)} h_k}{k^2} = \frac{45}{8}\zeta(6) - \frac{49}{32}\zeta(3)^2 + \frac{1}{2}\zeta(2)\sum^{\infty}_{k=1}\frac{h_k}{k^3} +
\frac{3}{8}\sum^{\infty}_{k=1}\frac{h_k}{k^5}
\end{eqnarray}
and
\begin{eqnarray}
\sum^{\infty}_{k=1} \frac{H_k h_k h_k^{(2)}}{k^2} = \frac{135}{128}\zeta(6) + \frac{63}{32}\zeta(3)^2 + \frac{1}{4}\zeta(2)\sum^{\infty}_{k=1}\frac{h_k}{k^3} -
\frac{1}{4}\sum^{\infty}_{k=1}\frac{h_k}{k^5}~.
\end{eqnarray}

\subsection{Second family}
Now we introduce the nonlinear two-valued help function in order to calculate the first member of the second family.
\begin{eqnarray}
\sum^{\infty}_{i=1} \frac{H_i^2}{(2i-1)(k+i)} &=& \frac{3\zeta(3)}{2k+1} + \frac{1}{2k+1} \sum^{\infty}_{i=1}\frac{H_i^2}{i(2i-1)} + \zeta(2)\frac{H_{k-1}}{2k+1} + \frac{2}{3}\frac{H_{k-1}^{(3)}}{2k+1}
+ \frac{H_{k-1} H_{k-1}^{(2)}}{2k+1} + \frac{1}{3}\frac{H_{k-1}^3}{2k+1} \nonumber \\ &-& \frac{1}{2k+1}\sum^{k-1}_{i=1} \frac{H_i}{i^2}
\end{eqnarray}
which again represents the solution of the corresponding inhomogeneous first order difference equation. With this we get:
\begin{eqnarray}
\sum^{\infty}_{k=1} \frac{1}{(2k+1)^2} \left( \sum^{\infty}_{i=1} \frac{H_i^2}{(2i-1)(k+i)} \right) &=&  \sum^{\infty}_{i=1} \frac{H_i^2}{2i-1} \left( \sum^{\infty}_{k=1} \frac{1}{(2k+1)^2(i+k)} \right)
\end{eqnarray} 
or
\begin{eqnarray}
\sum^{\infty}_{k=1} \frac{1}{(2k+1)^2} \left( \sum^{\infty}_{i=1} \frac{H_i^2}{(2i-1)(k+i)} \right) &=& \sum^{\infty}_{k=1} \frac{2}{(2k+1)^2} \left( \sum^{\infty}_{i=1} \frac{H_i^2}{(2i-1)^2} \right)
- \sum^{\infty}_{k=1} \frac{1}{(2k+1)(i+k)} \left( \sum^{\infty}_{i=1} \frac{H_i^2}{(2i-1)^2} \right) \nonumber \\ &&
\sum^{\infty}_{k=1} \frac{2}{(2k+1)^2} \left( \sum^{\infty}_{i=1} \frac{H_i^2}{(2i-1)^2} \right) + \left( 2 - 2ln(2) \right) \sum^{\infty}_{k=1}\frac{H_i^2}{(2i-1)^3} \nonumber \\ &&
- \sum^{\infty}_{k=1}\frac{H_i^3}{(2i-1)^3}~.
\end{eqnarray} 
Inserting Eq.~(194) in Eq.~(196) and rearranging the different terms it follows after an elongate algebraic procedure
\begin{eqnarray}
\sum^{\infty}_{k=1} \frac{H_k^3}{(2k-1)^3} &=& \frac{405}{128}\zeta(6) - \frac{93}{4}ln(2)\zeta(5) + \frac{93}{4}\zeta(5) + \frac{135}{8}\left(ln(2)\right)^2\zeta(4) - 
\frac{135}{4}ln(2)\zeta(4) - \frac{7}{32}\zeta(3)^2 + \frac{21}{4}ln(2)\zeta(2)\zeta(3) \nonumber \\ &-& \frac{21}{4}\zeta(2)\zeta(3) - 7\left(ln(2)\right)^3\zeta(3)
+ 21\left(ln(2)\right)^2\zeta(3) + 21ln(2)\zeta(3) - \frac{63}{2}\zeta(3) - 18\left(ln(2)\right)^2\zeta(2) \nonumber \\ &+& 30ln(2)\zeta(2) - 6\zeta(2) 
+ 8\left(ln(2)\right)^3	 - 36\left(ln(2)\right)^2 + 48ln(2) - 6\sum^{\infty}_{k=1}\frac{h_k}{k^3} + \frac{3}{4}\sum^{\infty}_{k=1}\frac{h_k}{k^5}~.
\end{eqnarray}

The second member of the second family results from the following ansatz:
\begin{eqnarray}
\sum^{\infty}_{k=1} \frac{H_k^{(2)}}{2k+1} \left(\sum^{\infty}_{i=1} \frac{H_i}{(2i-1)(k+i)} \right) = \sum^{\infty}_{i=1} \frac{H_i}{(2i-1)}
\left( \sum^{\infty}_{k=1} \frac{H_k^{(2)}}{(2k+1)(k+i)} \right)
\end{eqnarray}
or
\begin{eqnarray}
&& \left( 4ln(2) - 2\left(ln(2)\right)^2 + \zeta(2) \right) \sum^{\infty}_{k=1} \frac{H_{k-1}^{(2)}}{(2k-1)^2} - \sum^{\infty}_{k=1} \frac{H_k H_{k}^{(2)}}{k(2k+1)^2} + \frac{1}{2}
\sum^{\infty}_{k=1} \frac{H_{k-1}^{(2)} H_{k-1}^{(2)}}{(2k-1)^2} + \frac{1}{2}\sum^{\infty}_{k=1} \frac{H_{k-1}^2 H_{k-1}^2}{(2k-1)^2} = \nonumber \\ &&  \left( 2ln(2)\zeta(2) -\frac{3}{2}\zeta(3) \right)
\sum^{\infty}_{i=1} \frac{H_i}{(2i-1)^2} - \zeta(2) \sum^{\infty}_{i=1} \frac{H_{i}}{i(2i-1)^2} + \zeta(2) \sum^{\infty}_{i=1} \frac{H_i^2}{(2i-1)^2} - \sum^{\infty}_{i=1} \frac{H_i}{(2i-1)^2}
 \left( \sum^{i}_{k=1} \frac{H_k}{k^2} \right)~.
\end{eqnarray}
Only the last double sum on the right side of Eq.~(199) in unkwnon. Using the ansatz:
\begin{eqnarray}
\sum^{\infty}_{i=1} \frac{H_i}{i} \left(\sum^{\infty}_{k=1} \frac{h_k^{(2)}}{k(i+k)} \right) = \sum^{\infty}_{k=1} \frac{h_k^{(2)}}{k} \left( \sum^{\infty}_{i=1} \frac{H_i}{i(k+i)} \right) 
\end{eqnarray}
or
\begin{eqnarray}
&& \zeta(2)\sum^{\infty}_{k=1} \frac{h_k^{(2)}}{k^2} - \sum^{\infty}_{k=1} \frac{H_k h_k^{(2)}}{k^3} + \frac{1}{2}\sum^{\infty}_{k=1} \frac{H_k^{(2)} h_k^{(2)}}{k^2} + 
\frac{1}{2}\sum^{\infty}_{k=1} \frac{H_k^2 h_k^{(2)}}{k^2} = \nonumber \\ && \frac{3}{2}\zeta(2)\sum^{\infty}_{k=1} \frac{H_k h_k}{k^2} - 2ln(2)\zeta(2)\sum^{\infty}_{k=1} \frac{H_k h_k^{(2)}}{k^2} 
- \sum^{\infty}_{k=1} \frac{H_k}{k^2} \left( \sum^{k}_{i=1} \frac{H_{i-1}}{(2i-1)^2} \right)
\end{eqnarray}
the corresponding nested double sum in  Eq.~(199) follows explicitly by changing the summation order in the last nested double sum on the right side of Eq.~(201). Inserting the corresponding result
in Eq.~(199) the second member follows explicitly after a variety of standard algebraic manipulations
\begin{eqnarray}
\sum^{\infty}_{k=1} \frac{H_k^2 H_k^{(2)}}{(2k-1)^2} &=& -\frac{3593}{64}\zeta(6) + \frac{589}{8}ln(2)\zeta(5) - \frac{589}{8}\zeta(5) - \frac{121}{4}\left(ln(2)\right)^2\zeta(4) + 
\frac{121}{2}ln(2)\zeta(4) + \frac{65}{4}\zeta(4) - \frac{49}{32}\zeta(3)^2 \nonumber \\ &-& 7ln(2)\zeta(2)\zeta(3) + 7\zeta(2)\zeta(3) + \frac{28}{3}\left(ln(2)\right)^3\zeta(3)
- 28\left(ln(2)\right)^2\zeta(3) - 42ln(2)\zeta(3) + 36\zeta(3) \nonumber \\ &-& \frac{4}{3}\left(ln(2)\right)^4\zeta(2) +  \frac{16}{3}\left(ln(2)\right)^3\zeta(2) +
16\left(ln (2)\right)^2\zeta(2) - 24ln(2)\zeta(2) + \frac{4}{45}\left(ln(2)\right)^6 - \frac{8}{15}\left(ln(2)\right)^5  \nonumber \\ &-& \frac{32}{3}\left(ln(2)\right)^3 + 48\left(ln(2)\right)^2 - 64ln(2) +
10\sum^{\infty}_{k=1}\frac{h_k}{k^3} - \frac{19}{8}\sum^{\infty}_{k=1}\frac{h_k}{k^5} + 64Li_6\left(\frac{1}{2}\right) + 64Li_5\left(\frac{1}{2}\right)~.
\end{eqnarray}

\subsection{Third family}
In order to calculate the first member of the third family we use the ansatz:
\begin{eqnarray}
\sum^{\infty}_{k=1} \frac{H_k^{(2)}}{k(2k-1)} \left( \sum^{k}_{i=1} \frac{H_i^{(2)}}{i(2i-1)}\right) = \frac{1}{2} \left(\sum^{\infty}_{k=1} \frac{H_k^{(2)}}{k(2k-1)}\right)^2 +
\frac{1}{2}\sum^{\infty}_{k=1} \frac{H_k^{(2)} H_k^{(2)}}{k^2(2k-1)^2} = \nonumber \\
2\sum^{\infty}_{k=1} \frac{H_k^{(2)}}{k(2k-1)} \left( \sum^{k}_{i=1} \frac{H_i^{(2)}}{2i-1}\right) - \sum^{\infty}_{k=1} \frac{H_k^{(2)}}{k(2k-1)} \left( \sum^{k}_{i=1} \frac{H_i^{(2)}}{i}\right) =
\nonumber \\ 2\sum^{\infty}_{k=1} \frac{H_k^{(2)}}{k(2k-1)} \left( \sum^{k}_{i=1} \frac{H_i^{(2)}}{2i-1}\right) - \sum^{\infty}_{k=1} \frac{H_k^{(2)}}{k(2k-1)} \left( H_k H_k^{(2)} + H_k^{(3)} -
\frac{H_i}{i^2}\right) = \nonumber \\ 2\sum^{\infty}_{k=1} \frac{H_k^{(2)}}{k(2k-1)} \left( \sum^{k}_{i=1} \frac{H_i^{(2)}}{2i-1}\right) - \sum^{\infty}_{k=1} \frac{H_k^{(2)} H_k^{(3)}}{k(2k-1)} 
- \sum^{\infty}_{k=1} \frac{H_k H_k^{(2)} H_k^{(2)}}{k(2k-1)} + \sum^{\infty}_{k=1} \frac{H_k^{(2)}}{k(2k-1)} \left( \sum^{k}_{i=1} \frac{H_i}{i^2}\right) = \nonumber \\
2 \sum^{\infty}_{k=1} \frac{H_k^{(2)} H_k^{(2)} h_k}{k(2k-1)} - 2\sum^{\infty}_{k=1} \frac{H_k^{(2)}}{k(2k-1)} \left( \sum^{k}_{i=1} \frac{h_i}{i^2}\right) - \sum^{\infty}_{k=1} \frac{H_k^{(2)} H_k^{(3)}}{k(2k-1)}
- \sum^{\infty}_{k=1} \frac{H_k H_k^{(2)} H_k^{(2)}}{k(2k-1)} + \sum^{\infty}_{k=1} \frac{H_k^{(2)}}{k(2k-1)} \left( \sum^{k}_{i=1} \frac{H_i}{i^2}\right)~. \nonumber \\
\end{eqnarray}
The first sum in the last row of Eq.~(203) can be identified as the second member of the seventh family. In addition two nested double sums must be explicitly calculated.
As a first step We use the fifth member of the first family (Eq.~(192)) in order to write: 
\begin{eqnarray}
\sum^{\infty}_{k=1} \frac{H_k H_k^{(2)} h_k}{k^2} &=& \sum^{\infty}_{k=1} \frac{H_k}{k^2} \left( \sum^{k}_{i=1} \frac{h_i}{i^2} + \sum^{k}_{i=1} \frac{H_{i-1}^{(2)}}{2i-1} \right)~.
\end{eqnarray}
By use of Eq.~(177) we get from this:
\begin{eqnarray}
\sum^{\infty}_{k=1} \frac{H_k}{k^2} \left( \sum^{k}_{i=1} \frac{H_{i-1}^{(2)}}{2i-1} \right) &=& \sum^{\infty}_{k=1} \frac{H_k H_k^{(2)} h_k}{k^2} - \sum^{\infty}_{k=1} \frac{H_k h_k}{k^4}
+  \sum^{\infty}_{k=1} \frac{h_k^2}{k^4} - \zeta(2)\sum^{\infty}_{k=1} \frac{h_k}{k^3} + \frac{1}{2}\sum^{\infty}_{k=1} \frac{H_k^{(2)} h_k}{k^3} - \sum^{\infty}_{k=1} \frac{H_k^2 h_k}{k^3}
\end{eqnarray}
or explicitly
\begin{eqnarray}
\sum^{\infty}_{k=1} \frac{H_k}{k^2} \left( \sum^{k}_{i=1} \frac{H_{i-1}^{(2)}}{2i-1} \right) &=& \frac{315}{64}\zeta(6) - \frac{49}{16}\zeta(3)^2 + \frac{1}{2}\zeta(2)\sum^{\infty}_{k=1}\frac{h_k}{k^3}~.
\end{eqnarray}
As a next step we calculate the nested double sum
\begin{eqnarray}
\sum^{\infty}_{k=1} \frac{H_k}{k^2} \left( \sum^{k}_{i=1} \frac{H_{i-1}^{(2)}}{i} \right) &=& \sum^{\infty}_{k=1} \frac{H_k}{k^2} \left( H_k H_k^{(2)} + H_k^{(3)} - \sum^{k}_{i=1} \frac{H_{i}}{i^2} \right)
= \sum^{\infty}_{k=1} \frac{H_k^2 H_k^{(2)}}{k^2} - \sum^{\infty}_{k=1} \frac{H_k}{k^2} \left( \sum^{k}_{i=1} \frac{H_i}{i^2} \right) \nonumber \\ &=&
\frac{67}{48}\zeta(6) + \zeta(3)^2~.
\end{eqnarray}
Combining Eq.~(206) and Eq.~(207) it follows:
\begin{eqnarray}
\sum^{\infty}_{k=1} \frac{H_k}{k^2} \left( \sum^{k}_{i=1} \frac{H_{i-1}^{(2)}}{i(2i-1)} \right) = \frac{811}{96}\zeta(6) - \frac{57}{8}\zeta(3)^2 + \zeta(2)\sum^{\infty}_{k=1}\frac{h_k}{k^3}~.
\end{eqnarray}
Rearranging the summations in Eq.~(208) the nested double sum appearing in the last row of Eq.~(203) follows explicitly. It remains to calculate the second nested double sum. We start the calculation
as follows:
\begin{eqnarray}
\sum^{\infty}_{k=1} \frac{H_k^{(2)}}{k(2k-1)} \left( \sum^{k}_{i=1} \frac{h_i}{i^2 } \right) &=& \sum^{\infty}_{k=1} \frac{H_k^{(2)}}{k(2k-1)} \left( \sum^{\infty}_{i=1} \frac{h_i}{i^2} \right)
+ \sum^{\infty}_{k=1} \frac{H_k^{(2)} h_k}{k^3(2k-1)} - \sum^{\infty}_{k=1} \frac{h_k}{k^2} \left( \sum^{k}_{i=1}\frac{H_i^{(2)}}{i(2i-1)}\right)~.
\end{eqnarray}
Here only the last nested double sum is unknown. It follows by partial fraction decomposition:
\begin{eqnarray}
\sum^{\infty}_{k=1} \frac{h_k}{k^2} \left( \sum^{k}_{i=1}\frac{H_i^{(2)}}{i(2i-1)}\right) &=& 2\sum^{\infty}_{k=1} \frac{h_k}{k^2} \left( \sum^{k}_{i=1}\frac{H_i^{(2)}}{(2i-1)}\right) -
\sum^{\infty}_{k=1} \frac{h_k}{k^2} \left( \sum^{k}_{i=1}\frac{H_i^{(2)}}{i}\right) \nonumber \\ &=& 2\sum^{\infty}_{k=1} \frac{h_k}{k^2} \left( \sum^{k}_{i=1}\frac{H_i^{(2)}}{(2i-1)}\right) -
\sum^{\infty}_{k=1} \frac{h_k}{k^2} \left( H_k H_k^{(2)} + H_k^{(3)} - \sum^{k}_{i=1}\frac{H_i}{i^2} \right)~.
\end{eqnarray}
With help of Eq.~(177) and Eq.~(192) the last three sums on the right side of Eq.~(210) can be calculated explicitly. This way only the first nested double sum is unknown. It follows: 
\begin{eqnarray}
\sum^{\infty}_{k=1} \frac{h_k}{k^2} \left( \sum^{k}_{i=1}\frac{H_i^{(2)}}{(2i-1)}\right) &=& \sum^{\infty}_{k=1} \frac{h_k}{k^2} \left( \sum^{k}_{i=1}\frac{H_{i-1}^{(2)}}{(2i-1)}\right) +
\sum^{\infty}_{k=1} \frac{h_k}{k^2} \left( \sum^{k}_{i=1}\frac{1}{i^2(2i-1)}\right) \nonumber \\ &=&
\sum^{\infty}_{k=1} \frac{h_k}{k^2} \left( \sum^{k}_{i=1}\frac{1}{i^2(2i-1)}\right) + \sum^{\infty}_{k=1} \frac{H_k^{(2)} h_k^2}{k^2} - \sum^{\infty}_{k=1} \frac{h_k}{k^2} \left(\sum^{k}_{i=1} \frac{h_i}{i^2}\right)~.
\end{eqnarray}
From this the second nested double sum 
\begin{eqnarray}
\sum^{\infty}_{k=1} \frac{H_k^{(2)}}{k(2k-1)} \left( \sum^{k}_{i=1} \frac{h_i}{i^2 } \right) &=& \sum^{\infty}_{k=1} \frac{H_k^{(2)}}{k(2k-1)} \left( \sum^{\infty}_{i=1} \frac{h_i}{i^2} \right)
+ \sum^{\infty}_{k=1} \frac{H_k^{(2)} h_k}{k^3(2k-1)} - 2\sum^{\infty}_{k=1} \frac{h_k}{k^2} \left( \sum^{k}_{i=1}\frac{1}{i^2(2i-1)}\right) \nonumber \\ &+& \left(  \sum^{\infty}_{k=1} \frac{h_k}{k^2} \right)^2
+ \sum^{\infty}_{k=1} \frac{h_k^2}{k^4} - 2\sum^{\infty}_{k=1} \frac{H_k^{(2)} h_k^2}{k^2} + \sum^{\infty}_{k=1}\frac{H_k H_k^{(2)} h_k}{k^2} + \sum^{\infty}_{k=1}\frac{H_k^{(3)} h_k}{k^2} \nonumber \\ &-&
\sum^{\infty}_{k=1} \frac{h_k}{k^2} \left( \sum^{k}_{i=1}\frac{H_i^2}{i^2}\right)
\end{eqnarray}
follows as all other sums are known explicitly in Eq.~(212). This way the first member of the third family follows explicitly to:
\begin{eqnarray}
\sum^{\infty}_{k=1} \frac{H_k H_k^{(2)} H_k^{(2)}}{k(2k-1)} &=& -\frac{47}{16}\zeta(6) + \frac{31}{2}ln(2)\zeta(5) - \frac{31}{2}\zeta(5) - 5\left(ln(2)\right)^2\zeta(4) + 
10ln(2)\zeta(4) - 26\zeta(4) - \frac{37}{8}\zeta(3)^2 \nonumber \\ &+& 56ln(2)\zeta(3) - 40\zeta(3) - 16\left(ln(2)\right)^2\zeta(2) 
+ 32ln(2)\zeta(2) - 32\left(ln(2)\right)^2 + 64ln(2) \nonumber \\ &-& 8\sum^{\infty}_{k=1}\frac{h_k}{k^3} - \frac{1}{2}\sum^{\infty}_{k=1}\frac{h_k}{k^5}~.
\end{eqnarray}

The second member can be calculated explicitly by using the ansatz 
\begin{eqnarray}
\sum^{\infty}_{i=1} \frac{H_i}{2i-1}\left( \sum^{\infty}_{k=1} \frac{H_k^{(3)}}{k(i+k)} \right) = \sum^{\infty}_{k=1} \frac{H_k^{(3)}}{k} \left(\sum^{\infty}_{i=1} \frac{H_i}{(i+k)(2i-1)} \right)  
\end{eqnarray}
or
\begin{eqnarray}
\zeta(4)\sum^{\infty}_{i=1} \frac{H_i}{i(2i-1)} + \zeta(3)\sum^{\infty}_{i=1} \frac{H_i H_{i-1}}{i(2i-1)} - \zeta(2)\sum^{\infty}_{i=1} \frac{H_i H_{i-1}^{(2)}}{i(2i-1)} + \sum^{\infty}_{i=1} \frac{H_i}{i(2i-1)}
\left( \sum^{i-1}_{k=1} \frac{H_k}{k^3} \right) = \nonumber \\
\left( 4ln(2) - 2\left(ln(2)\right)^2 + \zeta(2) \right) \sum^{\infty}_{k=1} \frac{H_k^{(3)}}{k(2k+1)} - \sum^{\infty}_{k=1} \frac{H_k H_k^{(3)}}{k^2(2k+1)} +
\frac{1}{2} \sum^{\infty}_{k=1} \frac{H_k^{(2)} H_k^{(3)}}{k(2k+1)} + \frac{1}{2} \sum^{\infty}_{k=1} \frac{H_k^2 H_k^{(3)}}{k(2k+1)}~.
\end{eqnarray}
Only the last sum on the left side of Eq.(215) is unknown. It follows:
\begin{eqnarray}
\sum^{\infty}_{i=1} \frac{H_i}{i(2i-1)} \left( \sum^{i-1}_{k=1} \frac{H_k}{k^3} \right) &=& \sum^{\infty}_{i=1} \frac{H_i}{i(2i-1)} \left( \sum^{\infty}_{k=1} \frac{H_k}{k^3} \right) -
\sum^{\infty}_{k=1} \frac{H_k}{k^3} \left( \sum^{k}_{i=1} \frac{H_i}{i(2i-1)} \right)  \nonumber \\ &=&  \sum^{\infty}_{i=1} \frac{H_i}{i(2i-1)} \left( \sum^{\infty}_{k=1} \frac{H_k}{k^3} \right) -
2\sum^{\infty}_{k=1} \frac{H_k}{k^3} \left( \sum^{k}_{i=1} \frac{H_i}{2i-1} \right) - \sum^{\infty}_{k=1} \frac{H_k}{k^3} \left( \sum^{k}_{i=1} \frac{H_i}{i} \right)  \nonumber \\ &=&
\sum^{\infty}_{i=1} \frac{H_i}{i(2i-1)} \left( \sum^{\infty}_{k=1} \frac{H_k}{k^3} \right) + \frac{1}{2}\sum^{\infty}_{k=1} \frac{H_k H_k^{(2)}}{k^3} + \frac{1}{2}\sum^{\infty}_{k=1} \frac{H_k^3}{k^3}
- 4\sum^{\infty}_{k=1} \frac{H_k h_k}{k^3} + 2\sum^{\infty}_{k=1} \frac{H_k^2}{k^3} \nonumber \\ &-& 2\sum^{\infty}_{k=1} \frac{H_k}{k^3} \left( \sum^{k}_{i=1} \frac{H_{i-1}}{2i-1} \right)~.
\end{eqnarray}
Again only the last nested double sum on the right side of Eq.~(216) is unknown. The corresponding calculation starts with:
\begin{eqnarray}
\sum^{\infty}_{k=1} \frac{H_k}{k^2} \left(\sum^{\infty}_{i=1} \frac{h_i}{i(k+i)} \right) = \sum^{\infty}_{i=1} \frac{h_i}{i} \left( \sum^{\infty}_{k=1} \frac{H_k}{k^2(i+k)} \right)
\end{eqnarray}
or
\begin{eqnarray}
2ln(2) \sum^{\infty}_{k=1} \frac{H_k h_k}{k^3} + \sum^{\infty}_{k=1} \frac{H_k}{k^3} \left( \sum^{k}_{i=1} \frac{H_{i-1}}{2i-1} \right) = \sum^{\infty}_{k=1} \frac{h_k}{k^2} \left( \sum^{\infty}_{i=1} \frac{H_i}{i^2} \right) 
- \sum^{\infty}_{i=1} \frac{h_i}{i^2} \left( \sum^{\infty}_{k=1} \frac{H_k}{k(i+k)} \right)~.
\end{eqnarray}
From this we get:
\begin{eqnarray}
\sum^{\infty}_{k=1} \frac{H_k}{k^3} \left( \sum^{k}_{i=1} \frac{H_{i-1}}{2i-1} \right) &=&  \sum^{\infty}_{k=1} \frac{h_k}{k^2} \left( \sum^{\infty}_{i=1} \frac{H_i}{i^2} \right) -
2ln(2) \sum^{\infty}_{k=1} \frac{H_k h_k}{k^3} - \zeta(2)\sum^{\infty}_{k=1} \frac{h_k}{k^3} + \sum^{\infty}_{k=1} \frac{H_k h_k}{k^4} - \frac{1}{2}\sum^{\infty}_{k=1} \frac{H_k^{(2)} h_k}{k^3}
\nonumber \\ &-& \frac{1}{2}\sum^{\infty}_{k=1} \frac{H_k^2 h_k}{k^3}~.
\end{eqnarray}
Inserting this result in Eq.~(216) and inserting Eq.~(215) in Eq.~(215) the second member results after a lengthly calculational procedure.

\begin{eqnarray}
\sum^{\infty}_{k=1} \frac{H_k^2 H_k^{(3)}}{k(2k-1)} &=& \frac{1735}{16}\zeta(6) - \frac{279}{2}ln(2)\zeta(5) + \frac{279}{2}\zeta(5) + 53\left(ln(2)\right)^2\zeta(4) - 
106ln(2)\zeta(4) - 12\zeta(4) + \frac{49}{8}\zeta(3)^2 \nonumber \\ &+& 12ln(2)\zeta(2)\zeta(3) - 12\zeta(2)\zeta(3) - 16\left(ln(2)\right)^3\zeta(3)
+ 48\left(ln(2)\right)^2\zeta(3) + 8ln(2)\zeta(3) - 4\zeta(3)  \nonumber \\ &+& \frac{8}{3} \left(ln(2)\right)^4\zeta(2) - \frac{32}{3}\left(ln(2)\right)^3\zeta(2) - 16ln(2)\zeta(2) +
16\zeta(2) - \frac{8}{45}\left(ln(2)\right)^6 + \frac{16}{15}\left(ln(2)\right)^5 \nonumber \\ &+& \frac{64}{3}\left(ln(2)\right)^3 - 64\left(ln(2)\right)^2 + 64ln(2) - 8\sum^{\infty}_{k=1}\frac{h_k}{k^3}
+ \frac{9}{2}\sum^{\infty}_{k=1}\frac{h_k}{k^5} - 128Li_6\left(\frac{1}{2}\right) - 128Li_5\left(\frac{1}{2}\right)~.
\end{eqnarray}

\subsection{Fourth family}
The first member of the fourth family follows using the identity defined by Eq.~(140). The explicit calculation results in: 
\begin{eqnarray}
\sum^{\infty}_{k=1} \frac{h_k^3}{(2k-1)^3} &=& \frac{2823}{512}\zeta(6) - \frac{651}{128}ln(2)\zeta(5) + \frac{249}{64}\left(ln(2)\right)^2\zeta(4) +
\frac{133}{512}\zeta(3)^2 - \frac{9}{32}ln(2)\zeta(2)\zeta(3) + \frac{1}{8}\left(ln(2)\right)^4\zeta(2) \nonumber \\ &-& \frac{1}{120}\left(ln(2)\right)^6 + \frac{21}{128}\sum^{\infty}_{k=1}\frac{h_k}{k^5}
- 6Li_6\left(\frac{1}{2}\right)~.
\end{eqnarray}

As the second member of the fourth family we present the following identity:
\begin{eqnarray}
\sum^{\infty}_{k=1} \frac{h_k^2 h_k^{(2)}}{(2k-1)^2} &=& \frac{399}{1024}\zeta(6) + \frac{93}{128}ln(2)\zeta(5) + \frac{75}{64}\left(ln(2)\right)^2\zeta(4) + \frac{7}{512}\zeta(3)^2 + \frac{3}{16}ln(2)\zeta(2)\zeta(3)
\nonumber \\ &-& \frac{3}{128}\sum^{\infty}_{k=1}\frac{h_k}{k^5}~.
\end{eqnarray}
This result follows from the ansatz
\begin{eqnarray}
\sum^{\infty}_{k=1}\frac{H_k}{k} \left(\sum^{\infty}_{i=1}\frac{h_i^{(2)}}{(2i-1)(2i+2k-1)} \right) &=& \frac{7}{16}\zeta(3) \sum^{\infty}_{k=1} \frac{H_k}{k^2} + \frac{3}{16} \zeta(2)\sum^{\infty}_{k=1}
\frac{H_k H_{k-1}}{k^2} - \frac{1}{8}\sum^{\infty}_{k=1} \frac{H_k}{k^2} \left( \sum^{k-1}_{i=1} \frac{h_i}{i^2} \right)
\end{eqnarray}
or
\begin{eqnarray}
\sum^{\infty}_{k=1} \frac{H_k}{k^2} \left( \sum^{k-1}_{i=1} \frac{h_i}{i^2} \right)  &=& \frac{7}{2}\zeta(3) \sum^{\infty}_{k=1} \frac{H_k}{k^2} + \frac{3}{2} \zeta(2)\sum^{\infty}_{k=1}
\frac{H_k H_{k-1}}{k^2} -8 \sum^{\infty}_{i=1}\frac{h_i^{(2)}}{(2i-1)} \left( \sum^{\infty}_{k=1}\frac{H_k}{k(2i+2k-1} \right)
\end{eqnarray}
or
\begin{eqnarray}
\sum^{\infty}_{k=1} \frac{H_k}{k^2} \left( \sum^{k-1}_{i=1} \frac{h_i}{i^2} \right)  &=& \frac{7}{2}\zeta(3) \sum^{\infty}_{k=1} \frac{H_k}{k^2} + \frac{3}{2} \zeta(2)\sum^{\infty}_{k=1}
\frac{H_k H_{k-1}}{k^2} - 8\sum^{\infty}_{i=1}\frac{h_i^{(2)}}{(2i-1)} \Bigg( \frac{2\left(ln(2)\right)^2}{2i-1} + \frac{4ln(2)}{(2i-1)^2} \nonumber \\ &-& 4ln(2)\frac{h_i}{(2i-1)} - 4ln(2)\frac{h_i}{(2i-1)^2} +
\frac{2h_i^{(2)}}{(2i-1)} + \frac{2h_i^2}{(2i-1)}  \Bigg )~.
\end{eqnarray}
Comparing with Eq.~(177) we have an identity between the two cubic Euler sums 
\begin{eqnarray}
\sum^{\infty}_{i=1}\frac{H_i h_i^2}{i^3} \sim \sum^{\infty}_{i=1}\frac{h_i^2 h_i^{(2)}}{(2i-1)^2}
\end{eqnarray}
where Eq.~(222) follows explicitly if Eq.~(177) is used together with Eq.~(189).

\subsection{Fifth family}
The first and second members of the fifth family follow using the two identities defined by Eq.~(132) and Eq.~(133). The explicit calculations result in: 
\begin{eqnarray}
\sum^{\infty}_{k=1} \frac{H_k^2 h_k}{(2k-1)^3} &=& \frac{45}{16}\zeta(6) - \frac{155}{32}ln(2)\zeta(5) + \frac{279}{32}\zeta(5) - \frac{15}{4}\left(ln(2)\right)^2\zeta(4) + 
\frac{15}{8}ln(2)\zeta(4) - \frac{15}{4}\zeta(4)- \frac{105}{64}\zeta(3)^2 \nonumber \\ &+& \frac{21}{8}ln(2)\zeta(2)\zeta(3) - \frac{7}{2}\zeta(2)\zeta(3) + \frac{7}{2}\left(ln(2)\right)^3\zeta(3)
- 7\left(ln(2)\right)^2\zeta(3) + \frac{3}{4}\zeta(3) + 6\left(ln(2)\right)^2\zeta(2) \nonumber \\ &-& 10ln(2)\zeta(2) + 6\zeta(2) + \frac{1}{8}\zeta(2)\sum^{\infty}_{k=1}\frac{h_k}{k^3} -
\frac{1}{2}\left(ln(2)\right)^2\sum^{\infty}_{k=1}\frac{h_k}{k^3} + ln(2)\sum^{\infty}_{k=1}\frac{h_k}{k^3} + \frac{1}{2}\sum^{\infty}_{k=1}\frac{h_k}{k^3}  \nonumber \\ &+&
\frac{1}{8}\sum^{\infty}_{k=1}\frac{h_k}{k^5}
\end{eqnarray}
and
\begin{eqnarray}
\sum^{\infty}_{k=1} \frac{H_k h_k^2}{(2k-1)^3} &=& \frac{495}{512}\zeta(6) + \frac{31}{4}ln(2)\zeta(5) - \frac{217}{64}\zeta(5) - \frac{121}{32}\left(ln(2)\right)^2\zeta(4) + 
\frac{83}{16}ln(2)\zeta(4) - \frac{15}{16}\zeta(4) - \frac{133}{256}\zeta(3)^2 \nonumber \\ &-& \frac{25}{16}ln(2)\zeta(2)\zeta(3) - \frac{3}{16}\zeta(2)\zeta(3)
- \frac{7}{4}ln(2)\zeta(3) + \frac{3}{8}\zeta(3) - \frac{1}{3}\left(ln(2)\right)^4\zeta(2) + \frac{1}{3}\left(ln(2)\right)^3\zeta(2) \nonumber \\ &-& \frac{3}{2}\left(ln(2)\right)^2\zeta(2) +
\frac{3}{2}ln(2)\zeta(2) + \frac{1}{30}\left(ln(2)\right)^6 - \frac{1}{30}\left(ln(2)\right)^5 + \frac{1}{16}\zeta(2)\sum^{\infty}_{k=1}\frac{h_k}{k^3} +
\frac{1}{4}\sum^{\infty}_{k=1}\frac{h_k}{k^3}  \nonumber \\ &-& \frac{9}{64}\sum^{\infty}_{k=1}\frac{h_k}{k^5} - 4ln(2)Li_5\left(\frac{1}{2}\right) + 4Li_5\left(\frac{1}{2}\right)~.
\end{eqnarray}

In analogy to the first member of the second family the third member of the fifth family can be calculated explicitly by using the following ansatz:
\begin{eqnarray}
\sum^{\infty}_{i=1} \frac{h_i^{(2)}}{2i-1} \left( \sum^{\infty}_{k=1} \frac{H_k}{(2k+1)(i+k)} \right) &=& \left( \zeta(2) - 2\left(ln(2)\right)^2 \right) \sum^{\infty}_{i=1} \frac{h_i^{(2)}}{(2i-1)^2}
+ \frac{1}{2} \sum^{\infty}_{i=1} \frac{H_{i-1}^{(2)} h_i^{(2)}}{(2i-1)^2} + \frac{1}{2} \sum^{\infty}_{i=1} \frac{H_{i-1}^2 h_i^{(2)}}{(2i-1)^2}
\end{eqnarray}
or
\begin{eqnarray}
\sum^{\infty}_{k=1} \frac{H_k}{2k+1} \left( \sum^{\infty}_{i=1} \frac{h_i^{(2)}}{(2i-1)(k+i)} \right) &=& \left( \frac{21}{8}\zeta(3) - \frac{3}{2}ln(2)\zeta(2) \right) \sum^{\infty}_{k=1} \frac{H_k}{(2k+1)^2}
+ \frac{3}{2}\zeta(2)\sum^{\infty}_{k=1} \frac{H_k h_k}{(2k+1)^2} \nonumber \\ &-& 2ln(2) \sum^{\infty}_{k=1} \frac{H_k h_k^{(2)}}{(2k+1)^2} - \sum^{\infty}_{k=1} \frac{H_k}{(2k+1)^2} \left( \sum^{k}_{i=1}
\frac{H_{i-1}}{(2i-1)^2} \right)
\end{eqnarray}
with
\begin{eqnarray}
\sum^{\infty}_{k=1} \frac{H_k}{(2k+1)^2} \left( \sum^{k}_{i=1} \frac{H_{i-1}}{(2i-1)^2} \right) &=& \frac{1}{2}  \left( \sum^{\infty}_{k=1} \frac{H_{k-1}}{(2k-1)^2} \right)^2 - 
\frac{1}{2} \sum^{\infty}_{k=1} \frac{H_{k-1}^2}{(2k-1)^4}~. 
\end{eqnarray}
This way Eq.~(232) results to:

\begin{eqnarray}
\sum^{\infty}_{k=1} \frac{H_k^2 h_k^{(2)}}{(2k-1)^2} &=& \frac{765}{128}\zeta(6) - \frac{217}{16}ln(2)\zeta(5) + \frac{217}{16}\zeta(5) + \frac{75}{16}\left(ln(2)\right)^2\zeta(4) - 
\frac{75}{8}ln(2)\zeta(4) - \frac{15}{8}\zeta(4) - \frac{35}{32}\zeta(3)^2 \nonumber \\ &+& \frac{21}{8}ln(2)\zeta(2)\zeta(3) - \frac{21}{8}\zeta(2)\zeta(3) +
\frac{35}{2}ln(2)\zeta(3) - \frac{21}{2}\zeta(3) - 6\left(ln(2)\right)^2\zeta(2) + 6ln(2)\zeta(2) \nonumber \\ &-& 2\sum^{\infty}_{k=1}\frac{h_k}{k^3} + \frac{1}{4}\sum^{\infty}_{k=1}\frac{h_k}{k^5}~.
\end{eqnarray}

A similar approach allows to calculate explicitly the fourth member. It follows:
\begin{eqnarray}
\sum^{\infty}_{i=1} \frac{H_i}{i} \left( \sum^{\infty}_{k=1} \frac{H_k^{(2)}}{(2k-1)(2i+2k-1)} \right)  &=& \left( 4ln(2) - \zeta(2) \right) \sum^{\infty}_{i=1} \frac{H_i}{i(2i-1)} - 4ln(2)
\sum^{\infty}_{i=1} \frac{H_i}{i(2i-1)^2} + \frac{1}{2}\zeta(2)\sum^{\infty}_{i=1} \frac{H_i h_i}{i^2} \nonumber \\ &+& 2ln(2)\sum^{\infty}_{i=1} \frac{H_i h_i^{(2)}}{i^2} - 2\sum^{\infty}_{i=1} \frac{H_i}{i^2} \left(
\sum^{i-1}_{k=1} \frac{h_k}{(2k-1)^2} \right)
\end{eqnarray}
or
\begin{eqnarray}
\sum^{\infty}_{k=1} \frac{H_k^{(2)}}{2k-1} \left( \sum^{\infty}_{i=1} \frac{H_i}{i(2i+2k-1)} \right) &=& 2\left(ln(2)\right)^2\sum^{\infty}_{k=1} \frac{H_k^{(2)}}{(2k-1)^2} + 4ln(2)\sum^{\infty}_{k=1} \frac{H_k^{(2)}}{(2k-1)^3}
- 4\sum^{\infty}_{k=1} \frac{H_k^{(2)} h_k}{(2k-1)^3} \nonumber \\ &-& 4ln(2)\sum^{\infty}_{k=1} \frac{H_k^{(2)} h_k}{(2k-1)^2} + 2\sum^{\infty}_{k=1} \frac{H_k^{(2)} h_k^{(2)}}{(2k-1)^2} + 
2\sum^{\infty}_{k=1} \frac{H_k^{(2)} h_k^2}{(2k-1)^2}~.
\end{eqnarray}
This way Eq.~(237) follows if the nested double sum in Eq.~(233) is known explicitly. We make the ansatz:
\begin{eqnarray}
\sum^{\infty}_{i=1} \frac{H_i}{i} \left( \sum^{i}_{k=1} \frac{h_k^{(2)}}{k(i+k)} \right) &=& 2ln(2)\sum^{\infty}_{i=1} \frac{H_i h_i^{(2)}}{i^2} - \frac{3}{2}\zeta(2)\sum^{\infty}_{i=1} \frac{H_i h_i}{i^2}
+ \sum^{\infty}_{i=1} \frac{H_i h_i^{(3)}}{i^2} + \sum^{\infty}_{i=1} \frac{H_i h_i h_i^{(2)}}{i^2} + \frac{1}{2}\sum^{\infty}_{i=1} \frac{H_i^2 h_i^{(2)}}{i^2} \nonumber \\ &+& \zeta(2)\sum^{\infty}_{i=1} \frac{h_i^{(2)}}{i^2}
- \sum^{\infty}_{i=1} \frac{H_i h_i^{(2)}}{i^3} + \frac{1}{2} \sum^{\infty}_{i=1} \frac{H_i^{(2)} h_i^{(2)}}{i^2} + \frac{1}{2} \sum^{\infty}_{i=1} \frac{H_i^2 h_i^{(2)}}{i^2} \nonumber \\ &-&
2\sum^{\infty}_{i=1} \frac{H_i}{i^2}\left( \sum^{i}_{k=1} \frac{h_k}{(2k-1)^2} \right)
\end{eqnarray}
where the two-valued help function appears as the solution of the corresponding inhomogeneous first order difference equation. Rearranging the summations on the left side of Eq.~(235) we get:
\begin{eqnarray}
\sum^{\infty}_{k=1} \frac{H_k}{k^2}\left( \sum^{k}_{i=1} \frac{h_i}{(2i-1)^2} \right) &=& \frac{1}{2}\sum^{\infty}_{k=1} \frac{h_k h_k^{(2)}}{k^3} - \sum^{\infty}_{k=1} \frac{H_k h_k^{(2)}}{k^3} +
ln(2)\sum^{\infty}_{k=1} \frac{H_k h_k^{(2)}}{k^2} - \frac{3}{4}\zeta(2)\sum^{\infty}_{k=1} \frac{H_k h_k}{k^2} + \frac{1}{2}\sum^{\infty}_{k=1} \frac{H_k h_k^{(3)}}{k^2} \nonumber \\ &+& \frac{1}{2}
\sum^{\infty}_{k=1} \frac{H_k^{(2)} h_k^{(2)}}{k^2} + \frac{3}{4}\sum^{\infty}_{k=1} \frac{H_k^2 h_k^{(2)}}{k^2}~. 
\end{eqnarray}

This way Eq.~(237) results explicitly with the help of Eq.~(191) and Eq.~(193).
\begin{eqnarray}
\sum^{\infty}_{k=1} \frac{H_k^{(2)} h_k^2}{(2k-1)^2} &=& -\frac{5847}{128}\zeta(6) + \frac{279}{16}ln(2)\zeta(5) - \frac{121}{16}\left(ln(2)\right)^2\zeta(4) + \frac{75}{16}\zeta(4) - \frac{21}{128}\zeta(3)^2 +
\frac{63}{8}ln(2)\zeta(2)\zeta(3) \nonumber \\ &+& \frac{7}{2}ln(2)\zeta(3) - \frac{3}{2}\zeta(3) + \frac{1}{3}\left(ln(2)\right)^4\zeta(2) + 3\left(ln(2)\right)^2\zeta(2) - 6ln(2)\zeta(2) 
- \frac{1}{15}\left(ln(2)\right)^6 \nonumber \\ &-& \frac{3}{8}\zeta(2)\sum^{\infty}_{k=1}\frac{h_k}{k^3} - \frac{1}{2}\sum^{\infty}_{k=1}\frac{h_k}{k^3} - \frac{29}{32}\sum^{\infty}_{k=1}\frac{h_k}{k^5}
+ 48Li_6\left(\frac{1}{2}\right) + 16ln(2)Li_5\left(\frac{1}{2}\right)~.
\end{eqnarray}
In order to calculate the fifth member we start the computation with the following ansatz by use of a two-valued help function introduced in \cite{bra22} 
\begin{eqnarray}
\sum^{\infty}_{k=1} \frac{H_k}{(2k-1)^2}\left(\sum^{k}_{i=1} \frac{H_i}{(2i+2k-1)^2} \right) &=& \sum^{\infty}_{k=1} \frac{H_k}{(2k-1)^2}\left(\sum^{\infty}_{i=1} \frac{H_{i-1}}{(2i-1)^2} \right)
+ \frac{3}{2}\zeta(2)\sum^{\infty}_{k=1} \frac{H_k h_{k-1}}{(2k-1)^2} + 2ln(2)\sum^{\infty}_{k=1} \frac{H_k h_{k-1}^{(2)}}{(2k-1)^2} \nonumber \\ &-& 2\sum^{\infty}_{k=1} \frac{H_k h_{k-1}^{(3)}}{(2k-1)^2} -
\sum^{\infty}_{k=1} \frac{H_k h_{k-1} h_{k-1}^{(2)}}{(2k-1)^2}
\end{eqnarray}
or
\begin{eqnarray}
\sum^{\infty}_{i=1} H_i \left(\sum^{k}_{i=1} \frac{H_k}{(2k-1)^2(2i+2k-1)^2} \right) &=& \frac{1}{4}\sum^{\infty}_{i=1} \frac{H_i}{i^2}\left(\sum^{\infty}_{k=1} \frac{H_{k}}{(2k-1)^2} \right)
- \frac{1}{2}\sum^{\infty}_{k=1} \frac{H_i}{i^2} \left(\sum^{\infty}_{k=1} \frac{H_{k}}{(2k-1)(2i+2k-1)} \right)  \nonumber \\ &+& \frac{1}{4} \sum^{\infty}_{i=1}\frac{H_i}{i^2}\left( \sum^{\infty}_{k=1}
\frac{H_k}{(2i+2k-1)^2} \right)~. 
\end{eqnarray}

From this it follows:
\begin{eqnarray}
\sum^{\infty}_{k=1} \frac{H_k h_{k-1} h_{k-1}^{(2)}}{(2k-1)^2} &=& \frac{1}{2} \sum^{\infty}_{k=1} \frac{H_{k}}{(2k-1)^2} \left( \sum^{\infty}_{i=1} \frac{H_{i}}{(2i-1)^2} \right) +
\frac{3}{4}\zeta(2)\sum^{\infty}_{k=1} \frac{H_k h_{k-1}}{(2k-1)^2}+ ln(2)\sum^{\infty}_{k=1} \frac{H_k h_{k-1}^{(2)}}{(2k-1)^2} - \sum^{\infty}_{k=1} \frac{H_k h_{k-1}^{(3)}}{(2k-1)^2}
\nonumber \\ &-& \frac{1}{8}\sum^{\infty}_{i=1} \frac{H_i}{i^2}\left(\sum^{\infty}_{k=1} \frac{H_{k}}{(2k-1)^2} \right) + \frac{1}{2}ln(2) \sum^{\infty}_{k=1} \frac{H_{k}}{k^2(2k-1)} - \frac{1}{4}ln(2)
\sum^{\infty}_{k=1} \frac{H_{k} h_k}{k^3} - \frac{1}{4}\sum^{\infty}_{k=1} \frac{H_{k} h_k}{k^3(2k-1)} \nonumber \\ &-& \frac{1}{8}\sum^{\infty}_{k=1} \frac{H_{k}}{k^2}
\left( \sum^{\infty}_{i=1} \frac{H_{i-1}}{(2i-1)^2} \right) + \frac{1}{8} \sum^{\infty}_{k=1} \frac{H_{k} h_k^{(2)}}{k^3} + \frac{1}{8} \sum^{\infty}_{k=1} \frac{H_{k} h_k^2}{k^3} -
\frac{3}{16}\zeta(2)\sum^{\infty}_{k=1} \frac{H_k h_{k-1}}{k^2} \nonumber \\ &-& \frac{1}{4}ln(2)\sum^{\infty}_{k=1} \frac{H_k h_{k-1}^{(2)}}{k^2}
+ \frac{1}{4}\sum^{\infty}_{k=1} \frac{H_k h_{k-1}^{(3)}}{k^2} + \frac{1}{4}\sum^{\infty}_{k=1} \frac{H_k h_{k-1} h_{k-1}^{(2)}}{k^2}~.
\end{eqnarray}

All Euler sums are known explicitly and the last cubic Euler sum can be taken from Eq.~(193). This way Eq.~(241) results again after intensive algebraic manipulations.
\begin{eqnarray}
\sum^{\infty}_{k=1} \frac{H_k h_k h_k^{(2)}}{(2k-1)^2} &=& \frac{1935}{1024}\zeta(6) + \frac{341}{128}ln(2)\zeta(5) + \frac{93}{128}\zeta(5) - \frac{75}{32}\left(ln(2)\right)^2\zeta(4) + 
\frac{75}{32}ln(2)\zeta(4) - \frac{105}{64}\zeta(4) - \frac{119}{256}\zeta(3)^2 \nonumber \\ &-& \frac{27}{32}ln(2)\zeta(2)\zeta(3) + \frac{3}{16}\zeta(2)\zeta(3)
- \frac{7}{8}ln(2)\zeta(3) + \frac{3}{4}\left(ln(2)\right)^2\zeta(2) + \frac{3}{32}\zeta(2)\sum^{\infty}_{k=1}\frac{h_k}{k^3} +
\frac{1}{8}\sum^{\infty}_{k=1}\frac{h_k}{k^3} \nonumber \\ &-& \frac{7}{64}\sum^{\infty}_{k=1}\frac{h_k}{k^5}~.
\end{eqnarray}

The last member of the fifth family follows from the following ansatz:
\begin{eqnarray}
\sum^{\infty}_{k=1} \frac{H_k^{(2)}}{2k+1} \left(\sum^{\infty}_{i=1} \frac{h_i}{(2i-1)(k+i)} \right) = \sum^{\infty}_{i=1} \frac{h_i}{(2i-1)}
\left( \sum^{\infty}_{k=1} \frac{H_k^{(2)}}{(2k+1)(k+i)} \right)
\end{eqnarray}
or
\begin{eqnarray}
&& 2ln(2)\sum^{\infty}_{k=1} \frac{H_{k-1}^{(2)} h_{k-1}}{(2k-1)^2} + \zeta(2)\sum^{\infty}_{k=1} \frac{H_{k-1}^{(2)}}{(2k-1)^2} + \sum^{\infty}_{k=1} \frac{H_{k-1} H_{k-1}^{(2)} h_{k-1}}{(2k-1)^2}
- \sum^{\infty}_{k=1} \frac{H_{k-1}^{(2)}}{(2k-1)^2} \left( \sum^{k-1}_{i=1} \frac{h_{i}}{i} \right) = \nonumber \\ && \left( 2ln(2)\zeta(2) - \frac{3}{2}\zeta(3)\right)
\sum^{\infty}_{i=1} \frac{h_i}{(2i-1)^2} - \zeta(2) \sum^{\infty}_{i=1} \frac{h_i}{i(2i-1)^2} + \zeta(2)\sum^{\infty}_{i=1} \frac{H_i h_i}{(2i-1)^2} - \sum^{\infty}_{i=1} \frac{h_i}{(2i-1)^2}
\left( \sum^{i-1}_{k=1} \frac{H_{k}}{k^2} \right)~.
\end{eqnarray}
Both nested double sums are known explicitly from Eq.~(282) and Eq.~(236), with
\begin{eqnarray}
\sum^{\infty}_{i=1} \frac{h_i}{(2i-1)^2} \left( \sum^{i-1}_{k=1} \frac{H_{k}}{k^2} \right) = \sum^{\infty}_{i=1} \frac{h_i}{(2i-1)^2} \left( \sum^{\infty}_{k=1} \frac{H_{k}}{k^2} \right) -
\sum^{\infty}_{k=1} \frac{H_{k}}{k^2} \left( \sum^{k}_{i=1} \frac{h_i}{(2i-1)^2} \right)~.
\end{eqnarray}
This way last member of the fifth family results explicitly to:

\begin{eqnarray}
\sum^{\infty}_{k=1} \frac{H_k H_k^{(2)} h_k}{(2k-1)^2} &=& \frac{9789}{256}\zeta(6) - \frac{1147}{32}ln(2)\zeta(5) + \frac{279}{16}\zeta(5) + \frac{87}{4}\left(ln(2)\right)^2\zeta(4) - 
\frac{227}{8}ln(2)\zeta(4) + \frac{185}{16}\zeta(4) + \frac{889}{128}\zeta(3)^2 \nonumber \\ &-& \frac{49}{8}ln(2)\zeta(2)\zeta(3) + \frac{63}{8}\zeta(2)\zeta(3) -
7\left(ln(2)\right)^3\zeta(3) + 14\left(ln(2)\right)^2\zeta(3) - \frac{3}{2}\zeta(3) + \left(ln(2)\right)^4\zeta(2) \nonumber \\ &-& \frac{8}{3}\left(ln(2)\right)^3\zeta(2) - 6\left(ln(2)\right)^2\zeta(2) +
14ln(2)\zeta(2) - 12\zeta(2) - \frac{1}{15}\left(ln(2)\right)^6 + \frac{4}{15}\left(ln(2)\right)^5 - \frac{1}{4}\zeta(2)\sum^{\infty}_{k=1}\frac{h_k}{k^3}
\nonumber \\ &+& \left(ln(2)\right)^2\sum^{\infty}_{k=1}\frac{h_k}{k^3} - 2ln(2)\sum^{\infty}_{k=1}\frac{h_k}{k^3} - \frac{3}{2}\sum^{\infty}_{k=1}\frac{h_k}{k^3} + \frac{45}{32}\sum^{\infty}_{k=1}\frac{h_k}{k^5}
- 48Li_6\left(\frac{1}{2}\right) - 32Li_5\left(\frac{1}{2}\right)~.
\end{eqnarray}

\subsection{Sixth family}

The first member of the sixth family can be calculated by use of the following ansatz:
\begin{eqnarray}
\sum^{\infty}_{k=1} \frac{H_k}{k} \left( \sum^{\infty}_{i=1} \frac{h_{i}^{(3)}}{i(2k+2i-1)} \right) &=& \sum^{\infty}_{i=1} \frac{h_{i}^{(3)}}{i} \left( \sum^{\infty}_{k=1} \frac{H_k}{k(2k+2i-1)} \right)
\end{eqnarray}
or
\begin{eqnarray}
&&2\left(ln(2)\right)^2 \sum^{\infty}_{k=1} \frac{h_{k}^{(3)}}{k(2k-1)} - 4\sum^{\infty}_{k=1} \frac{h_k h_{k}^{(3)}}{k(2k-1)^2} - 4ln(2)\sum^{\infty}_{k=1} \frac{h_k h_{k}^{(3)}}{k(2k-1)} + 4ln(2)
\sum^{\infty}_{k=1} \frac{h_{k}^{(3)}}{k(2k-1)^2} + 2\sum^{\infty}_{k=1} \frac{h_k^{(2)} h_{k}^{(3)}}{k(2k-1)} + \nonumber \\ && 2\sum^{\infty}_{k=1} \frac{h_k^2 h_{k}^{(3)}}{k(2k-1)}  =
\sum^{\infty}_{k=1} \frac{H_k}{k(2k-1)} \left( \sum^{\infty}_{i=1} \frac{h_{i}^{(3)}}{i(2i+1)} \right) + \frac{7}{8}\zeta(3) \sum^{\infty}_{k=1} \frac{H_k H_{k-1}}{k(2k-1)} - \frac{3}{8}\zeta(2)
\sum^{\infty}_{k=1} \frac{H_k H_{k-1}^{(2)}}{k(2k-1)} + \nonumber \\ && \frac{1}{4} \sum^{\infty}_{k=1} \frac{H_k}{k(2k-1)} \left( \sum^{k-1}_{i=1} \frac{h_i}{i^3} \right)~.
\end{eqnarray}
Here only the last nested double sum is unknown. For an explicit calculation of these sum we use the ansatz:
\begin{eqnarray}
&&\sum^{\infty}_{k=1} \frac{H_k}{k(2k-1)} \left( \sum^{k}_{i=1} \frac{h_i}{i^3} \right) = \sum^{\infty}_{k=1} \frac{H_k}{k(2k-1)} \left( \sum^{\infty}_{i=1} \frac{h_i}{i^3} \right) +
\sum^{\infty}_{k=1} \frac{H_k h_k}{k^4(2k-1)} - \sum^{\infty}_{k=1} \frac{h_k}{k^3} \left( \sum^{k}_{i=1} \frac{H_i}{i(2i-1)} \right) =  \nonumber \\ && \sum^{\infty}_{k=1} \frac{H_k}{k(2k-1)}
\left( \sum^{\infty}_{i=1} \frac{h_i}{i^3} \right) + \sum^{\infty}_{k=1} \frac{H_k h_k}{k^4(2k-1)} - 2\sum^{\infty}_{k=1} \frac{H_k h_k^2}{k^3} + 2\sum^{\infty}_{k=1} \frac{h_k}{k^3}
\left( \sum^{k}_{i=1} \frac{h_i}{i} \right) - 4\sum^{\infty}_{k=1} \frac{h_k^2}{k^3} + 2\sum^{\infty}_{k=1} \frac{H_k h_k}{k^3} \nonumber \\ &&
\frac{1}{2}\sum^{\infty}_{k=1} \frac{H_k^{(2)} h_k}{k^3} + \frac{1}{2}\sum^{\infty}_{k=1} \frac{H_k^2 h_k}{k^3}~.
\end{eqnarray}
The nested double sum is defined by Eq.~(149). All other sums are known explicitly. Therefore it follows:
\begin{eqnarray}
\sum^{\infty}_{k=1} \frac{H_k}{k(2k-1)} \left( \sum^{k}_{i=1} \frac{h_i}{i^3} \right) &=& -\frac{1999}{64}\zeta(6) - \frac{31}{4}ln(2)\zeta(5) + \frac{31}{4}\zeta(5) + \frac{53}{4}\left(ln(2)\right)^2\zeta(4) - 
\frac{45}{4}\zeta(4) - \frac{77}{32}\zeta(3)^2 + 7ln(2)\zeta(2)\zeta(3) \nonumber \\ &-& 7\zeta(2)\zeta(3) - \frac{28}{3}\left(ln(2)\right)^3\zeta(3) + 10\zeta(3) + 2\left(ln(2)\right)^4\zeta(2) -
16ln(2)\zeta(2) + 16\zeta(2) - \frac{2}{9}\left(ln(2)\right)^6 \nonumber \\ &+& 4ln(2)\sum^{\infty}_{k=1}\frac{h_k}{k^3} - 2\sum^{\infty}_{k=1}\frac{h_k}{k^3} - \frac{9}{8}\sum^{\infty}_{k=1}\frac{h_k}{k^5}
+ 32Li_6\left(\frac{1}{2}\right) +  32ln(2)Li_5\left(\frac{1}{2}\right)~. 
\end{eqnarray}
This way the first member of the sixth family follows to:
\begin{eqnarray}
\sum^{\infty}_{k=1} \frac{h_k^2 h_k^{(3)}}{k(2k-1)} &=& -\frac{1567}{512}\zeta(6) + \frac{279}{64}ln(2)\zeta(5) - \frac{23}{32}\left(ln(2)\right)^2\zeta(4) - \frac{7}{256}\zeta(3)^2 +
\frac{9}{16}ln(2)\zeta(2)\zeta(3) \nonumber \\ &-& \frac{1}{12}\left(ln(2)\right)^4\zeta(2) + \frac{1}{180}\left(ln(2)\right)^6 - \frac{9}{64}\sum^{\infty}_{k=1}\frac{h_k}{k^5} + 4Li_6\left(\frac{1}{2}\right)~.
\end{eqnarray}

In order to calculate the second member of the sixth family we start with the ansatz:
\begin{eqnarray}
\sum^{\infty}_{k=1} \frac{H_k}{k(2k-1)} \left(\sum^{\infty}_{i=1} \frac{h_{i}^{(2)}}{(2i+2k-1)^2} \right) &=& \sum^{\infty}_{i=1} h_{i}^{(2)} \left(\sum^{\infty}_{k=1} \frac{H_k}{k(2k-1)(2i+2k-1)^2} \right)  
\end{eqnarray}
or
\begin{eqnarray}
&&\sum^{\infty}_{k=1} \frac{H_k}{k(2k-1)} \left(\sum^{\infty}_{i=1} \frac{h_{i-1}^{(2)}}{(2i-1)^2} \right) - \frac{3}{8}\zeta(2)\sum^{\infty}_{k=1} \frac{H_k H_{k-1}^{(2)}}{k(2k-1)} + \frac{1}{4}
\sum^{\infty}_{k=1} \frac{H_k}{k(2k-1)} \left(\sum^{k-1}_{i=1} \frac{h_i}{i^3} \right) + \frac{1}{4}\sum^{\infty}_{k=1} \frac{H_k}{k(2k-1)} \left(\sum^{k-1}_{i=1} \frac{h_i^{(2)}}{i^2} \right) = \nonumber \\ &&
2ln(2)\sum^{\infty}_{k=1} \frac{h_k^{(2)}}{k(2k-1)} - ln(2)\sum^{\infty}_{k=1} \frac{h_k h_k^{(2)}}{k^2} - \sum^{\infty}_{k=1} \frac{h_k h_k^{(2)}}{k^2(2k-1)} + \frac{1}{2}\sum^{\infty}_{k=1}
\frac{h_k^{(2)} h_k^{(2)}}{k^2} + \frac{1}{2}\sum^{\infty}_{k=1} \frac{h_k^2 h_k^{(2)}}{k^2} - 2\left(ln(2)\right)^2\sum^{\infty}_{k=1} \frac{h_k^{(2)}}{(2k-1)^2} -  \nonumber \\ &&
4ln(2)\sum^{\infty}_{k=1} \frac{h_k^{(2)}}{(2k-1)^3} + 4ln(2)\sum^{\infty}_{k=1} \frac{h_k h_k^{(2)}}{(2k-1)^2} + 4\sum^{\infty}_{k=1} \frac{h_k h_k^{(2)}}{(2k-1)^3} - 2\sum^{\infty}_{k=1}
\frac{h_k^{(2)} h_k^{(2)}}{(2k-1)^2} - 2\sum^{\infty}_{k=1} \frac{h_k^2 h_k^{(2)}}{(2k-1)^2} + \nonumber \\ && \sum^{\infty}_{k=1} \frac{H_{k-1}}{(2k-1)^2} \left(\sum^{\infty}_{i=1} \frac{h_i^{(2)}}{i(2i-1)} \right)
+ \frac{3}{2}\zeta(2) \sum^{\infty}_{k=1} \frac{h_{k-1} h_k^{(2)}}{k(2k-1)} - 2ln(2)\sum^{\infty}_{k=1} \frac{h_{k-1}^{(2)} h_k^{(2)}}{k(2k-1)} - 2\sum^{\infty}_{k=1} \frac{h_{k-1}^{(3)} h_k^{(2)}}{k(2k-1)} -
\nonumber \\ && 2\sum^{\infty}_{k=1} \frac{h_{k-1} h_{k-1}^{(2)} h_k^{(2)}}{k(2k-1)}~.
\end{eqnarray}
Only the last nested double sum in the first row of Eq.~(252) is unknown. We start the corresponding calculation using the following ansatz:
\begin{eqnarray}
\sum^{\infty}_{k=1} \frac{H_k}{k(2k-1)} \left( \sum^{k-1}_{i=1} \frac{h_i^{(2)}}{i^2} \right) = \sum^{\infty}_{k=1} \frac{H_k}{k(2k-1)} \left(\sum^{\infty}_{i=1} \frac{h_i^{(2)}}{i^2} \right) -
\sum^{\infty}_{i=1} \frac{h_i^{(2)}}{i^2} \left( \sum^{\infty}_{k=1} \frac{H_k}{k(2k-1)} \right)~.   
\end{eqnarray}
From this we get by partial fraction decomposition of the last sum:
\begin{eqnarray}
\sum^{\infty}_{k=1} \frac{H_k}{k(2k-1)} \left( \sum^{k-1}_{i=1} \frac{h_i^{(2)}}{i^2} \right) &=& \sum^{\infty}_{k=1} \frac{H_k}{k(2k-1)} \left(\sum^{\infty}_{i=1} \frac{h_i^{(2)}}{i^2} \right)
- 2\sum^{\infty}_{k=1} \frac{H_k h_k h_k^{(2)}}{k^2} + \frac{1}{2}\sum^{\infty}_{k=1} \frac{H_k^{(2)} h_k^{(2)}}{k^2} + \frac{1}{2}\sum^{\infty}_{k=1} \frac{H_k^2 h_k^{(2)}}{k^2} \nonumber \\ &-& 4
\sum^{\infty}_{k=1} \frac{h_k^{(2)} h_k^{(2)}}{k^2} + 2\sum^{\infty}_{k=1} \frac{H_k h_k^{(2)}}{k^2} + 2\sum^{\infty}_{k=1} \frac{h_k^{(2)}}{k^2} \left( \sum^{k}_{i=1} \frac{h_i}{i} \right)~.
\end{eqnarray}
The last sum in Eq.~(254) is defined by Eq.~(281). All other sums are known explicitly. This way it results 
\begin{eqnarray}
\sum^{\infty}_{k=1} \frac{H_k}{k(2k-1)} \left(\sum^{k}_{i=1} \frac{h_{i}^{(2)}}{i^2} \right) &=& \frac{47}{2}\zeta(6) + \frac{31}{2}ln(2)\zeta(5) - \frac{31}{2}\zeta(5) - \frac{151}{8}\left(ln(2)\right)^2\zeta(4) + 
\frac{45}{4}ln(2)\zeta(4) + \frac{15}{2}\zeta(4) - \frac{7}{64}\zeta(3)^2 \nonumber \\ &-& \frac{7}{2}ln(2)\zeta(2)\zeta(3) + \frac{7}{2}\zeta(2)\zeta(3) + \frac{28}{3}\left(ln(2)\right)^3\zeta(3) - 35ln(2)\zeta(3) + 
21\zeta(3) - 2\left(ln(2)\right)^4\zeta(2) \nonumber \\ &+& 12\left(ln(2)\right)^2\zeta(2) - 12ln(2)\zeta(2) + \frac{2}{9}\left(ln(2)\right)^6 + \frac{3}{4}\zeta(2)\sum^{\infty}_{k=1}\frac{h_k}{k^3} 
- 4ln(2)\sum^{\infty}_{k=1}\frac{h_k}{k^3} + 6\sum^{\infty}_{k=1}\frac{h_k}{k^3} \nonumber \\ &+& \frac{19}{16}\sum^{\infty}_{k=1}\frac{h_k}{k^5} - 32Li_6\left(\frac{1}{2}\right) - 32ln(2)Li_5\left(\frac{1}{2}\right)
\end{eqnarray}
and we get for the second member of the sixth family: 
\begin{eqnarray}
\sum^{\infty}_{k=1} \frac{h_k h_k^{(2)} h_k^{(2)}}{k(2k-1)} &=& \frac{831}{512}\zeta(6) - \frac{31}{64}ln(2)\zeta(5) - \frac{45}{32}\left(ln(2)\right)^2\zeta(4) + \frac{91}{256}\zeta(3)^2 +
\frac{9}{16}ln(2)\zeta(2)\zeta(3) + \frac{1}{64}\sum^{\infty}_{k=1}\frac{h_k}{k^5}~.
\end{eqnarray}

\subsection{Seventh family}
In order to compute explicitly the first member of the seventh family, following Adegoke \cite{ade16}, we use the ansatz:
\begin{eqnarray}
\sum^{\infty}_{k=1} \frac{H_k h_k^{(2)} h_k^{(2)}}{k(2k-1)} &=& 2\sum^{\infty}_{k=1} \frac{H_k}{k(2k-1)} \left( \sum^{k}_{i=1} \frac{h_{i}^{(2)}}{(2i-1)^2} \right) - \sum^{\infty}_{k=1}
\frac{H_k h_k^{(4)}}{k(2k-1)} = \nonumber 2\sum^{\infty}_{k=1} \frac{H_k}{k(2k-1)} \left( \sum^{\infty}_{i=1} \frac{h_{i}^{(2)}}{(2i-1)^2} \right) \nonumber \\ &+& 2\sum^{\infty}_{k=1}
\frac{H_k h_k^{(2)}}{k(2k-1)^3} - \sum^{\infty}_{k=1} \frac{H_k h_k^{(4)}}{k(2k-1)} - 2\sum^{\infty}_{k=1} \frac{h_{k}^{(2)}}{(2k-1)^2} \left( \sum^{k}_{i=1} \frac{H_i}{i(2i-1)} \right)~. 
\end{eqnarray}
All Euler sums are known explicitly, except the last nested double sum needs an explicit computation. We start with:
\begin{eqnarray}
\sum^{\infty}_{k=1} \frac{h_{k}^{(2)}}{(2k-1)^2} \left( \sum^{k}_{i=1} \frac{H_i}{i(2i-1)} \right) &=& 2\sum^{\infty}_{k=1} \frac{h_{k}^{(2)}}{(2k-1)^2} \left( \sum^{k}_{i=1} \frac{H_{i-1}}{(2i-1)} \right)-
\sum^{\infty}_{k=1} \frac{h_{k}^{(2)}}{(2k-1)^2} \left( \sum^{k}_{i=1} \frac{H_i}{i}  \right) \nonumber \\ &-& 2\sum^{\infty}_{k=1} \frac{h_{k}^{(2)}}{(2k-1)^2} \left( \sum^{k}_{i=1} \frac{1}{i(2i-1)} \right)
\end{eqnarray}
or
\begin{eqnarray}
\sum^{\infty}_{k=1} \frac{h_{k}^{(2)}}{(2k-1)^2} \left( \sum^{k}_{i=1} \frac{H_i}{i(2i-1)} \right) &=& 2\sum^{\infty}_{k=1} \frac{H_k h_k h_{k}^{(2)}}{(2k-1)^2} - 2\sum^{\infty}_{k=1} \frac{h_{k}^{(2)}}{(2k-1)^2}
\left( \sum^{k}_{i=1} \frac{h_i}{i} \right) - \frac{1}{2}\sum^{\infty}_{k=1} \frac{H_{k}^{(2)} h_{k}^{(2)}}{(2k-1)^2} \nonumber \\ &-& \frac{1}{2}\sum^{\infty}_{k=1} \frac{H_{k}^2 h_{k}^{(2)}}{(2k-1)^2}
- 4\sum^{\infty}_{k=1} \frac{h_k h_{k}^{(2)}}{(2k-1)^2} + 2\sum^{\infty}_{k=1} \frac{H_k h_{k}^{(2)}}{(2k-1)^2}~.
\end{eqnarray}
All Euler sums are known where two members of the sixth family appear. Therefore, only the nested double sum is unknown. The corresponding explicit calculation starts with:
\begin{eqnarray}
\sum^{\infty}_{i=1} \frac{h_{i}^{(2)}}{(2i-1)} \left( \sum^{\infty}_{k=1} \frac{h_k}{k(2i+2k-1)} \right) &=& \frac{1}{2}\zeta(2)\sum^{\infty}_{i=1} \frac{h_{i}^{(2)}}{(2i-1)^2} +
\sum^{\infty}_{i=1} \frac{h_{i}^{(2)}}{(2i-1)^2} \left( \sum^{i-1}_{k=1} \frac{h_k}{k} \right)
\end{eqnarray}
or
\begin{eqnarray}
\sum^{\infty}_{k=1} \frac{h_k}{k} \left( \sum^{\infty}_{i=1} \frac{h_{i}^{(2)}}{(2i-1)(2i+2k-1)} \right) &=& \frac{7}{16}\zeta(3)\sum^{\infty}_{k=1} \frac{h_k}{k^2} + \frac{3}{16}\zeta(2)
\sum^{\infty}_{k=1} \frac{H_{k-1} h_k}{k} - \frac{1}{8}\sum^{\infty}_{k=1} \frac{h_k}{k^2} \left( \sum^{k-1}_{i=1} \frac{h_i}{i^2} \right)~.
\end{eqnarray}
From this we get:
\begin{eqnarray}
\sum^{\infty}_{k=1} \frac{h_{k}^{(2)}}{(2k-1)^2} \left( \sum^{k}_{i=1} \frac{h_i}{i} \right) &=& \frac{525}{1024}\zeta(6) + \frac{93}{128}\zeta(5) + \frac{75}{32}ln(2)\zeta(4) - \frac{105}{64}\zeta(4)
+ \frac{49}{256}\zeta(3)^2 + \frac{3}{16}\zeta(2)\zeta(3) - \frac{7}{8}ln(2)\zeta(3) \nonumber \\ &+& \frac{3}{4}\left(ln(2)\right)^2\zeta(2) - \frac{3}{32}\zeta(2)\sum^{\infty}_{k=1}\frac{h_k}{k^3}
+ \frac{1}{8}\sum^{\infty}_{k=1}\frac{h_k}{k^3}
\end{eqnarray}
and this way Eq.~(263) follows to:
\begin{eqnarray}
\sum^{\infty}_{k=1} \frac{H_k h_k^{(2)} h_k^{(2)}}{k(2k-1)} &=& \frac{585}{256}\zeta(6) - \frac{341}{32}ln(2)\zeta(5) + \frac{155}{32}\zeta(5) + \frac{45}{8}\left(ln(2)\right)^2\zeta(4) 
-\frac{45}{8}ln(2)\zeta(4) + \frac{161}{64}\zeta(3)^2  \nonumber \\ &-& \frac{9}{8}ln(2)\zeta(2)\zeta(3) + \frac{9}{8}\zeta(2)\zeta(3) - \frac{3}{8}\zeta(2)\sum^{\infty}_{k=1}\frac{h_k}{k^3}
+ \frac{3}{16}\sum^{\infty}_{k=1}\frac{h_k}{k^5}~.
\end{eqnarray}

The calculational procedure concerning the second member starts with the ansatz:
\begin{eqnarray}
\sum^{\infty}_{i=1} \frac{H_i^{(2)}}{i} \left(\sum^{\infty}_{k=1} \frac{h_k}{k(i+k)} \right) = \sum^{\infty}_{k=1} \frac{h_k}{k} \left(\sum^{\infty}_{i=1} \frac{H_i^{(2)}}{i(k+k)} \right) \nonumber \\ 
\end{eqnarray}
or
\begin{eqnarray}
2ln(2)\sum^{\infty}_{i=1} \frac{H_i^{(2)} h_i}{i^2} + \sum^{\infty}_{i=1} \frac{H_i H_i^{(2)} h_i}{i^2} - \sum^{\infty}_{i=1} \frac{H_i^{(2)}}{i^2} \left(\sum^{i}_{k=1} \frac{h_k}{k} \right) =  \nonumber \\
\zeta(3)\sum^{\infty}_{k=1} \frac{h_k}{k^2} - \zeta(2)\sum^{\infty}_{k=1} \frac{h_k}{k^3} + \zeta(2)\sum^{\infty}_{k=1} \frac{H_k h_k}{k^2} - \sum^{\infty}_{k=1} \frac{h_k}{k^2}
\left(\sum^{k-1}_{i=1} \frac{H_i}{i^2} \right) = \nonumber \\
\zeta(3)\sum^{\infty}_{k=1} \frac{h_k}{k^2} - \zeta(2)\sum^{\infty}_{k=1} \frac{h_k}{k^3} + \zeta(2)\sum^{\infty}_{k=1} \frac{H_k h_k}{k^2} - \sum^{\infty}_{k=1} \frac{h_k}{k^2}
\left(\sum^{\infty}_{i=1} \frac{H_i}{i^2}\right) + \sum^{\infty}_{k=1} \frac{H_k}{k^2}\left(\sum^{k}_{i=1} \frac{h_i}{i^2}\right)~.
\end{eqnarray}
As the nested double sum in the last row of Eq.~(265) has been defined by Eq.~(177) the nested double in the first row of Eq.~(265) results because all other sums are known explicitly. It follows:
\begin{eqnarray}
\sum^{\infty}_{k=1} \frac{H_k^{(2)}}{k^2} \left(\sum^{k}_{i=1} \frac{h_i}{i} \right) &=& \frac{93}{8}ln(2)\zeta(5) - \frac{21}{16}\zeta(3)^2 - \frac{7}{2}ln(2)\zeta(2)\zeta(3) +
\zeta(2)\sum^{\infty}_{k=1}\frac{h_k}{k^3}~.
\end{eqnarray}
Now one can write:
\begin{eqnarray}
\sum^{\infty}_{k=1} \frac{H_k^{(2)}}{k^2} \left(\sum^{k}_{i=1} \frac{h_i}{i} \right) = 2\sum^{\infty}_{k=1} \frac{H_k^{(2)}}{k^2} \left( \sum^{k}_{i=1} \frac{h_i}{2i-1} \right) -
\sum^{\infty}_{k=1} \frac{H_k^{(2)}}{k^2} \left(\sum^{k}_{i=1} \frac{h_i}{i(2i-1)} \right) =  \nonumber \\ \sum^{\infty}_{k=1} \frac{H_k^{(2)} h_k^{(2)}}{k^2} + 
\sum^{\infty}_{k=1} \frac{H_k^{(2)} h_k^2}{k^2} - \sum^{\infty}_{k=1} \frac{H_k^{(2)}}{k^2} \left(\sum^{\infty}_{i=1} \frac{h_i}{i(2i-1)} \right) - \sum^{\infty}_{k=1} \frac{H_k^{(2)} h_k}{k^3(2k-1)}
+ \sum^{\infty}_{k=1} \frac{h_k}{k(2k-1)} \left(\sum^{k}_{i=1} \frac{H_i^{(2)}}{i^2} \right) =  \nonumber \\
\sum^{\infty}_{k=1} \frac{H_k^{(2)} h_k^{(2)}}{k^2} + \sum^{\infty}_{k=1} \frac{H_k^{(2)} h_k^2}{k^2} - \sum^{\infty}_{k=1} \frac{H_k^{(2)}}{k^2} \left(\sum^{\infty}_{i=1} \frac{h_i}{i(2i-1)} \right)
- \sum^{\infty}_{k=1} \frac{H_k^{(2)} h_k}{k^3(2k-1)} + \frac{1}{2}\sum^{\infty}_{k=1} \frac{H_k^{(4)} h_k}{k(2k-1)} +  \nonumber \\ \frac{1}{2}\sum^{\infty}_{k=1} \frac{H_k^{(2)} H_k^{(2)} h_k}{k(2k-1)}~.
\end{eqnarray}
This way after a variety of standard algebraic manipulations the second member follows explicitly:
\begin{eqnarray}
\sum^{\infty}_{k=1} \frac{H_k^{(2)} H_k^{(2)} h_k}{k(2k-1)} &=& \frac{185}{32}\zeta(6) + \frac{31}{4}ln(2)\zeta(5) - \frac{31}{4}\zeta(5) - 25\zeta(4)
- \frac{35}{4}\zeta(3)^2 + 28ln(2)\zeta(3) - 14\zeta(3) + 16\zeta(2) \nonumber \\ &+& \zeta(2)\sum^{\infty}_{k=1}\frac{h_k}{k^3} + 4\sum^{\infty}_{k=1}\frac{h_k}{k^3}
- \frac{3}{2}\sum^{\infty}_{k=1}\frac{h_k}{k^5}~. 
\end{eqnarray}

The explicit calculation of the third member is based on the following ansatz \cite{ade16}:
\begin{eqnarray}
\sum^{\infty}_{k=1} \frac{H_k}{k(2k-1)} \left( \sum^{k}_{i=1} \frac{h_i^{(2)}}{i^2} \right) &=& \sum^{\infty}_{k=1} \frac{H_k H_k^{(2)} h_k^{(2)}}{k(2k-1)} - \sum^{\infty}_{k=1} \frac{H_k}{k(2k-1)}
\left( \sum^{k}_{i=1} \frac{H_{i-1}^{(2)}}{(2i-1)^2} \right) \nonumber \\ &=& \sum^{\infty}_{k=1} \frac{H_k H_k^{(2)} h_k^{(2)}}{k(2k-1)} - \sum^{\infty}_{k=1} \frac{H_k}{k(2k-1)}
\left( \sum^{\infty}_{i=1} \frac{H_{i-1}^{(2)}}{(2i-1)^2} \right) - \sum^{\infty}_{k=1} \frac{H_k H_{k-1}^{(2)}}{k(2k-1)^3} \nonumber \\ &+& \sum^{\infty}_{k=1} \frac{H_{k-1}^{(2)}}{(2k-1)^2} 
\left( \sum^{k}_{i=1} \frac{H_i}{i(2i-1)} \right) \nonumber \\ &=& \sum^{\infty}_{k=1} \frac{H_k H_k^{(2)} h_k^{(2)}}{k(2k-1)} - \sum^{\infty}_{k=1} \frac{H_k}{k(2k-1)}
\left( \sum^{\infty}_{i=1} \frac{H_{i-1}^{(2)}}{(2i-1)^2} \right) - \sum^{\infty}_{k=1} \frac{H_k H_{k-1}^{(2)}}{k(2k-1)^3} \nonumber \\ &+& 4\sum^{\infty}_{k=1} \frac{H_{k-1}^{(2)} h_k}{(2k-1)^2}
- 2\sum^{\infty}_{k=1} \frac{H_{k-1}^{(2)} H_k}{(2k-1)^2} - \frac{1}{2}\sum^{\infty}_{k=1} \frac{H_k^{(2)} H_{k-1}^{(2)}}{(2k-1)^2} - \frac{1}{2}\sum^{\infty}_{k=1} \frac{H_k^2 H_{k-1}^{(2)}}{(2k-1)^2} 
\nonumber \\ &+& 2\sum^{\infty}_{k=1} \frac{H_{k-1}^{(2)}}{(2k-1)^2} \left( \sum^{k}_{i=1} \frac{H_{i-1}}{(2i-1)} \right)~. 
\end{eqnarray}
Using the last member of the fifth family it remains to calculate in Eq.~(269) only the first nested double sum on the left side and the last nested double sum on the right side. We calculate first the
last nested double sum on the right side of Eq.~(269). It follows 
\begin{eqnarray}
\sum^{\infty}_{k=1} \frac{H_{k-1}^{(2)}}{(2k-1)^2} \left( \sum^{k}_{i=1} \frac{H_{i-1}}{(2i-1)} \right) &=& \sum^{\infty}_{k=1} \frac{H_k H_{k-1}^{(2)} h_k}{(2k-1)^2} - \sum^{\infty}_{k=1} \frac{H_{k-1}^{(2)}}{(2k-1)^2}
\left( \sum^{k}_{i=1} \frac{h_i}{i} \right)~.
\end{eqnarray}
Again the explicitly known last member of the fifth family appears and a nested double sum which is known explicitly from Eq.~(286). This way the corresponding nested double sum is also known explicitly.
The first nested double sum on the left side of Eq.~(269) results from the ansatz:
\begin{eqnarray}
\sum^{\infty}_{k=1} \frac{H_k}{k(2k-1)} \left( \sum^{k}_{i=1} \frac{h_i^{(2)}}{i^2} \right) &=& \sum^{\infty}_{k=1} \frac{H_k}{k(2k-1)} \left( \sum^{\infty}_{i=1} \frac{h_i^{(2)}}{i^2} \right) +
\sum^{\infty}_{k=1} \frac{H_k h_k^{(2)}}{k^3(2k-1)} - 2\sum^{\infty}_{k=1} \frac{h_k^{(2)}}{k^2} \left(\sum^{k}_{i=1} \frac{H_i}{2i-1)} \right) +
\sum^{\infty}_{k=1} \frac{h_k^{(2)}}{k^2} \left(\sum^{k}_{i=1} \frac{H_i}{i} \right) \nonumber \\ &=& \sum^{\infty}_{k=1} \frac{H_k}{k(2k-1)} \left( \sum^{\infty}_{i=1} \frac{h_i^{(2)}}{i^2} \right) +
\sum^{\infty}_{k=1} \frac{H_k h_k^{(2)}}{k^3(2k-1)} + \frac{1}{2}\sum^{\infty}_{k=1} \frac{H_k^{(2)} h_k^{(2)}}{k^2} + \frac{1}{2}\sum^{\infty}_{k=1} \frac{H_k^2 h_k^{(2)}}{k^2} \nonumber \\ &-&
2\sum^{\infty}_{k=1} \frac{h_k^{(2)}}{k^2} \left(\sum^{k}_{i=1} \frac{1}{i(2i-1)} \right) - 2\sum^{\infty}_{k=1} \frac{H_k h_k h_k^{(2)}}{k^2} + 2\sum^{\infty}_{k=1} \frac{h_k^{(2)}}{k^2}
\left(\sum^{k}_{i=1} \frac{h_i}{i} \right)~.
\end{eqnarray}
As all sums including the last nested double sum (Eq.~(285)) are known explicitly the third member can be calculated. The result is presented in Eq.~(272). 

\begin{eqnarray}
\sum^{\infty}_{k=1} \frac{H_k H_k^{(2)} h_k^{(2)}}{k(2k-1)} &=& -\frac{1503}{16}\zeta(6) + \frac{341}{2}ln(2)\zeta(5) - \frac{2139}{16}\zeta(5) - \frac{605}{8}\left(ln(2)\right)^2\zeta(4) + 
\frac{189}{2}ln(2)\zeta(4) + \frac{15}{2}\zeta(4) - \frac{329}{64}\zeta(3)^2 \nonumber \\ &-& \frac{77}{4}ln(2)\zeta(2)\zeta(3) + \frac{63}{4}\zeta(2)\zeta(3) +
28\left(ln(2)\right)^3\zeta(3) - 42\left(ln(2)\right)^2\zeta(3) - 35ln(2)\zeta(3) + 21\zeta(3) \nonumber \\ &-& \frac{14}{3}\left(ln(2)\right)^4\zeta(2) + 8\left(ln(2)\right)^3\zeta(2) +
12\left(ln(2)\right)^2\zeta(2) - 12ln(2)\zeta(2) + \frac{2}{5}\left(ln(2)\right)^6 - \frac{4}{5}\left(ln(2)\right)^5 \nonumber \\ &+& \frac{9}{4}\zeta(2)\sum^{\infty}_{k=1}\frac{h_k}{k^3}
- 4\left(ln(2)\right)^2\sum^{\infty}_{k=1}\frac{h_k}{k^3} + 4ln(2)\sum^{\infty}_{k=1}\frac{h_k}{k^3} + 6\sum^{\infty}_{k=1}\frac{h_k}{k^3} - \frac{61}{16}\sum^{\infty}_{k=1}\frac{h_k}{k^5}
+ 96Li_6\left(\frac{1}{2}\right) \nonumber \\ &-& 32ln(2)Li_5\left(\frac{1}{2}\right) + 96Li_5\left(\frac{1}{2}\right)~.
\end{eqnarray}

In order to calculate explicitly the fourth member of the seventh family we start with the ansatz \cite{ade16}
\begin{eqnarray}
\sum^{\infty}_{k=1} \frac{H_k^{(2)} h_k h_k^{(2)}}{k(2k-1)} &=& \sum^{\infty}_{k=1} \frac{h_k}{k(2k-1)} \left( \sum^{k}_{i=1} \frac{h_i^{(2)}}{i^2} \right) + \sum^{\infty}_{k=1} \frac{h_k}{k(2k-1)}
\left( \sum^{k}_{i=1} \frac{H^{(2)}_{i-1}}{(2i-1)^2} \right) \nonumber \\ &=& \sum^{\infty}_{k=1} \frac{h_k}{k(2k-1)} \left( \sum^{k}_{i=1} \frac{h_i^{(2)}}{i^2} \right) + 
\sum^{\infty}_{k=1} \frac{H^{(2)}_{k-1}}{(2k-1)^2} \left( \sum^{\infty}_{i=1} \frac{h_i}{i(2i-1)} \right) + \sum^{\infty}_{k=1} \frac{H^{(2)}_{k-1} h_k}{k(2k-1)^3} \nonumber \\ &-&
\sum^{\infty}_{k=1} \frac{H^{(2)}_{k-1}}{(2k-1)^2} \left( \sum^{k}_{i=1} \frac{h_i}{i(2i-1)} \right) \nonumber \\ &=&
\sum^{\infty}_{k=1} \frac{h_k}{k(2k-1)} \left( \sum^{k}_{i=1} \frac{h_i^{(2)}}{i^2} \right) + \sum^{\infty}_{k=1} \frac{H^{(2)}_{k-1}}{k(2k-1)^2} \left( \sum^{\infty}_{i=1} \frac{h_i}{i(2i-1)} \right) +
\sum^{\infty}_{k=1} \frac{H^{(2)}_{k-1} h_k}{k(2k-1)^3} \nonumber \\ &-& \sum^{\infty}_{k=1} \frac{H^{(2)}_{k-1} h^{(2)}_k}{(2k-1)^2} - \sum^{\infty}_{k=1} \frac{H^{(2)}_{k-1} h^2_k}{(2k-1)^2} 
+ \sum^{\infty}_{k=1} \frac{H^{(2)}_{k-1}}{(2k-1)^2} \left(\sum^{k}_{i=1} \frac{h_i}{i} \right)~.
\end{eqnarray}
Here only the first and last nested double sums on the right side of Eq.~(273) are unknown. For the first sum it follows:
\begin{eqnarray}
\sum^{\infty}_{k=1} \frac{h_k}{k(2k-1)} \left( \sum^{k}_{i=1} \frac{h_i^{(2)}}{i^2} \right) &=& \zeta(2) \left( \sum^{\infty}_{k=1} \frac{h_k^{(2)}}{k^2} \right) + \sum^{\infty}_{k=1} \frac{h_k h_k^{(2)}}{k^3(2k-1)}
- \sum^{\infty}_{k=1} \frac{h_k^{(2)}}{k^2}  \left( \sum^{k}_{i=1} \frac{h_i}{i(2i-1)} \right) \nonumber \\ &=&
\zeta(2) \left( \sum^{\infty}_{k=1} \frac{h_k^{(2)}}{k^2} \right) + \sum^{\infty}_{k=1} \frac{h_k h_k^{(2)}}{k^3(2k-1)} - \sum^{\infty}_{k=1} \frac{h_k^{(2)} h_k^{(2)}}{k^2} -
\sum^{\infty}_{k=1} \frac{h_k^{(2)} h_k^2}{k^2} \nonumber \\ &+& \sum^{\infty}_{k=1} \frac{h_k^{(2)}}{k^2} \left( \sum^{k}_{i=1} \frac{h_i}{i} \right)~. 
\end{eqnarray}
Again only the last nested double sum on the right side of Eq.~(274) is unknown. We make the ansatz:
\begin{eqnarray}
\sum^{\infty}_{k=1} \frac{h_k^{(2)}}{k} \left( \sum^{\infty}_{i=1} \frac{h_i}{i(k+i)} \right) &=& 2ln(2)\sum^{\infty}_{k=1} \frac{h_k h_k^{(2)}}{k^2} + \sum^{\infty}_{k=1} \frac{H_k h_k h_k^{(2)}}{k^2} -
\sum^{\infty}_{i=1} \frac{h_k^{(2)}}{k^2} \left( \sum^{k}_{i=1} \frac{h_i}{i} \right) \nonumber \\ &=&
\sum^{\infty}_{i=1} \frac{h_i}{i} \left( \sum^{\infty}_{k=1} \frac{h_k^{(2)}}{k(i+k)} \right) \nonumber \\ &=&
\frac{3}{2}\zeta(2) \sum^{\infty}_{k=1} \frac{h_k^2}{k^2} - 2ln(2)\sum^{\infty}_{k=1} \frac{h_k h_k^{(2)}}{k^2} - \sum^{\infty}_{k=1} \frac{h_k}{k^2} \left( \sum^{k}_{i=1} \frac{H_{i-1}}{(2i-1)^2} \right)~.
\end{eqnarray}
And again only the last nested double sum on the right side of Eq.~(275) is unknown. We start the corresponding calculation as follows:
\begin{eqnarray}
\sum^{\infty}_{i=1} \frac{H_i}{2i-1} \left( \sum^{\infty}_{k=1} \frac{h_k^{(2)}}{k(2i+2k-1)} \right) &=& \sum^{\infty}_{k=1} \frac{h_k^{(2)}}{k} \left( \sum^{\infty}_{i=1} \frac{H_i}{(2i-1)(2i+2k-1)} \right)~.
\end{eqnarray}
From this we get:
\begin{eqnarray}
2ln(2) \sum^{\infty}_{k=1} \frac{h_k^{(2)}}{k(2k-1)} - ln(2)\sum^{\infty}_{k=1} \frac{h_k h_k^{(2)}}{k^2} - \sum^{\infty}_{k=1} \frac{h_k h_k^{(2)}}{k^2(2k-1)} + \frac{1}{2}\sum^{\infty}_{k=1}
\frac{h_k^{(2)} h_k^{(2)}}{k^2} + \frac{1}{2}\sum^{\infty}_{k=1} \frac{h_k^2 h_k^{(2)}}{k^2} = \nonumber \\
\left( \frac{3}{2}ln(2)\zeta(2) - \frac{7}{8}\zeta(3) \right) \sum^{\infty}_{i=1} \frac{H_i}{(2i-1)^2} + \frac{3}{4}\zeta(2) \sum^{\infty}_{i=1} \frac{H_i H_{i-1}}{(2i-1)^2} -
\frac{1}{2} \sum^{\infty}_{i=1} \frac{H_i}{(2i-1)^2}  \left( \sum^{i-1}_{k=1} \frac{h_k}{k^2} \right)  = \nonumber \\
\left( \frac{3}{2}ln(2)\zeta(2) - \frac{7}{8}\zeta(3) \right) \sum^{\infty}_{i=1} \frac{H_i}{(2i-1)^2} + \frac{3}{4}\zeta(2) \sum^{\infty}_{i=1} \frac{H_i H_{i-1}}{(2i-1)^2} -
\frac{1}{2} \sum^{\infty}_{i=1} \frac{H_i}{(2i-1)^2}  \left( \sum^{\infty}_{k=1} \frac{h_k}{k^2} \right) +  \nonumber \\
\frac{1}{2}\sum^{\infty}_{k=1} \frac{h_k}{k^2} \left( \sum^{k}_{i=1} \frac{H_{i-1}}{(2i-1)^2} \right)~.
\end{eqnarray}
As all sums are known explicitly the last nested double sum on the right side of Eq.~(277) can be computed, and with this the other nested double sums. This way it results for Eq.~(278):
\begin{eqnarray}
\sum^{\infty}_{k=1} \frac{h_k}{k(2k-1)} \left( \sum^{k}_{i=1} \frac{h_i^{(2)}}{i^2} \right) &=& - \frac{135}{128}\zeta(6) + \frac{31}{4}ln(2)\zeta(5) - \frac{31}{4}\zeta(5) + \frac{105}{16}\zeta(4) +
\frac{203}{96}\zeta(3)^2 - \frac{7}{4}ln(2)\zeta(2)\zeta(3) \nonumber \\ &+& \frac{7}{4}\zeta(2)\zeta(3) + \frac{7}{2}ln(2)\zeta(3) - 3\left(ln(2)\right)^2\zeta(2) - \frac{9}{4}\zeta(2)\sum^{\infty}_{k=1}\frac{h_k}{k^3}
- \frac{1}{2}\sum^{\infty}_{k=1}\frac{h_k}{k^3} - \frac{1}{4}\sum^{\infty}_{k=1}\frac{h_k}{k^5} \nonumber \\ &-& \frac{1}{3} \sum^{\infty}_{k=1} \frac{h_k^3}{k^3}~. 
\end{eqnarray}
This way it remains to calculate the last sum on the right side of Eq.~(273). With the following ansatz it results first:
\begin{eqnarray}
\sum^{\infty}_{k=1} \frac{H_k^{(2)}}{2k-1} \left( \sum^{\infty}_{i=1} \frac{h_i}{i(2i+2k-1)} \right) &=& \sum^{\infty}_{i=1} \frac{h_i}{i} \left( \sum^{\infty}_{k=1} \frac{H_k^{(2)}}{(2k-1)(2i+2k-1)} \right)
\end{eqnarray}
or
\begin{eqnarray}
\frac{1}{2}\zeta(2)\sum^{\infty}_{k=1} \frac{H_k^{(2)}}{(2k-1)^2} &+& \sum^{\infty}_{k=1} \frac{H_k^{(2)}}{(2k-1)^2} \left( \sum^{k-1}_{i=1} \frac{h_i}{i} \right) = \left( 4ln(2) - \zeta(2)\right)
\sum^{\infty}_{i=1} \frac{h_i}{i(2i-1)} - 4ln(2)\sum^{\infty}_{i=1} \frac{h_i}{i(2i-1)^2} \nonumber \\ &+& \frac{1}{2}\zeta(2)\sum^{\infty}_{i=1} \frac{h_i^2}{i^2} + 2ln(2)\sum^{\infty}_{i=1} \frac{h_i h_i^{(2)}}{i^2}
- 2\sum^{\infty}_{i=1} \frac{h_i}{i^2} \left( \sum^{i-1}_{k=1} \frac{h_k}{(2k-1)^2} \right)~.
\end{eqnarray}
Only the last nested double sum on the right side of Eq.~(280) is unknown. We make the ansatz:
\begin{eqnarray}
\sum^{\infty}_{k=1} \frac{h_k}{2k-1} \left( \sum^{\infty}_{i=1} \frac{h_i^{(2)}}{i(2i+2k-1)} \right) &=& \sum^{\infty}_{i=1} \frac{h_i^{(2)}}{i} \left( \sum^{\infty}_{k=1} \frac{h_k}{(2k-1)(2i+2k-1)} \right) 
\end{eqnarray}
or
\begin{eqnarray}
\left( \frac{3}{2}ln(2)\zeta(2) - \frac{7}{8}\zeta(3)\right) \sum^{\infty}_{k=1} \frac{h_k}{(2k-1)^2} + \frac{3}{4}\zeta(2) \sum^{\infty}_{k=1} \frac{H_{k-1}h_k}{(2k-1)^2} - \frac{1}{2} \sum^{\infty}_{k=1} \frac{h_k}{k^2}
\left( \sum^{\infty}_{i=1} \frac{h_i}{(2i-1)^2} \right) \nonumber \\ + \frac{1}{2} \sum^{\infty}_{i=1} \frac{h_k}{k^2} \left( \sum^{k}_{i=1} \frac{h_i}{(2i-1)^2} \right) =
 \frac{3}{8}\zeta(2)\sum^{\infty}_{i=1} \frac{h_i^{(2)}}{i^2} + \frac{1}{4}\sum^{\infty}_{i=1} \frac{h_i^{(2)}}{i^2} \left( \sum^{i-1}_{k=1} \frac{h_k}{k} \right)~.
\end{eqnarray}
Again only the last nested double sum on the right side of Eq.~(282) is unknown. With the following ansatz it results first: 
\begin{eqnarray}
\sum^{\infty}_{k=1} \frac{h_k}{k} \left( \sum^{\infty}_{i=1} \frac{h_i^{(2)}}{i(k+i)} \right) &=& \sum^{\infty}_{i=1} \frac{h_i^{(2)}}{i} \left( \sum^{\infty}_{k=1} \frac{h_k}{k(i+k}) \right) 
\end{eqnarray}
or
\begin{eqnarray}
2ln(2)\sum^{\infty}_{k=1} \frac{h_k h_k^{(2)}}{k^2} + \sum^{\infty}_{k=1} \frac{H_k h_k h_k^{(2)}}{k^2} - \sum^{\infty}_{k=1} \frac{h_k^{(2)}}{k^2} \left( \sum^{k-1}_{i=1} \frac{h_i}{i} \right) &=&
\frac{3}{2}\zeta(2)\sum^{\infty}_{k=1} \frac{h_k^2}{k^2} - 2ln(2)\sum^{\infty}_{k=1} \frac{h_k h_k^{(2)}}{k^2} \nonumber \\ &-& \sum^{\infty}_{k=1} \frac{h_k}{k^2} \left( \sum^{k}_{i=1} \frac{H_{i-1}}{(2i-1)^2} \right)~.
\end{eqnarray}
The last sum on the right side of Eq.~(284) is explicitly known from Eq.~(276). This way the last sum on the left side of Eq.~(284) can also be explicitly computed. It follows:
\begin{eqnarray}
\sum^{\infty}_{k=1} \frac{h_k^{(2)}}{k^2} \left( \sum^{k}_{i=1} \frac{h_i}{i} \right) = -\frac{405}{96}\zeta(6) + \frac{31}{4}ln(2)\zeta(5) + \frac{273}{96}\zeta(3)^2 - \frac{7}{4}ln(2)\zeta(2)\zeta(3)
+ \frac{1}{4}\zeta(2)\sum^{\infty}_{k=1}\frac{h_k}{k^3} - \frac{1}{4}\sum^{\infty}_{k=1}\frac{h_k}{k^5}~.
\end{eqnarray}
With this it follows: 
\begin{eqnarray}
\sum^{\infty}_{k=1} \frac{H^{(2)}_{k}}{(2k-1)^2} \left(\sum^{k-1}_{i=1} \frac{h_i}{i} \right) &=& \frac{15}{128}\zeta(6) + \frac{15}{8}\zeta(4) + \frac{7}{32}\zeta(3)^2 + \frac{7}{2}ln(2)\zeta(3) - 3\zeta(3)
-\sum^{\infty}_{k=1}\frac{h_k}{k^3} + \frac{1}{4}\sum^{\infty}_{k=1}\frac{h_k}{k^5}~.
\end{eqnarray}
Now are all sums necessary for an explicit calculation of the fourth member are known and it follows:
\begin{eqnarray}
\sum^{\infty}_{k=1} \frac{H_k^{(2)} h_k h_k^{(2)}}{k(2k-1)} &=& \frac{4857}{128}\zeta(6) - \frac{527}{16}ln(2)\zeta(5) - \frac{31}{4}\zeta(5) + \frac{121}{16}\left(ln(2)\right)^2\zeta(4) 
+ \frac{105}{16}\zeta(4) + \frac{525}{128}\zeta(3)^2 \nonumber \\ &+& \frac{21}{8}ln(2)\zeta(2)\zeta(3) + \frac{7}{4}\zeta(2)\zeta(3) + \frac{7}{2}ln(2)\zeta(3) - \frac{1}{3}\left(ln(2)\right)^4\zeta(2) 
- 3\left(ln(2)\right)^2\zeta(2) + \frac{1}{15}\left(ln(2)\right)^6 \nonumber \\ &-& \frac{3}{8}\zeta(2)\sum^{\infty}_{k=1}\frac{h_k}{k^3} - \frac{1}{2}\sum^{\infty}_{k=1}\frac{h_k}{k^3}
+ \frac{61}{32}\sum^{\infty}_{k=1}\frac{h_k}{k^5} - 48Li_6\left(\frac{1}{2}\right) -16ln(2)Li_5\left(\frac{1}{2}\right)~.
\end{eqnarray}

In order to calculate the fifth member we start with the ansatz:
\begin{eqnarray}
\sum^{\infty}_{i=1} \frac{H_i^{(3)}}{2i-1} \left( \sum^{\infty}_{k=1} \frac{h_k}{k(i+k)} \right) = \sum^{\infty}_{k=1} \frac{h_k}{k} \left( \sum^{\infty}_{i=1} \frac{H_i^{(3)}}{(2i-1)(k+i)} \right) 
\end{eqnarray}
or by use of the corresponding two-valued help functions:
\begin{eqnarray}
2ln(2)\sum^{\infty}_{i=1} \frac{H_i^{(3)}}{i(2i-1)} + \sum^{\infty}_{i=1} \frac{H_i H_i^{(3)} h_i}{i(2i-1)} - \sum^{\infty}_{i=1} \frac{H_i^{(3)}}{i(2i-1)} \left( \sum^{i}_{k=1} \frac{h_k}{k} \right) =
\sum^{\infty}_{k=1} \frac{h_k}{k(2k+1)} \left( \sum^{\infty}_{i=1} \frac{H_i^{(3)}}{i(2i-1)} \right) + \nonumber \\ \zeta(4)\sum^{\infty}_{k=1} \frac{h_k}{k(2k+1)} + \zeta(3)\sum^{\infty}_{k=1}
\frac{{H_{k-1}} h_k}{k(2k+1)} - \zeta(2)\sum^{\infty}_{k=1}\frac{{H_{k-1}^{(2)}} h_k}{k(2k+1)} + \sum^{\infty}_{k=1} \frac{h_k}{k(2k+1)} \left( \sum^{k-1}_{i=1} \frac{H_i}{i^3} \right)~. 
\end{eqnarray}
Here are the last nested  double sum on the left side and the last nested double sum on the right side are unknown. We start with the calculation of the last nested double sum appearing on the left side
of Eq.~(289). First we get:
\begin{eqnarray}
\sum^{\infty}_{i=1} \frac{H_i^{(3)}}{i} \left( \sum^{\infty}_{k=1} \frac{h_k}{k(2i+2k-1)} \right) = \sum^{\infty}_{k=1} \frac{h_k}{k} \left(\sum^{\infty}_{i=1} \frac{H_i^{(3)}}{i(2i+2k-1)} \right)
\end{eqnarray}
or
\begin{eqnarray} 
\frac{1}{2}\zeta(2)\sum^{\infty}_{i=1} \frac{H_i^{(3)}}{i(2i-1)} &+& \sum^{\infty}_{i=1} \frac{H_i^{(3)}}{i(2i-1)} \left( \sum^{i-1}_{k=1} \frac{h_k}{k} \right) = - \sum^{\infty}_{k=1}
\frac{h_k}{k(2k-1)}\left(\sum^{\infty}_{i=1} \frac{H_i^{(3)}}{i(2i+1)}\right) + 2\zeta(3)\sum^{\infty}_{k=1}\frac{h_k h_{k-1}}{k(2k-1)} \nonumber \\ &-& 4\zeta(2)\sum^{\infty}_{k=1}
\frac{h_k h_{k-1}^{(2)}}{k(2k-1)} -16ln(2)\sum^{\infty}_{k =1}\frac{h_k h_{k-1}^{(3)}}{k(2k-1)} + 16\sum^{\infty}_{k=1}\frac{h_k}{k(2k-1)} \left( \sum^{k-1}_{i=1} \frac{h_i}{(2i-1)^3} \right)~.  
\end{eqnarray}
Only the last nested double sum on the right side of Eq.~(291) is unknown. It follows first:
\begin{eqnarray}
\sum^{\infty}_{k=1}\frac{h_k}{k(2k-1)} \left( \sum^{k-1}_{i=1} \frac{h_i}{(2i-1)^3} \right) = \sum^{\infty}_{k=1}\frac{h_k}{k(2k-1)} \left( \sum^{\infty}_{i=1} \frac{h_i}{(2i-1)^3} \right) -
\sum^{\infty}_{i=1} \frac{h_i}{(2i-1)^3} \left( \sum^{i}_{k=1}\frac{h_k}{k(2k-1)} \right) = \nonumber \\ \sum^{\infty}_{k=1}\frac{h_k}{k(2k-1)} \left( \sum^{\infty}_{i=1} \frac{h_i}{(2i-1)^3} \right)
- \sum^{\infty}_{i=1}\frac{h_i h_i^{(2)}}{(2i-1)^3} - \sum^{\infty}_{i=1}\frac{h_i^3}{(2i-1)^3} + \sum^{\infty}_{i=1} \frac{h_i}{(2i-1)^3}  \left( \sum^{i}_{k=1}\frac{h_k}{k} \right)~. 
\end{eqnarray}
As the next to the last sum on the right side of Eq.~(292) is defined by Eq.~(221) only the last sum on the right side is unknown. It follows further:
\begin{eqnarray}
\sum^{\infty}_{i=1} \frac{h_i}{(2i-1)^2} \left( \sum^{\infty}_{k=1}\frac{h_k}{k(2k+2i-1)} \right) = \sum^{\infty}_{k=1}\frac{h_k}{k} \left(\sum^{\infty}_{i=1} \frac{h_i}{(2i-1)^2(2k+2i-1)} \right)
\end{eqnarray}
or
\begin{eqnarray}
\frac{1}{2}\zeta(2)\sum^{\infty}_{i=1} \frac{h_i}{(2i-1)^3} - \sum^{\infty}_{i=1} \frac{h_i^2}{i(2i-1)^3} &+& \sum^{\infty}_{i=1} \frac{h_i}{(2i-1)^3} \left( \sum^{i}_{k=1}\frac{h_k}{k} \right) =
\frac{1}{2}\sum^{\infty}_{k=1}\frac{h_k}{k^2} \left( \sum^{\infty}_{i=1} \frac{h_i}{(2i-1)^2} \right) \nonumber \\ &-& \frac{3}{16}\zeta(2)\sum^{\infty}_{k=1}\frac{h_k}{k^3} +
\frac{1}{8}\sum^{\infty}_{k=1}\frac{h_k^2}{k^4} - \frac{1}{8}\sum^{\infty}_{k=1}\frac{h_k}{k^3} \left( \sum^{k}_{i=1}\frac{h_i}{i} \right)~.
\end{eqnarray}
Here, the last nested double sum is known from Eq.~(149). All other sum are also known explicitly. This way the last nested double sum on the left side of Eq.~(289) follows explicitly. As a final step
we calculate the last nested double sum on the right side of Eq.~(289). It follows first by shifting the indices:
\begin{eqnarray}
\sum^{\infty}_{k=1} \frac{h_k}{k(2k+1)} \left( \sum^{k}_{i=1} \frac{H_i}{i^3} \right) = 2\sum^{\infty}_{k=1} \frac{H_k h_{k-1}}{k^3(2k-1)} + 2\sum^{\infty}_{k=1} \frac{1}{(2k-1)^2}
\left( \sum^{k}_{i=1} \frac{H_i}{i^3} \right) - \sum^{\infty}_{k=1} \frac{h_k}{k(2k-1)} \left( \sum^{k}_{i=1} \frac{H_i}{i^3} \right)~.
\end{eqnarray}
Therefore it remains to calculate the last sum on the right side of Eq.~(295) explicitly. It follows:
\begin{eqnarray}
\sum^{\infty}_{k=1} \frac{h_k}{k(2k-1)} \left( \sum^{k}_{i=1} \frac{H_i}{i^3} \right) &=& \zeta(2)\sum^{\infty}_{k=1} \frac{H_k}{k^3} + \sum^{\infty}_{k=1} \frac{H_k h_k}{k^4(2k-1)} -
\sum^{\infty}_{k=1} \frac{H_k}{k^3} \left( \sum^{k}_{i=1} \frac{h_i}{i(2i-1)} \right) \nonumber \\ &=& \zeta(2)\sum^{\infty}_{k=1} \frac{H_k}{k^3} + \sum^{\infty}_{k=1} \frac{H_k h_k}{k^4(2k-1)} -
\sum^{\infty}_{k=1}\frac{H_k k^{(2)}}{k^3} - \sum^{\infty}_{k=1}\frac{H_k k^2}{k^3} + \sum^{\infty}_{k=1} \frac{H_k}{k^3} \left( \sum^{k}_{i=1} \frac{h_i}{i} \right)~.
\end{eqnarray}
As the last sum in  Eq.~(296) is defined by Eq.~(219) the nested double sum can be calculated explicitly. This way the fifth member results to: 

\begin{eqnarray}
\sum^{\infty}_{k=1} \frac{H_k H_k^{(3)} h_k}{k(2k-1)} &=& \frac{45}{8}\zeta(6) - \frac{279}{4}ln(2)\zeta(5) + \frac{279}{4}\zeta(5) + \frac{53}{2}\left(ln(2)\right)^2\zeta(4) - 
\frac{53}{2}ln(2)\zeta(4) - \frac{45}{4}\zeta(4) - \frac{149}{32}\zeta(3)^2 \nonumber \\ &+& \frac{45}{2}ln(2)\zeta(2)\zeta(3) - \frac{45}{2}\zeta(2)\zeta(3) -
14\left(ln(2)\right)^3\zeta(3) + 14\left(ln(2)\right)^2\zeta(3) + 10\zeta(3) + \frac{8}{3}\left(ln(2)\right)^4\zeta(2) \nonumber \\ &-& \frac{8}{3}\left(ln(2)\right)^3\zeta(2) -
16ln(2)\zeta(2) + 16\zeta(2) - \frac{4}{15}\left(ln(2)\right)^6 + \frac{4}{15}\left(ln(2)\right)^5 - \frac{3}{2}\zeta(2)\sum^{\infty}_{k=1}\frac{h_k}{k^3} \nonumber \\ &+&
2\left(ln(2)\right)^2\sum^{\infty}_{k=1}\frac{h_k}{k^3} - 2\sum^{\infty}_{k=1}\frac{h_k}{k^3} + \frac{7}{8}\sum^{\infty}_{k=1}\frac{h_k}{k^5}
+ 32ln(2)Li_5\left(\frac{1}{2}\right) - 32Li_5\left(\frac{1}{2}\right)~.
\end{eqnarray}

The calculational procedure concerning the sixth member starts as follows:
\begin{eqnarray}
\sum^{\infty}_{k=1} \frac{h_k^{(3)}}{2k-1} \left( \sum^{\infty}_{i=1} \frac{h_i}{i(k+i)} \right) = \sum^{\infty}_{i=1} \frac{h_i}{i} \left( \sum^{\infty}_{k=1} \frac{h_k^{(3)}}{(2k-1)(k+i)} \right) 
\end{eqnarray}
or
\begin{eqnarray}
2ln(2)\sum^{\infty}_{k=1} \frac{h_k h_k^{(3)}}{k(2k-1)} + \sum^{\infty}_{k=1} \frac{H_k h_k h_k^{(3)}}{k(2k-1)} - \sum^{\infty}_{k=1} \frac{h_k^{(3)}}{k(2k-1)} \left( \sum^{k}_{i=1}\frac{h_i}{i} \right) =
\nonumber \\
\frac{7}{4}\zeta(3) \sum^{\infty}_{k=1} \frac{h_k^2}{k(2k+1)} - \frac{3}{2}\zeta(2) \sum^{\infty}_{k=1} \frac{h_k h_k^{(2)}}{k(2k+1)} + 2ln(2)\sum^{\infty}_{k=1} \frac{h_k h_k^{(3)}}{k(2k+1)}
- \sum^{\infty}_{k=1} \frac{h_k}{k(2k+1)} \left( \sum^{k}_{i=1}\frac{H_{i-1}}{(2i-1)^3} \right)~.
\end{eqnarray}
Here are the nested double sum on the left side and the nested double sum on the right side are unknown. First we calculate the nested double sum on the left side. It follows:
\begin{eqnarray}
\sum^{\infty}_{k=1}\frac{h_k}{k} \left( \sum^{\infty}_{i=1}\frac{h_i^{(3)}}{i(2i+2k-1)} \right) = \sum^{\infty}_{k=1}\frac{h_i^{(3)}}{i} \left(\sum^{\infty}_{k=1}\frac{h_k}{k(2i+2k-1)} \right) 
\end{eqnarray}
or
\begin{eqnarray}
\sum^{\infty}_{k=1} \frac{h_k}{k(2k-1)} \left( \sum^{\infty}_{i=1}\frac{h_i^{(3)}}{i(2i+1)} \right) + \frac{7}{8}\zeta(3) \sum^{\infty}_{k=1} \frac{H_{k-1} h_k}{k(2k-1)} - \frac{3}{8}\zeta(2)
\sum^{\infty}_{k=1} \frac{H_{k-1}^{(2)} h_k}{k(2k-1)} + \frac{1}{4}\sum^{\infty}_{k=1} \frac{h_k}{k(2k-1)} \left( \sum^{k-1}_{i=1}\frac{h_i}{i^3} \right) &=& \nonumber \\
\frac{1}{2}\zeta(2)\sum^{\infty}_{k=1} \frac{h_k^{(3)}}{k(2k-1)} + \sum^{\infty}_{k=1} \frac{h_k^{(3)}}{k(2k-1)} \left( \sum^{k-1}_{i=1}\frac{h_i}{i} \right)~.
\end{eqnarray}
In Eq.(301) only the nested double sum on the left side must be calculated in order to compute explicitly the double sum on the right side. It follows further:
\begin{eqnarray}
\sum^{\infty}_{k=1} \frac{h_k}{k(2k-1)} \left( \sum^{k-1}_{i=1}\frac{h_i}{i^3} \right) = \sum^{\infty}_{k=1} \frac{h_k}{k(2k-1)} \left( \sum^{\infty}_{i=1}\frac{h_i}{i^3} \right) 
- \sum^{\infty}_{k=1} \frac{h_k}{k^3} \left( \sum^{k}_{i=1}\frac{h_i}{i(2i-1)} \right)  \nonumber \\
\end{eqnarray}
or
\begin{eqnarray}
\sum^{\infty}_{k=1} \frac{h_k}{k(2k-1)} \left( \sum^{k-1}_{i=1}\frac{h_i}{i^3} \right) = \sum^{\infty}_{k=1} \frac{h_k}{k(2k-1)} \left( \sum^{\infty}_{i=1}\frac{h_i}{i^3} \right) 
- \sum^{\infty}_{k=1} \frac{h_k h_k^{(2)}}{k^3} - \sum^{\infty}_{k=1} \frac{h_k^3}{k^3} + \sum^{\infty}_{k=1} \frac{h_k}{k^3} \left( \sum^{k}_{i=1}\frac{h_i}{i} \right)~. 
\end{eqnarray}
The nested double sum on the right side of Eq.~(303) is explicitly known from Eq.~(149). This way the nested double sum on the left side of Eq.~(299) is explicitly known. Therefore it remains to calculate the
nested double sum on the right side of Eq.~(299). It follows: 
\begin{eqnarray}
\sum^{\infty}_{k=1} \frac{h_k}{k(2k+1)} \left( \sum^{k}_{i=1}\frac{H_{i-1}}{(2i-1)^3} \right) &=& \sum^{\infty}_{k=1} \frac{h_k}{k(2k+1)} \left( \sum^{\infty}_{i=1}\frac{H_{i-1}}{(2i-1)^3} \right)
+ \sum^{\infty}_{k=1} \frac{{H_{k-1} h_k}}{k(2k+1)(2k-1)^3} \nonumber \\ &-& \sum^{\infty}_{k=1}\frac{H_{k-1}}{(2k-1)^3} \left( \sum^{k}_{i=1}\frac{h_i}{i} - 2\sum^{k}_{i=1}\frac{h_i}{2i+1} \right) \nonumber \\ &=&
\sum^{\infty}_{k=1} \frac{h_k}{k(2k+1)} \left( \sum^{\infty}_{i=1}\frac{H_{i-1}}{(2i-1)^3} \right) + \sum^{\infty}_{k=1} \frac{{H_{k-1} h_k}}{k(2k+1)(2k-1)^3} + 2\sum^{\infty}_{k=1}
\frac{{H_{k-1} h_k}}{(2k+1)(2k-1)^3} \nonumber \\ &+& \sum^{\infty}_{k=1} \frac{H_{k-1} h_k^{(2)}}{(2k-1)^3} + \sum^{\infty}_{k=1} \frac{H_{k-1} h_k^2}{(2k-1)^3} - \sum^{\infty}_{k=1}\frac{H_{k-1}}{(2k-1)^3}
\left( \sum^{k}_{i=1}\frac{h_i}{i} \right)~.
\end{eqnarray}
Again, it remains to calculate the nested double sum on the right side of Eq.~(304). Starting with the ansatz: 
\begin{eqnarray}
\sum^{\infty}_{i=1} \frac{h_i}{i} \left( \sum^{\infty}_{k=1} \frac{H_k}{(2k-1)^2(2i+2k-1)} \right) &=& \sum^{\infty}_{k=1} \frac{H_k}{(2k-1)^2} \left( \sum^{\infty}_{i=1} \frac{h_i}{i(2i+2k-1)} \right) \nonumber \\ &=&
\frac{1}{2}\zeta(2) \sum^{\infty}_{k=1} \frac{H_k}{(2k-1)^3} + \sum^{\infty}_{k=1} \frac{H_k}{(2k-1)^3} \left( \sum^{k-1}_{i=1} \frac{h_i}{i}\right)
\end{eqnarray}
or 
\begin{eqnarray}
\sum^{\infty}_{i=1} \frac{h_i}{i} \left( \sum^{\infty}_{k=1} \frac{H_k}{(2k-1)^2(2i+2k-1)} \right) &=& \frac{1}{2} \sum^{\infty}_{i=1} \frac{h_i}{i^2} \left( \sum^{\infty}_{k=1} \frac{H_k}{(2k-1)^2} \right) -
ln(2) \sum^{\infty}_{k=1} \frac{h_k}{k^2(2k-1)} + \frac{1}{2}ln(2) \sum^{\infty}_{k=1} \frac{h_k^2}{k^3} \nonumber \\ &+& \frac{1}{2} \sum^{\infty}_{k=1} \frac{h_k^2}{k^3(2k-1)} - \frac{1}{4}
\frac{h_k h_k^{(2)}}{k^3} - \frac{1}{4}  \frac{h_k^3}{k^3}
\end{eqnarray}
a corresponding identity for this nested double sum follows. From this we get:
\begin{eqnarray}
\sum^{\infty}_{k=1} \frac{H_k}{(2k-1)^3} \left( \sum^{k}_{i=1} \frac{h_i}{i} \right) &=& \sum^{\infty}_{k=1} \frac{H_k^2 h_k}{(2k-1)^3}  -\frac{315}{128}\zeta(6) + \frac{31}{8}ln(2)\zeta(5) - \frac{217}{64}\zeta(5)
+ \frac{15}{4}\left(ln(2)\right)^2\zeta(4)
- \frac{15}{16}ln(2)\zeta(4) + \frac{15}{8}\zeta(4) \nonumber \\ &+& \frac{49}{32}\zeta(3)^2 - \frac{35}{16}ln(2)\zeta(2)\zeta(3) + \frac{21}{16}\zeta(2)\zeta(3) - \frac{7}{2}\left(ln(2)\right)^3\zeta(3) +
\frac{7}{2}\left(ln(2)\right)^2\zeta(3) \nonumber \\ &+& \frac{7}{4}ln(2)\zeta(3) - \frac{9}{4}\zeta(3) - 3\left(ln(2)\right)^2\zeta(2) + 2ln(2)\zeta(2) - \frac{1}{8}\zeta(2)\sum^{\infty}_{k=1}\frac{h_k}{k^3}
+ \frac{1}{2}\left(ln(2)\right)^2\sum^{\infty}_{k=1}\frac{h_k}{k^3} \nonumber \\ &-& \frac{1}{2}ln(2)\sum^{\infty}_{k=1}\frac{h_k}{k^3} - \frac{1}{2}\sum^{\infty}_{k=1}\frac{h_k}{k^3} -
\frac{1}{8}\sum^{\infty}_{k=1}\frac{h_k}{k^5}~.
\end{eqnarray}

This way the sixth member can be calculated explicitly. The result is:
\begin{eqnarray}
\sum^{\infty}_{k=1} \frac{H_k h_k h_k^{(3)}}{k(2k-1)} &=& \frac{135}{256}\zeta(6) - \frac{31}{16}ln(2)\zeta(5) + \frac{341}{64}\zeta(5) + \frac{1}{32}\left(ln(2)\right)^2\zeta(4) - 
\frac{23}{16}ln(2)\zeta(4) - \frac{7}{256}\zeta(3)^2 - \frac{23}{16}ln(2)\zeta(2)\zeta(3)  \nonumber \\ &+& \frac{9}{16}\zeta(2)\zeta(3)
+ \frac{1}{3}\left(ln(2)\right)^4\zeta(2) - \frac{1}{3}\left(ln(2)\right)^3\zeta(2) - \frac{1}{30}\left(ln(2)\right)^6 + \frac{1}{30}\left(ln(2)\right)^5 +
\frac{5}{16}\zeta(2)\sum^{\infty}_{k=1}\frac{h_k}{k^3} \nonumber \\ &-& \frac{7}{64}\sum^{\infty}_{k=1}\frac{h_k}{k^5}
+ 4ln(2)Li_5\left(\frac{1}{2}\right) - 4Li_5\left(\frac{1}{2}\right)~.
\end{eqnarray}

In order to calculate the seventh member we start with the identity
\begin{eqnarray}
\sum^{\infty}_{k=1}\frac{H_k}{k} \left( \sum^{\infty}_{i=1}\frac{H_i^{(3)}}{i(2i+2k-1)} \right) = \sum^{\infty}_{k=1}\frac{H_i^{(3)}}{i} \left(\sum^{\infty}_{k=1}\frac{H_k}{k(2i+2k-1)} \right) 
\end{eqnarray}
or
\begin{eqnarray}
-\sum^{\infty}_{i=1} \frac{H_i}{i(2i-1)} \left( \sum^{\infty}_{k=1}\frac{H_k^{(3)}}{k(2k+1)} \right) + 2\zeta(3)\sum^{\infty}_{i=1} \frac{H_i h_{i-1}}{i(2i-1)} - 4\zeta(2)
\sum^{\infty}_{i=1} \frac{H_i h_{i-1}^{(2)}}{i(2i-1)} - 16ln(2)\sum^{\infty}_{i=1} \frac{H_i h_{i-1}^{(3)}}{i(2i-1)} +  \nonumber \\ 16\sum^{\infty}_{i=1} \frac{H_i}{i(2i-1)}
\left( \sum^{i-1}_{k=1}\frac{h_k}{(2k-1)^3} \right) = 2\left(ln(2)\right)^2\sum^{\infty}_{k=1}\frac{H_k^{(3)}}{k(2k-1)} + 4ln(2)\sum^{\infty}_{k=1}\frac{H_k^{(3)}}{k(2k-1)^2} -
4\sum^{\infty}_{k=1}\frac{H_k^{(3)} h_k}{k(2k-1)^2} -  \nonumber \\ 4ln(2)\sum^{\infty}_{k=1}\frac{H_k^{(3)} h_k}{k(2k-1)} + 2\sum^{\infty}_{k=1}\frac{H_k^{(3)} h_k^{(2)}}{k(2k-1)} +
2\sum^{\infty}_{k=1}\frac{H_k^{(3)} h_k^2}{k(2k-1)}~.
\end{eqnarray}
The last Euler sum on the right side of Eq.~(310) can be identified as the seventh member of the seventh family. Therefore, only the nested double sum on the left side of Eq.~(310) is unknown.
We start the corresponding calculational procedure with:
\begin{eqnarray}
\sum^{\infty}_{i=1} \frac{H_i}{i(2i-1)} \left( \sum^{i-1}_{k=1}\frac{h_k}{(2k-1)^3} \right)  &=& \sum^{\infty}_{i=1} \frac{H_i}{i(2i-1)} \left( \sum^{\infty}_{k=1}\frac{h_k}{(2k-1)^3} \right) -
2\sum^{\infty}_{k=1}\frac{h_k}{(2k-1)^3} \left( \sum^{k}_{i=1} \frac{H_i}{(2i-1)} \right) + \nonumber \\ &+& \sum^{\infty}_{k=1}\frac{h_k}{(2k-1)^3} \left( \sum^{k}_{i=1} \frac{H_i}{i} \right) 
\end{eqnarray}
or
\begin{eqnarray}
\sum^{\infty}_{i=1} \frac{H_i}{i(2i-1)} \left( \sum^{i-1}_{k=1}\frac{h_k}{(2k-1)^3} \right)  &=& \sum^{\infty}_{i=1} \frac{H_i}{i(2i-1)} \left( \sum^{\infty}_{k=1}\frac{h_k}{(2k-1)^3} \right) -
2\sum^{\infty}_{k=1}\frac{h_k}{(2k-1)^3} \left( \sum^{k}_{i=1} \frac{H_{i-1}}{(2i-1)} \right) \nonumber \\ &-& 4\sum^{\infty}_{k=1}\frac{h_k^2}{(2k-1)^3} + 2\sum^{\infty}_{k=1}\frac{H_k h_k}{(2k-1)^3} 
+ \frac{1}{2}\sum^{\infty}_{k=1}\frac{H_k^{(2)} h_k}{(2k-1)^3} + \frac{1}{2}\sum^{\infty}_{k=1}\frac{H_k^2 h_k}{(2k-1)^3}
\end{eqnarray}
or
\begin{eqnarray}
\sum^{\infty}_{i=1} \frac{H_i}{i(2i-1)} \left( \sum^{i-1}_{k=1}\frac{h_k}{(2k-1)^3} \right)  &=& \sum^{\infty}_{i=1} \frac{H_i}{i(2i-1)} \left( \sum^{\infty}_{k=1}\frac{h_k}{(2k-1)^3} \right) -
2\sum^{\infty}_{k=1}\frac{H_k h_k^2}{(2k-1)^3} + 2\sum^{\infty}_{k=1}\frac{h_k}{(2k-1)^3} \left( \sum^{k}_{i=1} \frac{h_i}{i} \right)  \nonumber \\ &-& 
4\sum^{\infty}_{k=1}\frac{h_k^2}{(2k-1)^3} + 2\sum^{\infty}_{k=1}\frac{H_k h_k}{(2k-1)^3} + \frac{1}{2}\sum^{\infty}_{k=1}\frac{H_k^{(2)} h_k}{(2k-1)^3} +
\frac{1}{2}\sum^{\infty}_{k=1}\frac{H_k^2 h_k}{(2k-1)^3}~.
\end{eqnarray}
The second sum and the last sum have been defined in Eq.~(227) and Eq.~(229). The nested double sum on the right side of Eq.~(313) is explicitly known from Eq.~(151). This way Eq.~(314) results.

\begin{eqnarray}
\sum^{\infty}_{k=1} \frac{H_k^{(3)}  h_k^2}{k(2k-1)} &=& \frac{495}{64}\zeta(6) + \frac{217}{8}ln(2)\zeta(5) + \frac{31}{8}\zeta(5) - \frac{45}{4}\zeta(4)
- \frac{53}{8}\zeta(3)^2 - \frac{25}{2}ln(2)\zeta(2)\zeta(3) - \frac{7}{2}\zeta(2)\zeta(3) + 3\zeta(3) \nonumber \\ &+& 12ln(2)\zeta(2) + \frac{1}{2}\zeta(2)\sum^{\infty}_{k=1}\frac{h_k}{k^3}
- \frac{3}{2}\sum^{\infty}_{k=1}\frac{h_k}{k^5}~.
\end{eqnarray}

In order to compute the last cubic Euler sum of the seventh family we use of the following two-valued help function which again is the solution of the corresponding inhomogeneous difference 
equation:
\begin{eqnarray}
\sum^{\infty}_{i=1} \frac{H_i}{i} \left( \sum^{\infty}_{k=1} \frac{h_k^{(3)}}{(2k-1)(i+k)} \right) &=& \sum^{\infty}_{i=1} \frac{H_i}{i(2i+1)} \left( \sum^{\infty}_{k=1} \frac{h_k^{(3)}}{k(2k-1)} \right)
+ \frac{7}{4}\zeta(3) \sum^{\infty}_{i=1} \frac{H_i h_i}{i(2i+1)} - \frac{3}{2}\zeta(2) \sum^{\infty}_{i=1} \frac{H_i h_i^{(2)}}{i(2i+1)} \nonumber \\ &+& 2ln(2) \sum^{\infty}_{i=1} \frac{H_i h_i^{(3)}}{i(2i+1)}
+  \sum^{\infty}_{i=1} \frac{H_i}{i(2i+1)} \left( \sum^{i}_{k=1} \frac{H_{k-1}}{(2k-1)^3} \right)
\end{eqnarray}
or
\begin{eqnarray}
\sum^{\infty}_{k=1} \frac{h_k^{(3)}}{(2k-1)} \left( \sum^{\infty}_{i=1} \frac{H_i}{i(k+i)} \right) &=& \zeta(2)\sum^{\infty}_{k=1} \frac{h_k^{(3)}}{k(2k-1)} - \sum^{\infty}_{k=1} \frac{H_k h_k^{(3)}}{k^2(2k-1)}
+ \frac{1}{2}\sum^{\infty}_{k=1} \frac{H_k^{(2)} h_k^{(3)}}{k(2k-1)} + \frac{1}{2}\sum^{\infty}_{k=1} \frac{H_k^2 h_k^{(3)}}{k(2k-1)}~.
\end{eqnarray}
The last Euler sum on the right side of Eq.~(316) can be identified as the eight member of the seventh family. This way, only the double sum on the left side of Eq.~(316) is unknown. It follows first: 
\begin{eqnarray}
\sum^{\infty}_{i=1} \frac{H_i}{i(2i+1)} \left( \sum^{i}_{k=1} \frac{H_{k-1}}{(2k-1)^3} \right) &=& \sum^{\infty}_{i=1} \frac{H_i}{i(2i+1)} \left( \sum^{\infty}_{k=1} \frac{H_{k-1}}{(2k-1)^3} \right) + 
\sum^{\infty}_{i=1} \frac{H_i H_{i-1}}{i(2i+1)(2i-1)^3} \nonumber \\ &-& \sum^{\infty}_{k=1} \frac{H_{k-1}}{(2k-1)^3} \left( \sum^{k}_{i=1} \frac{H_i}{i} - 2 \sum^{k}_{i=1} \frac{H_i}{2i+1} \right)~.
\end{eqnarray} 
Here, it remains to compute the last two double sums on the right side of Eq.~(317). First we have:
\begin{eqnarray}
\sum^{\infty}_{k=1} \frac{H_{k-1}}{(2k-1)^3} \left( \sum^{k}_{i=1} \frac{H_i}{i} \right) &=& \frac{1}{2} \sum^{\infty}_{k=1} \frac{H_{k-1} H_k^{(2)}}{(2k-1)^3} + \frac{1}{2}
\sum^{\infty}_{k=1} \frac{H_{k-1} H_k^2}{(2k-1)^3}~.
\end{eqnarray}
Both Euler sums on the right side of Eq.~(318) are known explicitly. 
This way we continue with:
\begin{eqnarray}
\sum^{\infty}_{k=1} \frac{H_{k-1}}{(2k-1)^3} \left( \sum^{k}_{i=1} \frac{H_i}{2i+1} \right) &=& \sum^{\infty}_{k=1} \frac{H_{k-1}}{(2k-1)^3} \left( \sum^{k+1}_{i=1} \frac{H_{i-1}}{2i-1} \right) \nonumber \\ &=& 
\sum^{\infty}_{k=1} \frac{H_{k-1}}{(2k-1)^3} \left( \sum^{k}_{i=1} \frac{H_{i-1}}{2i-1} \right) + \sum^{\infty}_{k=1} \frac{H_{k-1} H_k}{(2k-1)^3(2k+1)}  \nonumber \\ &=&
\sum^{\infty}_{k=1} \frac{H_{k-1} H_k h_k}{(2k-1)^3} - \sum^{\infty}_{k=1} \frac{H_{k-1}}{(2k-1)^3} \left( \sum^{k}_{i=1} \frac{h_i}{i}\right) + \sum^{\infty}_{k=1} \frac{H_{k-1} H_k}{(2k-1)^3(2k+1)}~.
\end{eqnarray}
As the nested double sum on the right side of Eq.~(319) is known explicitly from Eq.~(307) it follows for the eight member of the seventh family:

\begin{eqnarray}
\sum^{\infty}_{k=1} \frac{H_k^2 h_k^{(3)}}{k(2k-1)} &=& \frac{135}{128}\zeta(6) - \frac{527}{16}ln(2)\zeta(5) + \frac{403}{16}\zeta(5) + \frac{165}{8}\left(ln(2)\right)^2\zeta(4) 
- 30ln(2)\zeta(4) + \frac{75}{8}\zeta(4) + \frac{119}{32}\zeta(3)^2 \nonumber \\ &+& 7ln(2)\zeta(2)\zeta(3) - \frac{21}{4}\zeta(2)\zeta(3) - 7\left(ln(2)\right)^3\zeta(3) + 
14\left(ln(2)\right)^2\zeta(3) - 7ln(2)\zeta(3) \nonumber \\ &-& \frac{7}{4}\zeta(2)\sum^{\infty}_{k=1}\frac{h_k}{k^3} + \left(ln(2)\right)^2\sum^{\infty}_{k=1}\frac{h_k}{k^3} -
2ln(2)\sum^{\infty}_{k=1}\frac{h_k}{k^3} + \sum^{\infty}_{k=1}\frac{h_k}{k^3} + \frac{3}{4}\sum^{\infty}_{k=1}\frac{h_k}{k^5}~.
\end{eqnarray}

\subsection{Eighth family}
The eighth family consists of two members. The first one results to
\begin{eqnarray}
\sum^{\infty}_{k=1} \frac{h_k^3}{k^3} = \frac{135}{128}\zeta(6) + \frac{7}{16}\zeta(3)^2~.
\end{eqnarray}

For the second one we found:
\begin{eqnarray}
\sum^{\infty}_{k=1} \frac{h_k^{(2)} h_k^2}{k^2} = \frac{135}{64}\zeta(6) + \frac{7}{8}\zeta(3)^2~.
\end{eqnarray}
Both members are known explicitly as discussed before. At last it should be mentioned here that 
\begin{eqnarray}
\sum^{\infty}_{k=1}\frac{h_k}{k^5} &=& -\frac{10669}{312}\zeta(6) + 31ln(2)\zeta(5) - \frac{120}{13}\left(ln(2)\right)^2\zeta(4) - \frac{7}{4}\zeta(3)^2 + \frac{6}{13}\left(ln(2)\right)^4\zeta(2)
- \frac{1}{13}\left(ln(2)\right)^6 \nonumber \\ &-& \frac{2016}{13}Li_6(1/2) - \frac{2592}{13}Li_6(-1/2) + \frac{16}{39}Li_6(-1/8)
\end{eqnarray}
holds \cite{Au24}, where this expression results from motivic techniques. This way the corresponding odd-typ linear Euler sums of order 6 can be expressed explicitly in terms of $Li_6$, although $Li_6(1/2)$,
$Li_6(-1/2)$ and $Li_6(-1/8)$ are needed for an explicit calculation.

\section{Summary}
We have introduced a generalized summation method that allows in combination with MZV theory \cite{au20} to calculate explicitly all eighth families of cubic Euler sums of orders 4 to 6 in terms of zeta
values and corresponding polylogarithmic values $Li_4(1/2)$,$Li_5(1/2)$, $Li_6(1/2)$, $Li_6(-1/2)$ and $Li_6(-1/8)$.

\section{Appendix A}
Here we present as an outlook few related biquadratic and quintic Euler sums of order six. In a forthcoming paper we will show that in the order six case also all biquadratic and quintic Euler sums can
be explicitly calculated as well in terms of zeta values and polylogarithmic values $Li_5(1/2)$, $Li_6(1/2)$, $Li_6(-1/2)$ and $Li_6(-1/8)$.

As an example, we present for the order 6 case the following biquadratic and quintic Euler sums:
\begin{eqnarray}
\sum^{\infty}_{k=1} \frac{h_k^4}{(2k-1)^2} &=& \frac{183}{16}\zeta(6) - \frac{713}{64}ln(2)\zeta(5) + \frac{249}{32}\left(ln(2)\right)^2\zeta(4) + \frac{161}{256}\zeta(3)^2 + \frac{3}{8}ln(2)\zeta(2)\zeta(3)
+ \left(ln(2)\right)^4\zeta(2) \nonumber \\ &-& \frac{1}{60}\left(ln(2)\right)^6 + \frac{23}{64} \sum^{\infty}_{k=1} \frac{h_k}{k^5} - 12Li_6\left(\frac{1}{2}\right)
\end{eqnarray}
and
\begin{eqnarray}
\sum^{\infty}_{k=1} \frac{H_k^3 h_k}{k^2} &=& -\frac{1851}{32}\zeta(6) + \frac{159}{4}\left(ln(2)\right)^2\zeta(4) - \frac{35}{32}\zeta(3)^2 - 28\left(ln(2)\right)^3\zeta(3) + 6\left(ln(2)\right)^4\zeta(2)
- \frac{2}{3} \left(ln(2)\right)^6 \nonumber \\ &-& \frac{3}{2}\zeta(2)\sum^{\infty}_{k=1}\frac{h_k}{k^3} + 6\left(ln(2)\right)^2 \sum^{\infty}_{k=1}\frac{h_k}{k^3} - \frac{23}{8} \sum^{\infty}_{k=1} \frac{h_k}{k^5}
+ 96Li_6\left(\frac{1}{2}\right) + 96ln(2)Li_5\left(\frac{1}{2}\right)
\end{eqnarray}

\begin{eqnarray}
\sum^{\infty}_{k=1} \frac{h_k^5}{k(2k-1)} &=& \frac{10829}{512}\zeta(6) - \frac{1209}{64}ln(2)\zeta(5) + \frac{565}{32}\left(ln(2)\right)^2\zeta(4) + \frac{245}{256}\zeta(3)^2 + \frac{45}{16}ln(2)\zeta(2)\zeta(3) 
+ \frac{25}{6}\left(ln(2)\right)^4\zeta(2) \nonumber \\ &-& \frac{1}{36}\left(ln(2)\right)^6 + \frac{39}{64}\sum^{\infty}_{k=1}\frac{h_k}{k^5} - 20Li_6\left(\frac{1}{2}\right)
\end{eqnarray}
and
\begin{eqnarray}
\sum^{\infty}_{k=1} \frac{H_k^3 h_k^2}{k(2k-1)} &=& \frac{16071}{32}\zeta(6) - \frac{2477}{8}ln(2)\zeta(5) + \frac{3683}{16}\zeta(5) + \frac{867}{8}\left(ln(2)\right)^2\zeta(4) - 171ln(2)\zeta(4)
+ \frac{57}{2}\zeta(4) + \frac{1115}{64}\zeta(3)^2 \nonumber \\ &+& \frac{25}{4}ln(2)\zeta(2)\zeta(3) - \frac{11}{4}\zeta(2)\zeta(3) - 3\left(ln(2)\right)^3\zeta(3) + 9\left(ln(2)\right)^2\zeta(3)
+ 33ln(2)\zeta(3) + 3\zeta(3) \nonumber \\ &-& 6\left(ln(2)\right)^4\zeta(2) + 20\left(ln(2)\right)^3\zeta(2) - 36\left(ln(2)\right)^2\zeta(2) + 12ln(2)\zeta(2) - \frac{2}{5}\left(ln(2)\right)^6
+ \frac{8}{5}\left(ln(2)\right)^5 \nonumber \\ &-& \frac{1}{4}\zeta(2)\sum^{\infty}_{k=1}\frac{h_k}{k^3} - 6\sum^{\infty}_{k=1}\frac{h_k}{k^3} +  \frac{159}{16}\sum^{\infty}_{k=1} \frac{h_k}{k^5}
- 288Li_6\left(\frac{1}{2}\right) - 192Li_5\left(\frac{1}{2}\right)~.
\end{eqnarray}

In the order 7 case the same calculational scheme can be applied for an explicit calculation of all corresponding Euler sums of higher degrees, although two additional multiple zeta values
$\zeta({\overline 5},1,1)$ and $\zeta({\overline 3},3,1)$ must be used in the computational scheme in addition to the corresponding zeta and polylogarithmic values.
As a first example we present the following cubic Euler sums:

\begin{eqnarray}
\sum^{\infty}_{k=1} \frac{h_k^3}{k^4} =  \frac{93}{4}\zeta(2)\zeta(5) - \frac{945}{32}\zeta(3)\zeta(4)
\end{eqnarray}

and
\begin{eqnarray}
\sum^{\infty}_{k=1} \frac{h_k^{(2)} h_k^2}{k^3} =  -\frac{93}{8}\zeta(2)\zeta(5) + \frac{525}{32}\zeta(3)\zeta(4)
\end{eqnarray}

and
\begin{eqnarray}
\sum^{\infty}_{k=1} \frac{h_k h_k^{(2)} h_k^{(2)}}{k^2} &=& \frac{1113}{64}ln(2)\zeta(6) - \frac{1023}{128}\zeta(2)\zeta(5) - \frac{105}{32}\zeta(3)\zeta(4) + 
\frac{5}{2}\left(ln(2)\right)^3\zeta(4) - \frac{21}{4}\left(ln(2)\right)^2\zeta(2)\zeta(3) \nonumber \\ &-& \frac{1}{10}\left(ln(2)\right)^5\zeta(2) +
\frac{3}{2}ln(2)\zeta(2)\sum^{\infty}_{k=1}\frac{h_k}{k^3} + 12\zeta(2) Li_5\left(\frac{1}{2}\right)
\end{eqnarray}

and
\begin{eqnarray}
\sum^{\infty}_{k=1} \frac{h_k^{(3)} h_k^2}{k^2} &=& -\frac{1113}{64}ln(2)\zeta(6) + \frac{1767}{128}\zeta(2)\zeta(5) - \frac{5}{2}\left(ln(2)\right)^3\zeta(4) +
\frac{21}{4}\left(ln(2)\right)^2\zeta(2)\zeta(3) + \frac{1}{10}\left(ln(2)\right)^5\zeta(2) \nonumber \\ &-& \frac{3}{2}ln(2)\zeta(2)\sum^{\infty}_{k=1}\frac{h_k}{k^3} - 
12\zeta(2) Li_5\left(\frac{1}{2}\right)
\end{eqnarray}

and
\begin{eqnarray}
\sum^{\infty}_{k=1} \frac{H_k^{(2)} h_k^2}{k^3} &=& \frac{77959}{544}\zeta(7) - \frac{9553}{272}\zeta(2)\zeta(5) - \frac{18481}{272}\zeta(3)\zeta(4) + \frac{7}{2}\zeta(3)
\sum^{\infty}_{k=1}\frac{h_k}{k^3} \nonumber \\ &+& \frac{1008}{17}\zeta({\overline 5},1,1) - \frac{168}{17}\zeta({\overline 3},3,1)~.
\end{eqnarray}

\end{document}